\pdfoutput=1
\documentclass[12pt]{elsarticle}
\usepackage{amsmath}
\usepackage{amssymb}
\usepackage{amsthm}
\usepackage{enumerate}
\usepackage[section]{algorithm}
\usepackage{algpseudocode}
\usepackage{varwidth}
\usepackage{setspace}
\usepackage{graphicx}
\usepackage{epstopdf}
\usepackage{subfigure}
\usepackage{tabularx}
\usepackage{relsize}
\usepackage[top=1in, bottom=1in, left=1in, right=1in]{geometry}
\usepackage{mathtools}

\graphicspath{{./Figures/}}

\usepackage[textsize=footnotesize]{todonotes}

\DeclareMathAlphabet{\mathpzc}{OT1}{pzc}{m}{it}

\newtheorem*{proposition}{Proposition}

\makeatletter
\renewcommand{\p@enumii}{\theenumi.}
\makeatother

\journal{Computer Methods in Applied Mechanics and Engineering}

\begin{document}

\begin{frontmatter}

\title{High-Order Finite Element Methods for Moving Boundary Problems with Prescribed Boundary Evolution}

\author[icme]{Evan S. Gawlik}
\ead{egawlik@stanford.edu}
\author[me,icme]{Adrian J. Lew}
\ead{lewa@stanford.edu}

\address[icme]{Computational and Mathematical Engineering, Stanford University}
\address[me]{Mechanical Engineering, Stanford University}

\begin{abstract}
  We introduce a framework for the design of finite element methods
  for two-dimensional moving boundary problems with prescribed boundary evolution that
  have arbitrarily high order of accuracy, both in space and in time.
  At the core of our approach is the use of a universal mesh: a
  stationary background mesh containing the domain of interest for all
  times that adapts to the geometry of the immersed domain by
  adjusting a small number of mesh elements in the neighborhood of the
  moving boundary.  The resulting method maintains an exact
  representation of the (prescribed) moving boundary at the discrete
  level, or an approximation of the appropriate order, yet is immune
  to large distortions of the mesh under large deformations of the
  domain.  The framework is general, making it possible to achieve any
  desired order of accuracy in space and time by selecting a preferred
  and suitable finite-element space on the universal mesh for the
  problem at hand, and a preferred and suitable time integrator for
  ordinary differential equations.  We illustrate our approach by
  constructing a particular class of methods, and apply them to a
  prescribed-boundary variant of the Stefan problem.  We present
  numerical evidence for the order of accuracy of our schemes in one
  and two dimensions.  
\end{abstract}


\end{frontmatter}

\section{Introduction}

\begin{figure}[t]
  \centering
  \subfigure[{$t=0$}]{
    \centering
    \includegraphics[trim = 0in 0.2in 0in 0.2in, clip=true, width=0.45\textwidth]{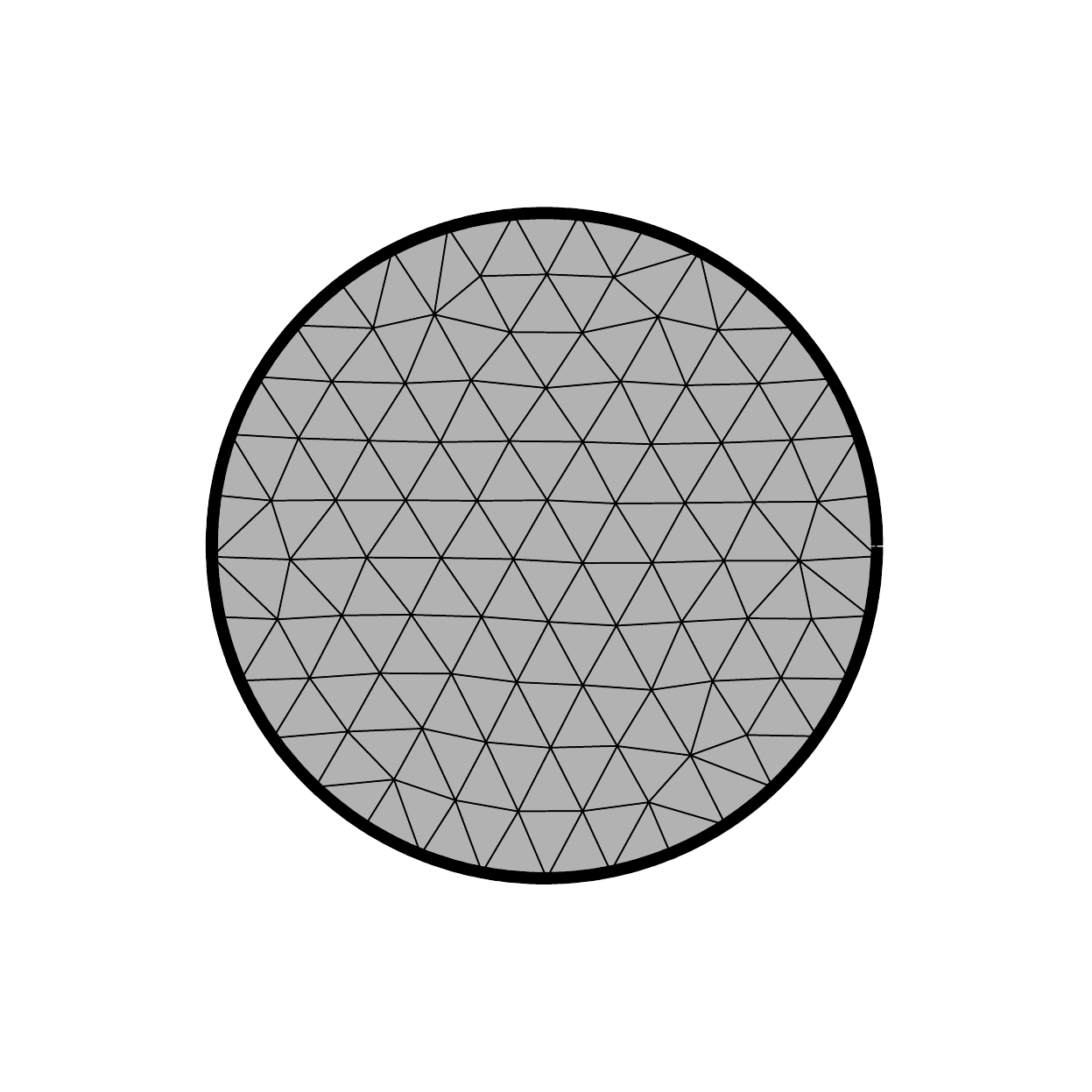}
    \label{fig:deformcirc}
  }
  \subfigure[{$t>0$}]{
    \centering
    \includegraphics[trim = 0in 0.2in 0in 0.2in, clip=true, width=0.45\textwidth]{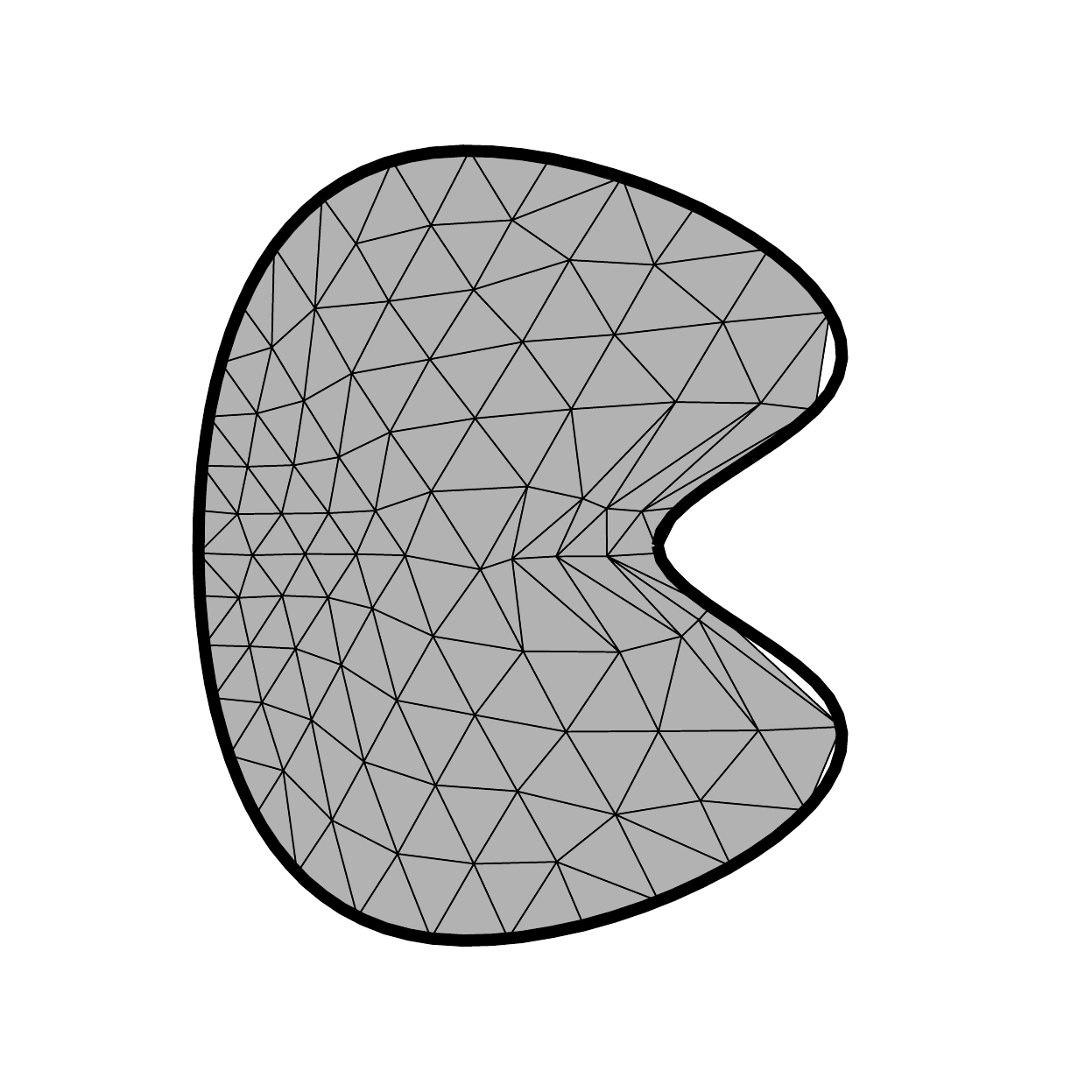}
    \label{fig:deformspline}
  }
  \caption{Schematic depiction of a deforming-mesh method.  Without a careful choice of nodal motions, elements can suffer unwanted distortions under large deformations of the moving domain.\protect\footnotemark}
  \label{fig:deform}
\end{figure}

\begin{figure}[t]
  \centering
  \subfigure[{$t=0$}]{
    \centering
    \includegraphics[trim = 0in 0.2in 0in 0.2in, clip=true, width=0.45\textwidth]{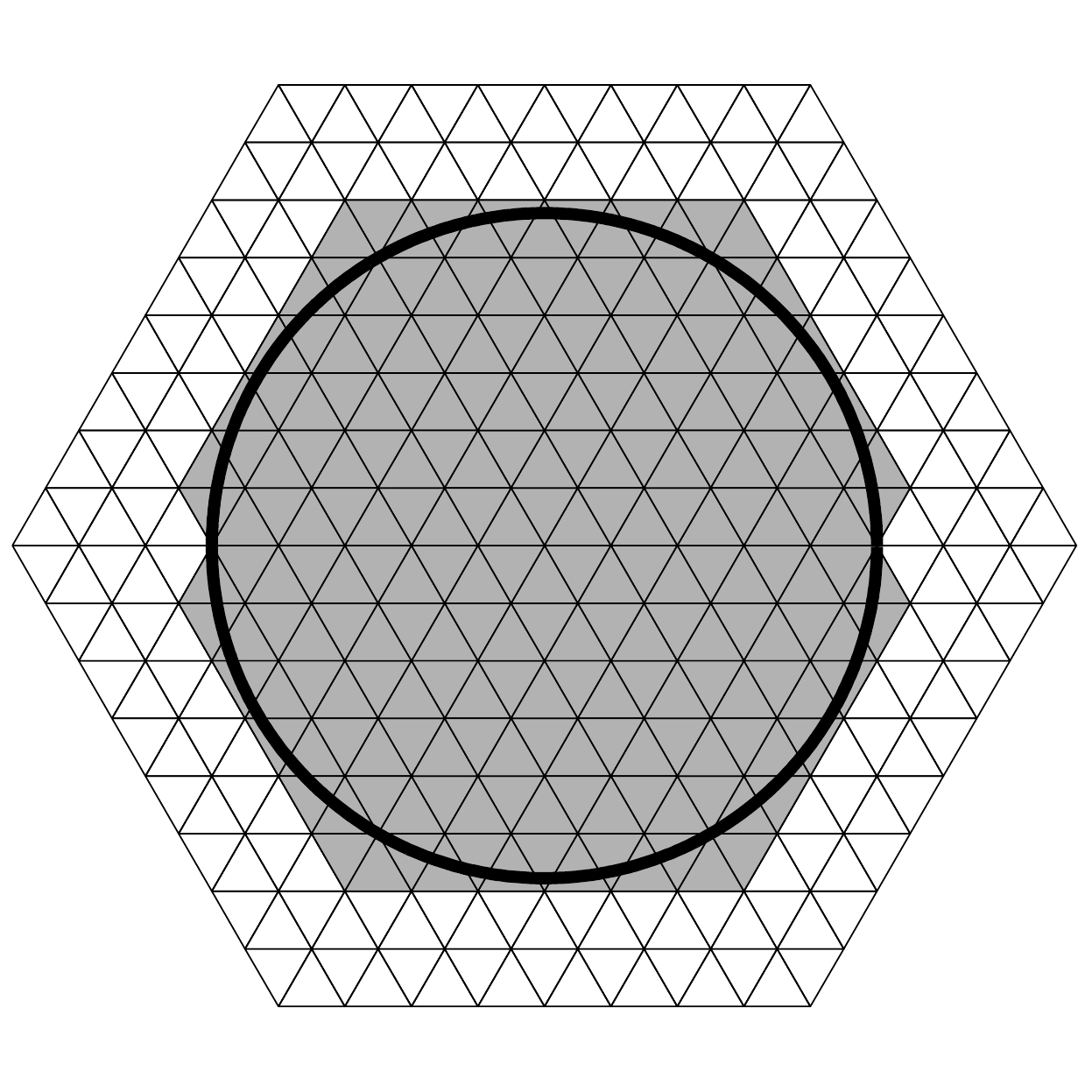}
    \label{fig:fixedcirc}
  }
  \subfigure[{$t>0$}]{
    \centering
    \includegraphics[trim = 0in 0.2in 0in 0.2in, clip=true, width=0.45\textwidth]{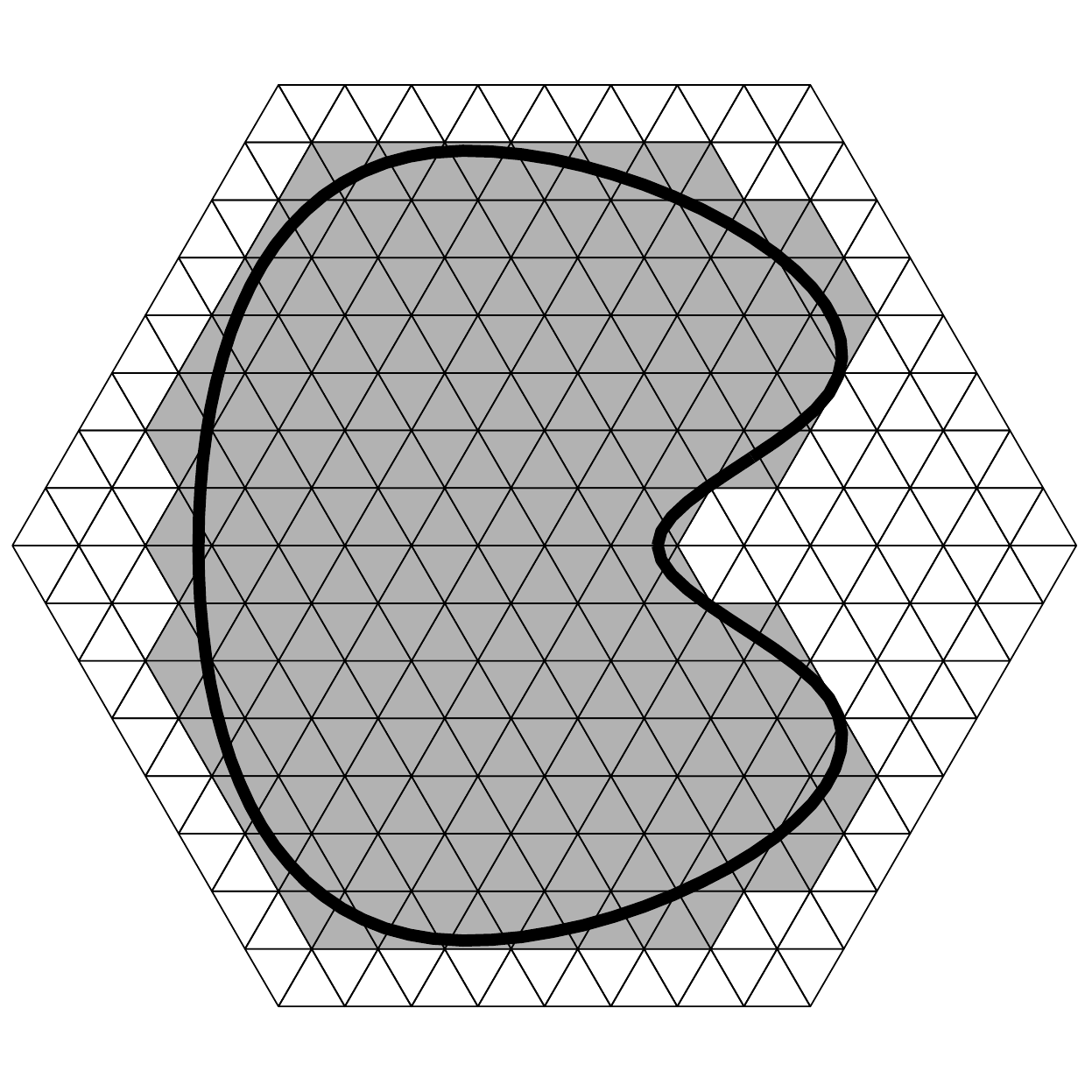}
    \label{fig:fixedspline}
  }      
  \caption{Schematic depiction of a fixed-mesh method.  Such methods employ a fixed background mesh which does not conform to the immersed domain.}
  \label{fig:fixed}
\end{figure}

\begin{figure}[t]
  \centering
  \subfigure[{$t=0$}]{
    \centering
    \includegraphics[trim = 0in 0.2in 0in 0.2in, clip=true, width=0.45\textwidth]{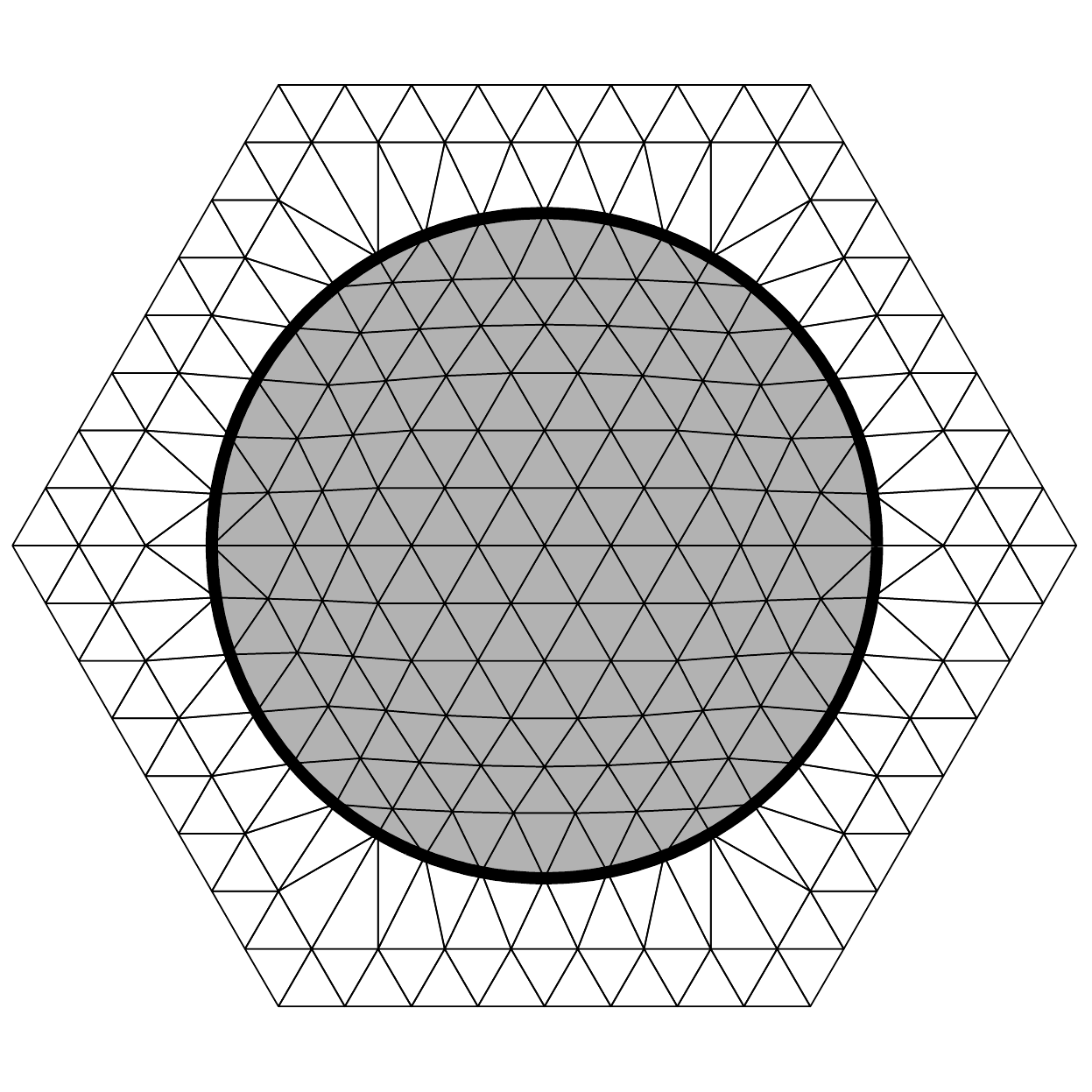}
    \label{fig:unicirc}
  }
  \subfigure[{$t>0$}]{
    \centering
    \includegraphics[trim = 0in 0.2in 0in 0.2in, clip=true, width=0.45\textwidth]{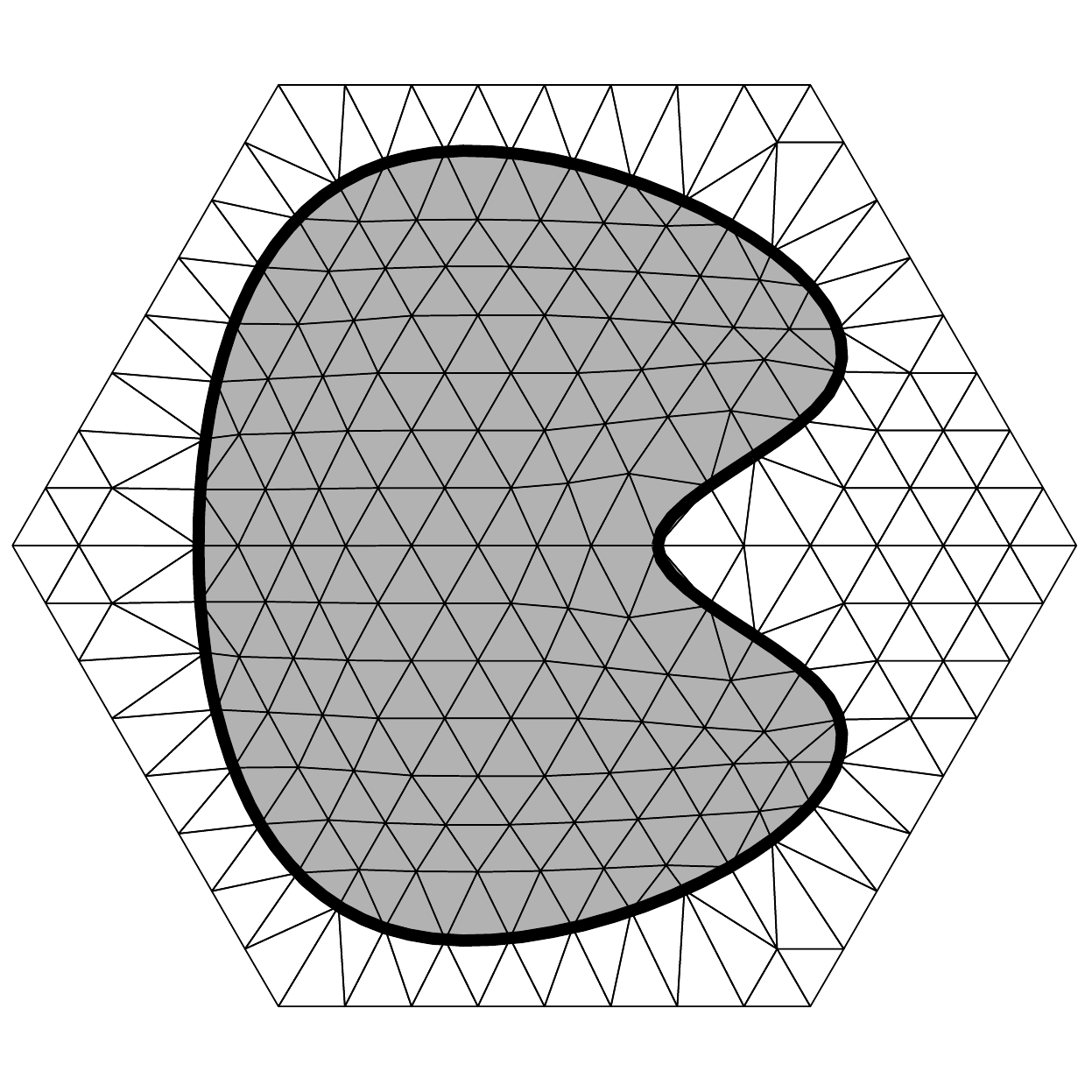}
    \label{fig:unispline}
  }
  \caption{Schematic depiction of a universal mesh.  By adapting the mesh to the immersed domain, one obtains a mesh that conforms to the domain exactly and is immune to large distortions of elements.}
  \label{fig:uni}
\end{figure}

Science and engineering are replete with instances of moving boundary problems: partial differential equations posed on domains that change with time.  Problems of this type, which arise in areas as diverse as fluid-structure interaction, multiphase flow physics, and fracture mechanics, are inherently challenging to solve numerically.

Broadly speaking, computational methods for moving boundary problems generally adhere to one of two paradigms.  \emph{Deforming-mesh} methods employ a computational mesh that deforms in concert with the moving domain, whereas \emph{fixed-mesh} methods employ a stationary background mesh in which the domain is immersed.  While the former approach can require that efforts be made to avoid distortions of the mesh under large deformations~\cite{Saksono2007}, the latter approach requires that special care be taken in order to account for any discrepancy between the exact boundary and element interfaces~\cite{Gibou2005,Nochetto1997}.  Figs.~\ref{fig:deform}-\ref{fig:fixed} illustrate these two paradigms schematically.

In this study, we eliminate these difficulties by employing a
\emph{universal mesh}: a stationary background mesh that adapts to the
geometry of the immersed domain by adjusting a small number of mesh
elements in the neighborhood of the moving boundary.  An example is
illustrated in Fig.~\ref{fig:uni}.  The resulting framework admits, in
a general fashion, the construction of methods that are of arbitrarily
high order of accuracy in space and time, without exhibiting the
aforementioned drawbacks of deforming-mesh and fixed-mesh
methods. This strategy was introduced for time-independent and
quasi-steady problems in~\cite{Rangarajan2011,Rangarajan2012a}. Here
we present its extension to time-dependent problems posed on moving
domains with prescribed evolution.  We relegate a discussion of
problems with unprescribed boundaries to future work, since the
treatment of unprescribed boundaries introduces its own set of
challenges -- approximation of the boundary, discretization of the
boundary evolution equations, and error analysis on approximate
domains -- that may have the undesired effect of blurring the focus of the present study.

In the process of deriving our method, we present a unified, geometric framework that puts our method and existing deforming-mesh methods on a common footing suitable for analysis.  The main idea is to recast the governing equations on a sequence of cylindrical spacetime slabs that span short intervals of time. The clarity brought about by this geometric viewpoint renders the analysis of numerical methods for moving-boundary problems more tractable, as it reduces the task to a standard analysis of fixed-domain problems with time-dependent PDE coefficients.

\footnotetext{Here, for purely illustrative purposes, we have employed a nodal mapping of the form $(r,\theta) \mapsto (f(\theta)r,\theta)$ in polar coordinates.}

\paragraph{Organization} This paper is organized as follows.  We begin
in \S\ref{sec:overview} by giving an informal overview of our
method, and illustrating the ideas by formulating the method for a
moving boundary problem in one spatial dimension.  We formulate a
two-dimensional model moving
boundary problem on a predefined, curved spacetime domain in \S\ref{sec:model_problem}, and proceed to derive its equivalent
reformulation on cylindrical spacetime slabs.  In
\S\ref{sec:discretization} we present, in an abstract manner,
the general form of a finite-element discretization of the same moving
boundary problem, as well as its reformulation on cylindrical
spacetime slabs.  This formalism will lead to a statement of the
general form of a numerical method for moving boundary problems with
prescribed boundary evolution that includes our method and
conventional deforming-mesh methods as special cases.  We finish
\S\ref{sec:discretization} by summarizing an error estimate for
methods of this form, referring the reader to our companion
paper~\cite{Gawlik2012b} for its proof.  In
\S\ref{sec:universal_meshes}, we present the key ingredient that
distinguishes our proposed method from standard approaches: the use of
a universal mesh.  We specialize the aforementioned error estimate to
this setting to deduce that the method's convergence rate is
suboptimal by half an order when the time step and mesh spacing scale
proportionately.  In \S\ref{sec:numerical_examples} we
demonstrate numerically our method's convergence rate on a
prescribed-boundary variant of a classic moving-boundary problem
called the Stefan problem, which asks for the evolution of a
solid-liquid interface during a melting process.  
Some concluding remarks are given in \S\ref{sec:conclusion}.

\paragraph{Previous work}

In what follows, we review some of the existing numerical methods for
moving-boundary problems, beginning with deforming-mesh methods and
finishing with fixed-mesh methods.

Deforming-mesh methods have enjoyed widespread success in the
scientific and engineering communities, where they are best known as
Arbitrary Lagrangian Eulerian (ALE) methods.  The appellation refers
to the fact that in prescribing a motion of the mesh, a kinematic
description of the physics is introduced that is neither Eulerian (in
which the domain moves over a fixed mesh) nor Lagrangian (in which the
domain does not move with respect to the mesh).  The resulting
formalism leads to governing equations that contain a term involving
the velocity of the prescribed mesh motion that is otherwise absent in
schemes on a fixed mesh~\cite{Lynch1982,Lesaint1989}.  Early
appearances of the ALE framework date back to the works of Hirt
et.~al.~\cite{Hirt1974}, Hughes et.~al.~\cite{Hughes1981}, and Donea
et.~al.~\cite{Donea1982}.  ALE methods have seen use in
fluid-structure
interaction~\cite{Souli2000,Geuzaine2003,Farhat2004,Farhat2006,Farhat2010,Takashi1994},
solid mechanics~\cite{Liu1988,Wang1997,Askes2004,Khoei2008},
thermodynamics~\cite{Sullivan1987,Zabaras1990,Albert1986,Beckett2001},
and other applications.

Relative to methods for problems with fixed domains, less attention
has been directed toward the development of ALE methods of high order
of accuracy and the associated error analysis.  Schemes of
second-order in time are
well-studied~\cite{Farhat2006,Farhat2010,Geuzaine2003,Boffi2004,Mackenzie2012,Formaggia2004,Formaggia1999,Gastaldi2001},
though the analysis of higher-order schemes has only recently been
addressed by Bonito and co-authors~\cite{Bonito2012a,Bonito2012b}, who
study the spatially continuous setting with discontinuous Galerkin
temporal discretizations.

One of the key challenges that ALE methods face is the maintenance of
a good-quality mesh during large deformations of the
domain~\cite{Aymone2004,Bar2001}.  Fig.~\ref{fig:deform} illustrates a
case where, using an intentionally naive choice of nodal motions, a
domain deformation can lead to triangles with poor aspect ratios.  In
more severe cases, element inversions can occur.  Such distortions are
detrimental both to the accuracy of the spatial discretization and to
the conditioning of the discrete governing equations~\cite{Masud2006}.
For this reason, it is common to use sophisticated mesh motion strategies that involve solving systems of equations (such as those of linear elasticity) for the positions of mesh nodes~\cite{Donea1983,Farhat1998,Johnson1994,Helenbrook2003}.

A related class of methods are spacetime methods
(e.g.,~\cite{Tezduyar2007}), where the spacetime domain swept out by
the moving spatial domain is discretized with straight or curved
elements.  These methods resemble deforming-mesh methods in the sense
that spatial slices of the spacetime mesh at fixed temporal nodes
constitute a mesh of the moving domain at those times.  Bonnerot and
Jamet~\cite{Bonnerot1974,Bonnerot1979} have used a spacetime framework
to construct high-order methods for the Stefan problem in one
dimension.  They require the use of curved elements along the moving
boundary to achieve the desired temporal accuracy.
Jamet~\cite{Jamet1978} provides a generalization of these high-order
methods to dimensions greater than one in the case that the boundary
evolution is prescribed in advance. More recently, Rhebergen and
Cockburn \cite{rhebergen2013space,rhebergen2012space} created 
 hybridizable-discontinuous-Galerkin-based spacetime methods for
advection-diffusion and incompressible flow problems with moving
domains.

At the other extreme are fixed-mesh methods, which cover a
sufficiently large domain with a mesh and evolve a numerical
representation of the boundary, holding the background mesh fixed~\cite{Peskin2002,Leveque1994,Fadlun2000,Udaykumar2001}.  A
variety of techniques can be used to represent the boundary, including
level sets~\cite{Sethian1999,Gibou2005}, marker
particles~\cite{Zhao2001}, and splines~\cite{Xu2006}.  Fixed mesh
methods require that special care be taken in constructing the
numerical partial differential operators in the neighborhood of the
moving boundary, so as to avoid losses in accuracy arising from the
disagreement between the moving boundary and element interfaces.  Some
authors~\cite{Palle1996,Nochetto1991} propose adaptively refining the
mesh in the neighborhood of the moving boundary to mitigate these
losses. 
In the special
case of a cartesian mesh, Gibou and Fedkiw~\cite{Gibou2005} have
developed a third-order method for the Stefan problem in two
dimensions using extrapolation to allow finite-difference stencils to
extend beyond the moving boundary.

The method presented in this paper classifies neither as a
deforming-mesh method nor as a fixed-mesh method, though it shares
attractive features from both categories.  It exhibits the immunity to
large mesh distortions enjoyed by fixed-mesh methods without
sacrificing the geometric conformity offered by deforming-mesh
methods.  Despite its conceptual simplicity, the method has not been
proposed in the literature.  An idea similar to ours, dubbed a
``fixed-mesh ALE'' method, has recently been proposed by Baiges and
Codina~\cite{Baiges2010,Baiges2011}, though there are several
important differences.  In particular, their method uses element
splitting to define intermediate meshes during temporal integration,
whereas our method leaves the connectivity of the mesh intact.
Second, they advocate imposing boundary conditions approximately to
improve efficiency; our method imposes boundary conditions exactly
without extra computational effort.  Finally, they focus only on
low-order schemes with piecewise linear approximations to the domain
deformation, while we derive schemes of arbitrarily high order.

\section{Overview of the method} \label{sec:overview}

\begin{figure}[t]
  \centering
  \subfigure[Submesh $\mathcal{S}_h^n$ at $t=t^{n-1}$]{
    \centering
    \includegraphics[trim = 0in 0.2in 0in 0.2in, clip=true, width=0.3\textwidth]{fixedspline_gray-eps-converted-to.pdf}
    \label{fig:step1}
  }
  \subfigure[Conforming mesh for $\Omega^{t^{n-1}}$]{
    \centering
    \includegraphics[trim = 0in 0.2in 0in 0.2in, clip=true, width=0.3\textwidth]{unispline_gray-eps-converted-to.pdf}
    \label{fig:step2}
  }
  \subfigure[Advancement to $t=t^n$]{
    \centering
    \includegraphics[trim = 0in 0.2in 0in 0.2in, clip=true, width=0.3\textwidth]{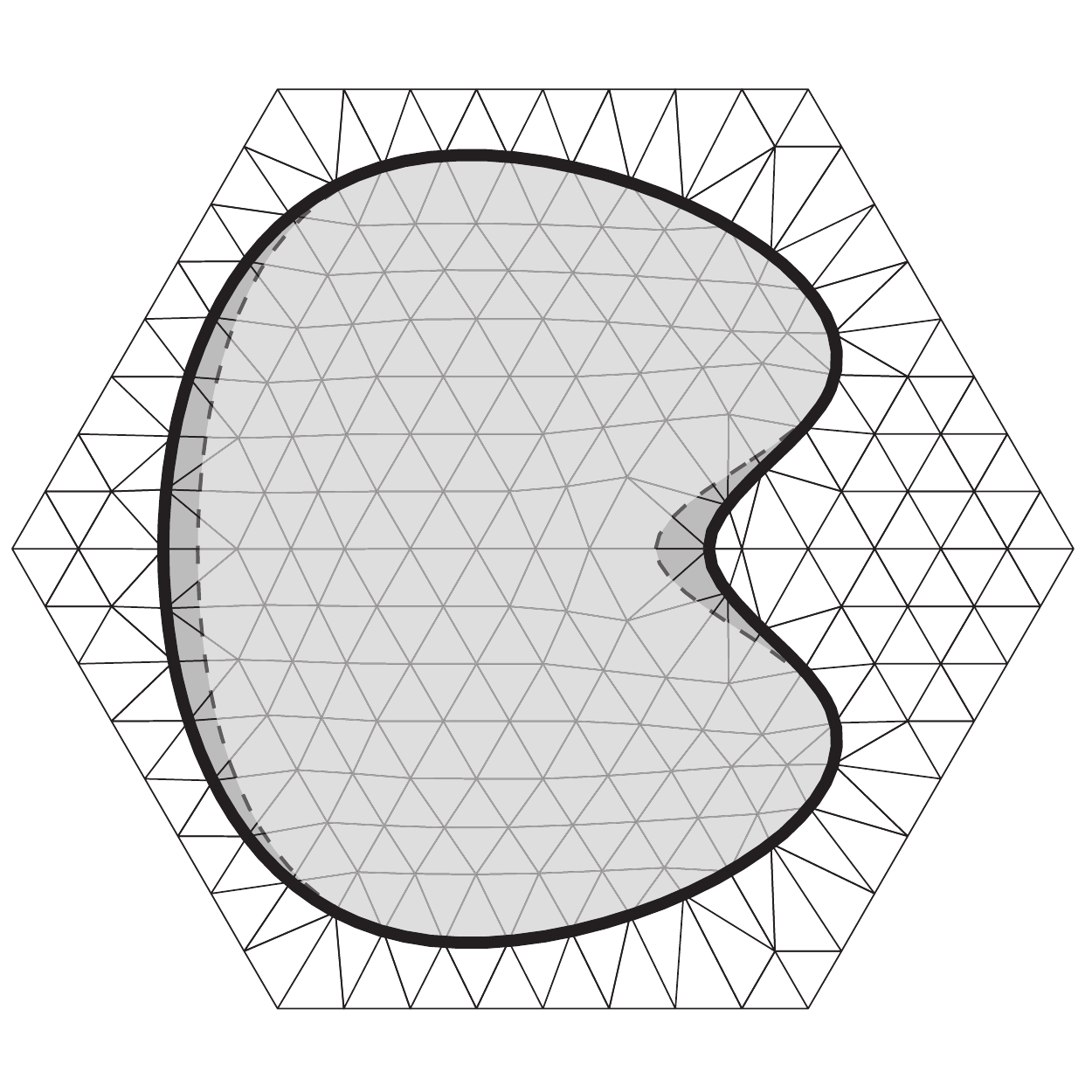}
    \label{fig:step3}
  }
  \subfigure[Submesh $\mathcal{S}_h^{n+1}$ at $t=t^n$]{
    \centering
    \includegraphics[trim = 0in 0.2in 0in 0.2in, clip=true, width=0.3\textwidth]{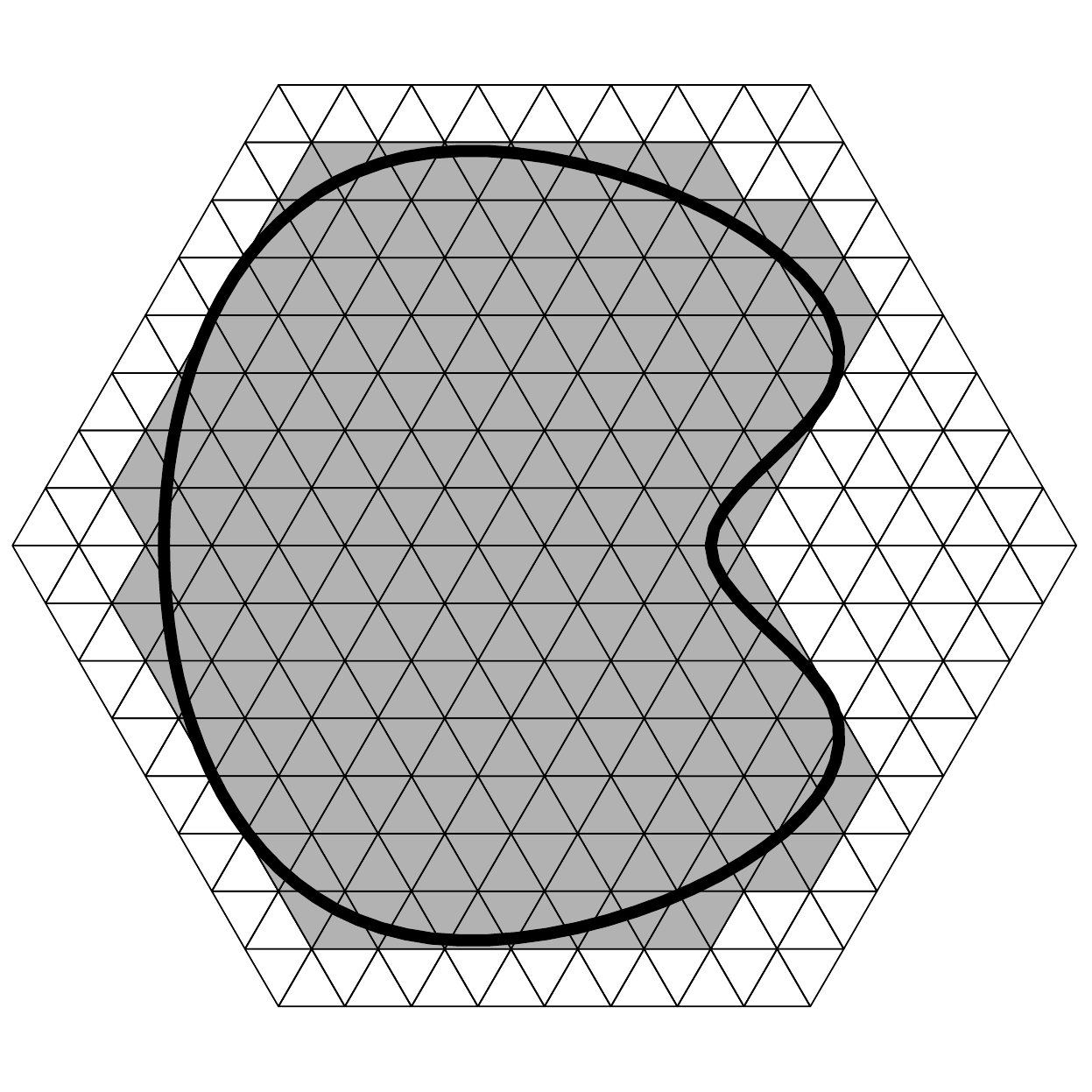}
    \label{fig:step4}
  }
  \subfigure[Conforming mesh for $\Omega^{t^n}$]{
    \centering
    \includegraphics[trim = 0in 0.2in 0in 0.2in, clip=true, width=0.3\textwidth]{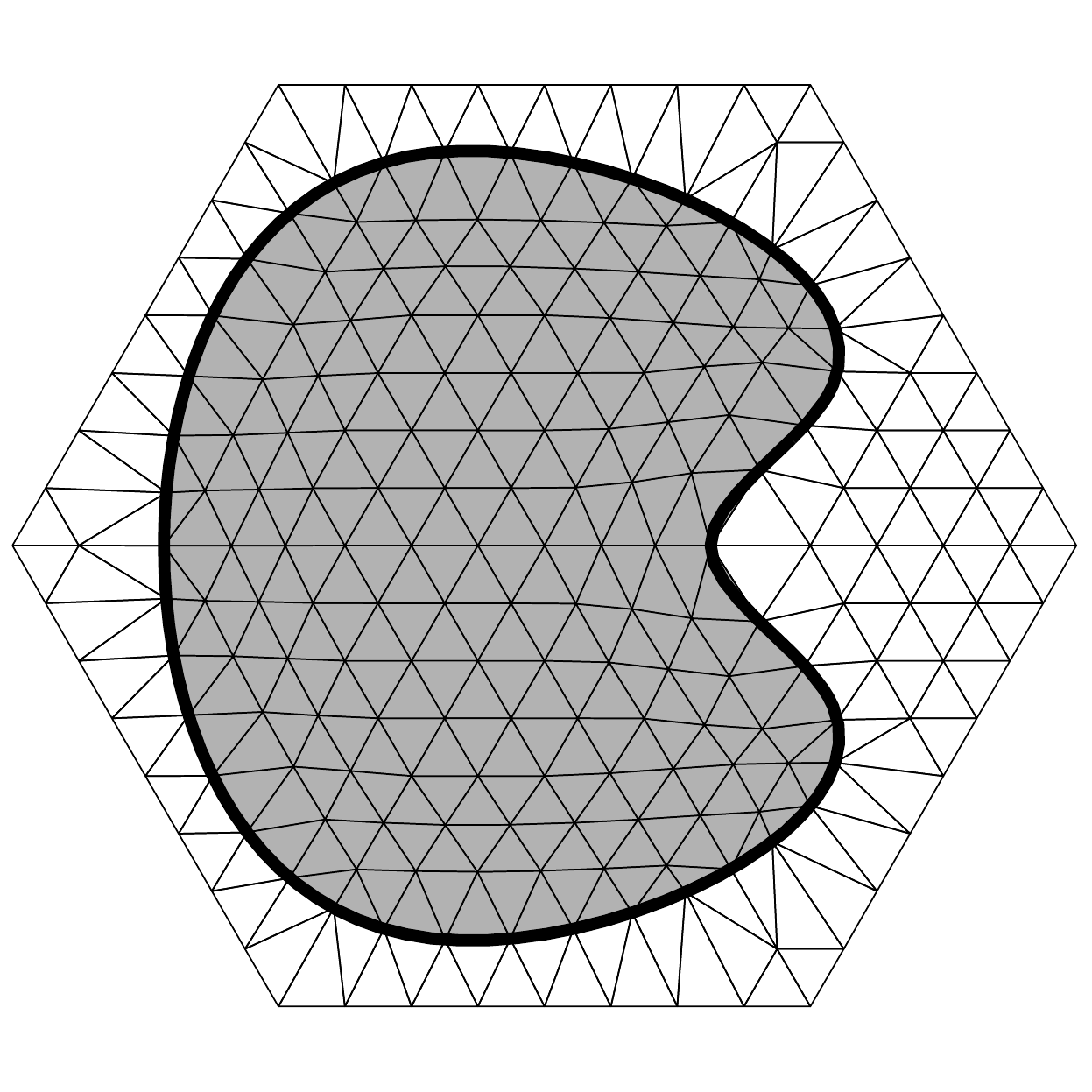}
    \label{fig:step5}
  }
  \subfigure[Advancement to $t=t^{n+1}$]{
    \centering
    \includegraphics[trim = 0in 0.2in 0in 0.2in, clip=true,    width=0.3\textwidth]{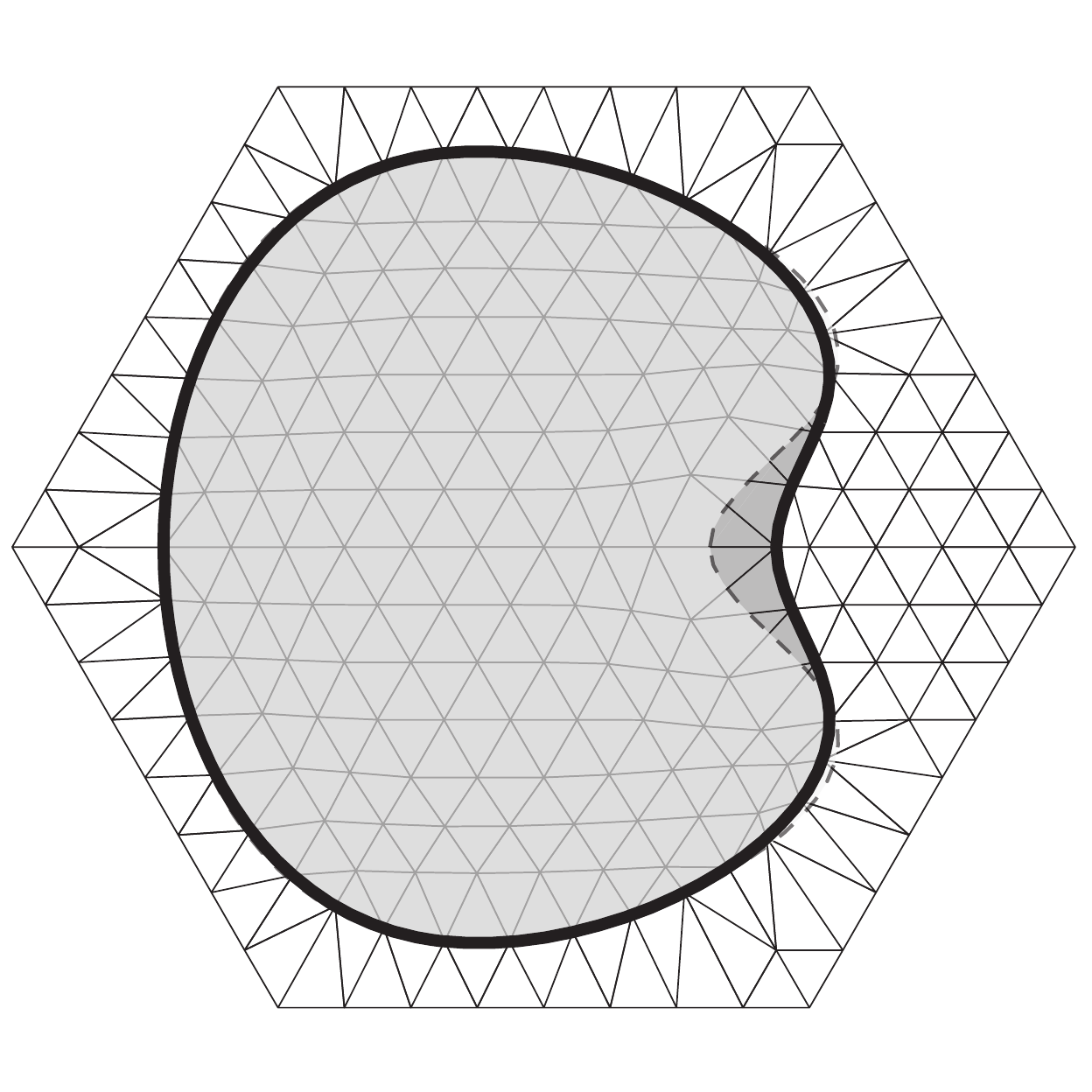}
    \label{fig:step6}
  }
  \caption{Sketch of how the reference domain is periodically
    redefined, and the mesh over it obtained. The triangles
    intersected by the domain in (a) are deformed through the
    universal mesh map to obtain a domain-matching discretization in
    (b). The evolution of the domain during $(t^{n-1},t^n]$ is then
    described through a map $\varphi^t$ defined over
    $\Omega^{t^{n-1}}$. The deformed mesh due to $\varphi^{t^n}$ is
    then shown in (c), where the reference domain $\Omega^{t^{n-1}}$
    is still depicted in light, transparent, gray. These steps are
    then repeated in (d), (e), and (f), for the interval
    $(t^n,t^{n+1}]$. The meshes in (c) and (e) both mesh
    $\Omega^{t^n}$, but since the two differ near the domain
    boundary, a projection of the solution is needed to continue  the integration in time. 
    }
  \label{fig:algorithm}
\end{figure}

There are three main difficulties to overcome in constructing
high-order methods for problems with moving domains: (a) Since the
domain is changing in time, approximations of the domain of the
appropriate order need to be constructed at all times at which the
time-integration scheme is evaluated, (b) the approximation space over
the evolving domain generally needs to evolve in time as well,
resulting in a changing set of degrees of freedom, and (c) the
approximation of time-derivatives of the solution near the evolving
boundary needs to be carefully constructed, since solution values at a
given spatial location may not be defined at all time instants within
a time step.

\paragraph{Pulling back to a reference domain}
A natural approach to sidestep  these issues is to reformulate
the problem as an evolution in a reference, fixed domain $\Omega^0$
through a diffeomorphism $\varphi^t\colon \Omega^0\to \Omega^t$ that
maps it to the evolving domain $\Omega^t$ at each time $t$. If the
solution sought is $u(x,t)$, defined over the domain $\Omega^t$ at each
time $t$, then this approach involves obtaining the partial
differential equation that the function $U(X,t) = u(\varphi^t(X),t)$,
defined over $\Omega^0$ at all times,
would satisfy.  The obvious advantage of this perspective is that any
of the standard numerical methods constructed for evolution problems
on fixed domains can now be applied, and hence high-order methods can
be easily formulated. 

With this idea, the issues associated with discretizing an evolving
domain are transformed into algorithmically constructing and computing
the map $\varphi^t$.  This is not too difficult when the changes in
the domain are small, i.e., when $\varphi^t$ is close to a rigid body
motion for all times. However, it becomes challenging when $\varphi^t$
induces large deformations of the domain. This is the typical problem
of Arbitrary Lagrangian-Eulerian methods: how to deform the mesh, or
alternatively, how to construct the ALE map (see
Fig. \ref{fig:deform}). In terms of the map $\varphi^t$ these same
problems materialize as a loss of local or global injectivity.

A restatement of this same idea from a different perspective is
to consider approximation spaces, such as a finite element spaces, that
evolve with the domain. This is precisely what is obtained if each
function in the approximation space over the reference configuration is pushed forward by
the map $\varphi^t$ at each time $t$. For example, for  finite
element spaces, each shape function  over
$\Omega^t$ has the form $n_a(\varphi^t(X)) = N_a(X)$, where $N_a$ is a
shape function in the finite element space over $\Omega^0$. We take
advantage of this equivalence throughout this manuscript.

\paragraph{Construction of maps}
One of the central ideas we introduce here is one way to construct
maps $\varphi^t$. To circumvent the problems that appear under large
deformations, we periodically redefine the reference domain to be
$\Omega^{t^n}$, $n=0,1,\ldots,N$, $t^n=n \tau$ for some $\tau
>0$, and accordingly $\varphi^t\colon \Omega^{t^{n}}\to \Omega^t$
for $t\in (t^{n},t^{n+1}]$.  

The combination of periodically redefining the reference configuration
and constructing a mesh over  it with the map proposed here is
illustrated in Fig. \ref{fig:algorithm}, for a two-dimensional moving
domain $\Omega^t \subset \mathbb{R}^2$.  Upon choosing a fixed
background triangulation $\mathcal{T}_h$ of a domain $\mathcal{D}
\subset \mathbb{R}^2$ that contains the domains $\Omega^t$ for all $t
\in [0,T]$, $T=N\tau$, the method proceeds as follows: (a) At each
temporal node $t^{n-1}$, a submesh $\mathcal{S}^n_h$ of
$\mathcal{T}_h$ that approximates $\Omega^{t^{n-1}}$
(Fig.~\ref{fig:step1}) is identified; (b)  The polygonal domain meshed by
$\mathcal{S}^n_h$ is deformed through the universal mesh map onto
$\Omega^{t^{n-1}}$ (Fig. \ref{fig:step2}); (c)  The map $\varphi^t$
for $t\in ( t^{n-1},t^{n}]$ is constructed as the identity
everywhere except over the elements with one edge over the moving
boundary. Over these elements $\varphi^t$ is defined as an
extension of the closest point projection of $\partial
\Omega^{t^{n-1}}$ to $\partial \Omega^{t}$. Fig. \ref{fig:step3} shows
the mesh over $\Omega^{t^n}$ obtained as $\varphi^{t^n}(\Omega^{t^{n-1}})$. These three steps are repeated
over $(t^n,t^{n+1}]$, as shown in Figs. \ref{fig:step4},
\ref{fig:step5}, and \ref{fig:step6}. 

\paragraph{Discretization and time integration}
As highlighted earlier, the introduction of the map $\varphi^t$
enables the construction of approximations of any order within each
interval $(t^{n-1},t^n]$, and we elaborate on this next.  

We denote the solution over $(t^{n-1},t^{n}]$ with $U^{n-1}(X,t)$, which takes
values over $\Omega^{t^{n-1}}$ at each time instant in this
interval. To obtain appropriate spatial accuracy, notice that a finite
element space of any order over $\Omega^{t^{n-1}}$
(Fig. \ref{fig:step2}) can be defined in a standard way, by composing
finite element functions over $\mathcal{S}^n_h$ with the universal
mesh map. The spatially discretized equations for $U^{n-1}$ over this
space form an ordinary system of differential equations whose unknowns
are the degrees of freedom for $U^{n-1}$, and hence any standard,
off-the-self integrator of any order can be adopted to approximate its
solution. 

The crucial role played by the universal mesh map is in full
display here, since for smooth domains it provides an exact
triangulation of $\Omega^{t^{n-1}}$. By ensuring that the mesh
conforms exactly to the moving domain at all times, the method is free
of geometric errors -- errors that result from discrepancies between
the exact domain and the computational approximation to the domain.

\paragraph{Projection}
To continue the time integration from the interval $(t^{n-1},t^n]$ to
the interval $(t^n,t^{n+1}]$, an initial condition at $t^n$ is needed,
based on the solution computed in $(t^{n-1},t^n]$. This initial
condition is $ U^{n}(x,t^{n}_+)=\lim_{t\searrow 
  t^{n}} U^{n}(x,t) = U^{n-1}(\left[\varphi^{t^n}\right]^{-1}(x),t)$, which is defined over $\Omega^{t^{n}}$.
 In general,
however, $U^{n}(x,t^{n}_+)$ does not belong to the discrete
approximation space over $\Omega^{t^{n}}$, so we project
$U^{n}(x,t^{n}_+)$ onto it through a suitably defined
projection operator; ideally an $L^2$-projection, but numerical
experiments with interpolation
have rendered very good results as well.

The introduction of
this projection $N$ times would generally have the
detrimental effect of reducing the order of convergence by
one if the spacing $\tau$ between temporal nodes $t^n$ is proportional to the time step $\Delta t$ adopted during integration over each interval $(t^{n-1},t^n]$. Nevertheless, one of the highlights of the map $\varphi^t$ we
construct is that it differs from the identity in a region of
thickness $O(h)$ from the domain boundary, where $h\sim \Delta t$ is
the spatial mesh size. This feature makes the net reduction of the
convergence rate due to the projection to be only of half an order (in the $L^2$ norm).

The implementation of this idea with finite element spaces is
facilitated by regarding this method as a way to construct
approximation spaces that evolve with the domain. This reduces the
effect of the map $\varphi^t$ to defining a ``curved'' mesh over
$\Omega^t$. By further interpolating the map $\varphi^t$ with the
finite element space, an isoparametric approximation of the domain is
obtained. In this way, standard finite element procedures can be
adopted to compute all needed quantities over either the exact or the
isoparametric approximation of $\Omega^t$. This curved mesh is
constructed at each stage of the time-integration scheme. 

\paragraph{Comparison with conventional ALE schemes}

In light of the preceding paragraph, the reader may recognize that our method resembles a conventional ALE scheme with a peculiar mesh motion strategy and regular, systematic ``remeshing.''  In particular, the mesh motion defined by $\varphi^t$ leaves all elements stationary except those with an edge on the moving boundary, and the ``remeshing'' entails the selection of a subtriangulation of a fixed background mesh and perturbing a few of its elements.  

The peculiarity of the approach endows it with several unique
features.  Since the mesh motion is restricted to boundary elements,
the lengths of the time intervals $(t^{n-1},t^n]$ between
``remeshing'' (and hence the time step $\Delta t$ adopted during time
integration over those intervals) are restricted by the mesh spacing; see Section~\ref{sec:error-estim-univ} for details.  An advantage of this strategy is that it easily handles large
domain deformations, and the nodal motions are independent and
explicitly defined.  However, for the reasons described earlier, the
theoretical convergence rate of the method is suboptimal by half an order in the $L^2$ norm.

\paragraph{Remarks}
Back to the difficulties highlighted at the beginning of this section,
it should be evident by now that the basic idea we just outlined
provides approximations of the domain of the proper order at all
times, and that at no point does the difficulty of dealing with nodes
that belong to $\Omega^t$ for only a fraction of the interval
$(t^{n-1},t^n]$ arise. The set of degrees of freedom in the
approximation space does generally change because of the periodic
redefinition of the reference configuration, a seemingly inevitable
step for large enough deformations of the domain, but the introduction
of the projection enables the continuation of the high-order
integration in time with a minimal accuracy loss. We should also
mention that a common difficulty for fixed-mesh methods, which is the
imposition of Dirichlet or Neumann boundary conditions, is handled in
a standard way with the approach in this manuscript. 

In the following, we construct the method in one spatial dimension, to
present some of the main ideas in a rigorous way, yet sidestepping the
notational and algorithmic difficulties introduced by domain
boundaries that are defined by curves instead of isolated points.

\subsection{Construction of the method in one spatial dimension}
\label{sec:constr-meth-one} 

Consider the moving boundary problem: Given a spacetime domain $\Omega=\{(x,t)\in
\mathbb{R}^2 \mid 0<x<s(t),\; 0<t<T\}$, find $u\colon \Omega\to
\mathbb R$ such that
\begin{subequations} \label{mbex}
\begin{align}
0=&\frac{\partial u}{\partial t} - \frac{\partial^2 u}{\partial x^2},
&& (x,t)\in \Omega\label{mbexa} \\
0=&u(0,t) = u(s(t),t), && 0 < t < T \label{mbexb}\\
 u^0(x) = &u(x,t), && 0<x<s(0)
\end{align}
\end{subequations}
where $s \colon [0,T] \rightarrow (0,1)$ is a smooth, prescribed function
of time, and $u^0\colon (0,s(0))\to \mathbb R$ is the initial condition.

For such a problem, it suffices to adopt a grid $0=X_0 < X_1 < \dots < X_M=1$ of the unit interval as the universal mesh -- a stationary background mesh that covers the domains $(0,s(t))$ for all times $0 \le t \le T$.  We shall also employ a partition $0 = t^0 < t^1 < \ldots < t^N = T$ of the time axis that is fine enough so that the change in $s(t)$ over a given interval $(t^{n-1},t^n]$ never exceeds the minimum mesh spacing.  That is,
\[
\max_{t \in (t^{n-1},t^n]} |s(t)-s(t^{n-1})| < \min_{0 < i \le M} (X_i-X_{i-1}).
\]

\begin{figure}[t]
  \centering
  \includegraphics[trim = 0in 0.7in 0in 1in, clip=true, width=0.6\textwidth]{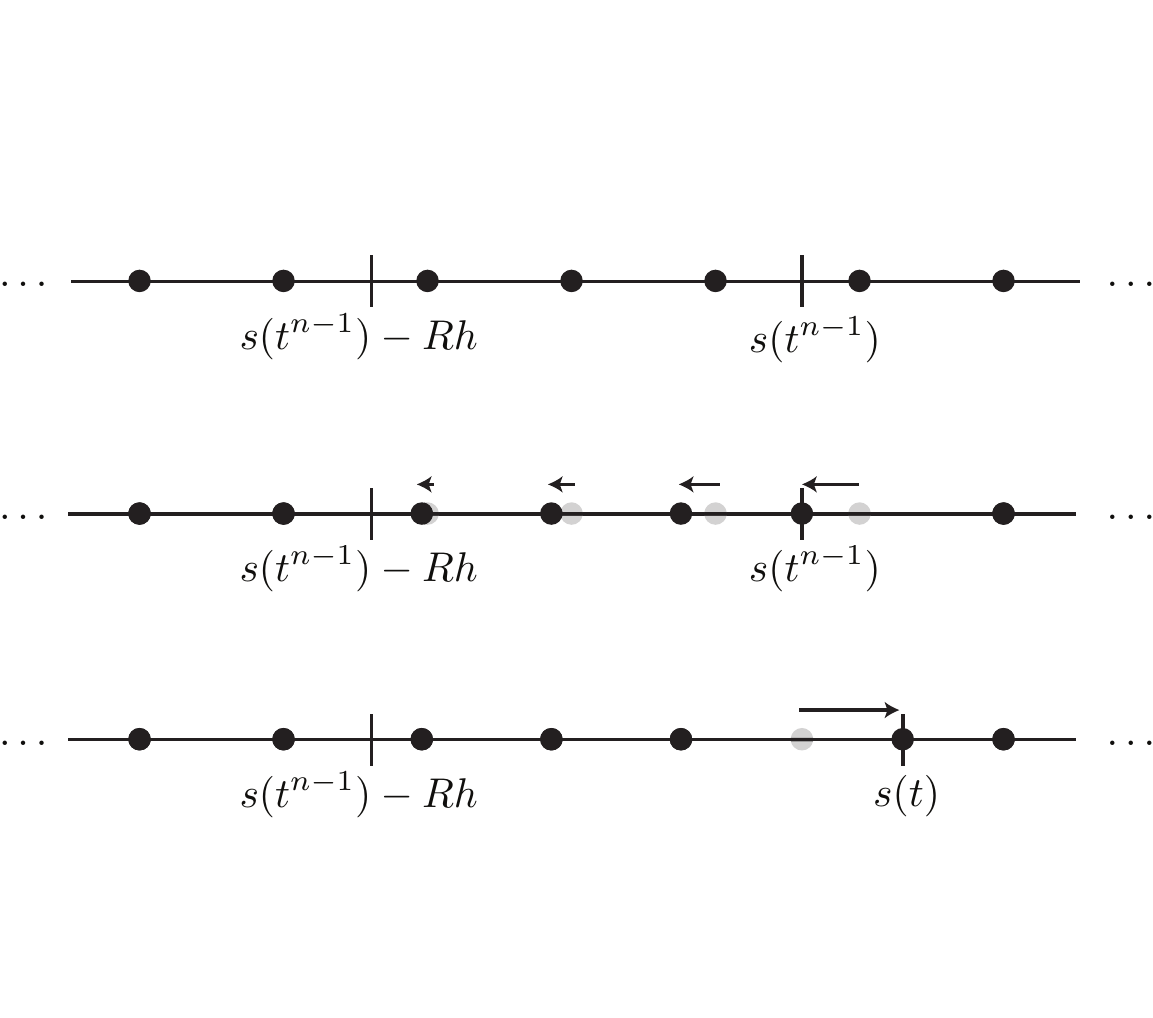}
  \caption{Illustration of the manner in which a one-dimensional universal mesh adapts to the immersed domain $(0,s(t))$ for $t \in (t^{n-1},t^n]$.  At $t=t^{n-1}_+$, the background mesh (top) is deformed by snapping the node that is closest to $s(t^{n-1})$ (among nodes outside the immersed domain) onto $s(t^{n-1})$ (middle).  In the process, the nodes between $s(t^{n-1})-Rh$ and $s(t^{n-1})$ are relaxed away from the boundary. At later times $t \in (t^{n-1},t^n]$ (bottom), the snapped node tracks the position of the boundary, while all other nodes remain in the positions they adopted at $t=t^{n-1}_+$.  Here, we used the map~(\ref{1dmap}) with $R=3$ and $\delta=0.3$. }
  \label{fig:unimesh1d}
\end{figure}

The universal mesh can be adapted to conform exactly to the domain
$(0,s(t))$ at any time $t$ by perturbing nodes in a small neighborhood
of $s(t)$.  A simple prescription  for $t \in (t^{n-1},t^n]$ is, for
each $i$,
\begin{equation} \label{1dmap}
x_i(t) = 
\begin{cases}
X_i- \delta h \left( 1 - \frac{s(t^{n-1}) - X_i}{Rh} \right) &\mbox{if } s(t^{n-1})-Rh \le X_i < s(t^{n-1}) \\
s(t) &\mbox{if } X_{i-1} < s(t^{n-1}) \le X_i \\
X_i &\mbox{otherwise}
\end{cases}
\end{equation}
where $R$ is a small positive integer, $\delta$ is a small positive
number, and $h = \max_{0 < i \le M} (X_i-X_{i-1})$.  See Fig.~\ref{fig:unimesh1d} for an illustration. In this case,
$\varphi^t(X)=\sum_{i=0}^M x_i(t) M_i(X)$, where $M_i$ is the standard
$P_1$ finite element shape function for node $i$: it is affine over
each element and satisfies $M_i(X_j)=\delta_{ij}$.

On this adapted mesh we may construct shape functions
$n_a(x,t) = N_a( (\varphi^t)^{-1}(x))$, where $N_a(X)$ are the shape
functions over the universal mesh. The shape functions $n_a$
are (for instance) piecewise polynomial in $x$ on each interval
$[x_{i-1}(t),x_i(t)]$ for any fixed $t$, and are continuous in $t \in (t^{n-1},t^n)$ for
each fixed $x$. For each
$t\in (t^{n-1},t^n)$, the shape functions $n_a$ satisfy that
\begin{equation}
  \label{eq:2}
  \frac{\partial n_a}{\partial t} (x,t)= -\frac{\partial n_a}{\partial x}(x,t)
  v_h(x,t), 
\end{equation}
where $v_h(\varphi^t(X),t)=\frac{\partial}{\partial t}
\varphi^t(X)$ is the (spatial/Eulerian) velocity of the adapted mesh. For $x_{i-1}(t) < x < x_i(t)$,
\begin{equation*}
v_h(x,t) = \dot{x}_{i}(t) \left( \frac{x-x_{i-1}(t)}{x_{i}(t)-x_{i-1}(t)} \right) + \dot{x}_{i-1}(t) \left( \frac{x_i(t)-x}{x_{i}(t)-x_{i-1}(t)} \right).
\end{equation*}
We then seek an approximate solution
\begin{equation*}
u_h(x,t) = \sum_{a=1}^{A} \mathbf{u}_a(t) n_a(x,t)
\end{equation*}
lying in the space of functions
\begin{equation*}
\mathcal{V}_h(t) = \mathrm{span} \{ n_a(\cdot,t) : n_a(x,t) = 0 \,\forall x>s(t)\}.
\end{equation*}
Here, $\mathbf{u}(t) = (\mathbf{u}_1(t),\mathbf{u}_2(t),\dots,\mathbf{u}_A(t))^T \in \mathbb{R}^A$ is a vector of time-dependent coefficients, which we allow to be discontinuous across the temporal nodes $t^n$.  We denote
\begin{equation*}
\mathbf{u}(t^n_+) = \lim_{t \searrow t^n} \mathbf{u}(t)
\end{equation*}
and similarly for other scalar- or vector-valued functions. To obtain
an equation for $u_h$, we perform a standard Galerkin projection of
(\ref{mbexa}) onto the space of functions $\mathcal{V}_h(t)$, which
leads to the following ordinary differential equation for
$\mathbf{u}$ at each $t\in (t^{n-1},t^n]$,
\begin{equation}
  \label{eq:3}
\mathbf{M}(t) \dot{\mathbf{u}}(t) - \mathbf{B}(t) \mathbf{u}(t)  + \mathbf{K}(t)\mathbf{u}(t) = 0.
\end{equation}
Here $\mathbf{M}(t) \in \mathbb{R}^{A \times A}$ is a mass matrix,
$\mathbf{K}(t) \in \mathbb{R}^{A \times A}$ is a stiffness matrix, and
$\mathbf{B}(t) \in \mathbb{R}^{A \times A}$ is an advection matrix,
constructed according to the following prescription.  For $a$ such
that $n_a(\cdot,t) \in \mathcal{V}_h(t)$,
\begin{align*}
\mathbf{M}_{ab}(t) &= \int_0^1 n_b(x,t) n_a(x,t) \, dx \\
\mathbf{B}_{ab}(t) &= \int_0^1 v_h(x,t) \frac{\partial n_b}{\partial x}(x,t) n_a(x,t) \, dx \\
\mathbf{K}_{ab}(t) &= \int_0^1 \frac{\partial n_b}{\partial x}(x,t) \frac{\partial n_a}{\partial x}(x,t) \, dx,
\end{align*}
while for $a$ such that $n_a(\cdot,t) \notin \mathcal{V}_h(t)$,
\begin{align*}
\mathbf{M}_{ab}(t) &= 0 \\
\mathbf{B}_{ab}(t) &= 0 \\
\mathbf{K}_{ab}(t) &= \delta_{ab}.
\end{align*}
These last values are set so that $u_h(x,t)=0$ for $x>s(t)$, which
follows from imposing (\ref{mbexb}). 
The algorithm then proceeds as follows:

\begin{algorithm}
  \caption{ Time integration for a universal mesh in one dimension. }
  \label{algorithm1d}
\begin{algorithmic}[1]
\Require Initial condition $u(x,0)=u^0(x)$.
\For{$n=1,2,\dots,N$}
\State Project the current numerical solution 
\[
u_h(x,t^{n-1}) = \sum_{a=1}^A \mathbf{u}_a(t^{n-1}) n_a(x,t^{n-1})
\]
\indent (or the initial condition $u(x,0)$ if $n=1$) onto
$\mathcal{V}_h(t^{n-1}_+)$ to obtain the vector of \indent coefficients $\mathbf{u}(t^{n-1}_+)$ in the expansion
\[u_h(x,t^{n-1}_+) = \sum_{a=1}^A \mathbf{u}_a(t^{n-1}_+) n_a(x,t^{n-1}_+).
\]
\State Numerically integrate
\[
\mathbf{M}(t) \dot{\mathbf{u}}(t) - \mathbf{B}(t) \mathbf{u}(t)  + \mathbf{K}(t)\mathbf{u}(t) = 0
\] 
\indent
 for $t\in (t^{n-1},t^n]$
with the initial condition $\mathbf{u}(t^{n-1}_+)$ and the constraints
induced by \indent  (\ref{mbexb}) to obtain $\mathbf{u}(t^n)$.
\EndFor
\State \Return{$u_h(x,t^N)$}
\end{algorithmic}
\end{algorithm}

Several salient features of the method should be evident at this point:
\begin{itemize}
\item \emph{The connectivity of the universal mesh never changes}
  during deformation -- only the nodal positions change.  As a
  consequence, the sizes and sparsity structures of various discrete
  quantities (the solution vector $\mathbf{u}$, the mass matrix
  $\mathbf{M}$, the stiffness matrix $\mathbf{K}$, and the advection
  matrix $\mathbf{B}$) can be held fixed, even though differing
  subsets of degrees of freedom may participate in the discrete
  equations at any interval $(t^{n-1},t^n]$.  One merely needs to
  impose ``homogeneous Dirichlet boundary conditions'' on the solution
  at nonparticipating degrees of freedom.
\item Large deformations of the domain pose no threat to the quality of the deformed mesh, provided $\max_{1 \le n \le N} (t^n-t^{n-1})$ is sufficiently small and the domain evolution is sufficiently regular.
\item In two dimensions, the nodal motions are independent and explicitly defined, rendering the mesh motion strategy low-cost and easily parallelizable.  See Section~\ref{sec:universal_meshes} for details.
\end{itemize}

\section{A Model Moving Boundary Problem} \label{sec:model_problem}

\subsection{The Continuous Problem}
\label{sec:continuous-problem}
\begin{figure}
\centering
\includegraphics[scale=0.35]{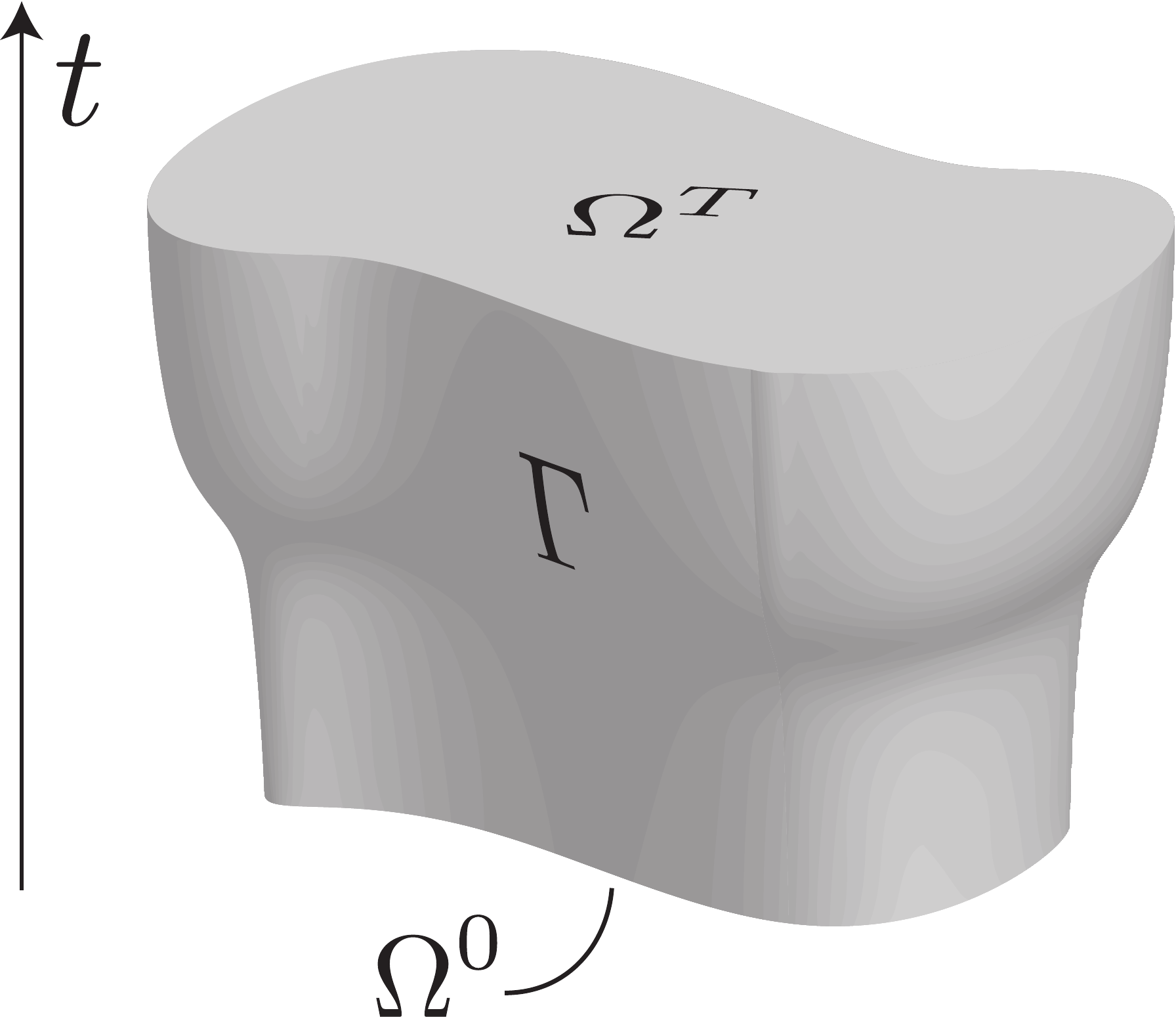}
\caption{Spacetime domain $\Omega$.}
\label{fig:spacetime_domain}
\end{figure}

Consider a moving boundary problem on a bounded spacetime domain $\Omega \subset \mathbb{R}^2 \times [0,T]$, as in Fig.~\ref{fig:spacetime_domain}.  For each $t \in [0,T]$, denote by $\Omega^t \subset \mathbb{R}^2$ the spatial component of the spacetime slice $\Omega \cap (\mathbb{R}^2 \times \{t\})$, and denote by $\Gamma^t$ the boundary of $\Omega^t$.  Finally, let $\Gamma = \bigcup_{0 < t < T} (\Gamma^t \times \{t\})$ denote the lateral boundary of the spacetime domain $\Omega$.  We assume that $\Omega^t$ is open in $\mathbb{R}^2$ for each $t$.
As a regularity requirement, we assume that for every $t \in [0,T]$,
the set $\Gamma^t$ can be expressed as the image of an embedding
$c(\cdot,t)$ of the unit circle $S^1$ into $\mathbb{R}^2$, where $c\in C^2(S^1 \times (0,T),\mathbb{R}^2)$.  

Now consider the following abstract moving boundary problem: Given
$f\colon\Omega\to \mathbb R$ and $u^0\colon \Omega^0\to\mathbb R$,
find $u\colon\Omega\to\mathbb R$ satisfying
\begin{subequations}
\begin{alignat}{3}
&\frac{\partial u}{\partial t} + \mathpzc{a}(u) = f &&\text{  in } \Omega \label{mba} \\
&u = 0 &&\text{  on } \Gamma \label{mbb} \\
&u = u^0 &&\text{  on } \Omega^0 \label{mbc},
\end{alignat} \label{mb}
\end{subequations}
where $\mathpzc{a}$ is a partial differential operator of the form
\begin{equation*}
\mathpzc{a}(u) = -\nabla_x \cdot (k_1\nabla_x u) + k_2 \cdot \nabla_x u + k_3 u
\end{equation*}
with coefficients $k_1(x,t) \in \mathbb{R}^{2 \times 2}$, $k_2(x,t) \in \mathbb{R}^2$, and $k_3(x,t) \in \mathbb{R}$ for every $(x,t) \in \Omega$.  We assume that $k_1$ is uniformly positive definite.  That is, there exists $C>0$ such that $v \cdot k_1(x,t) v \ge C v \cdot v$ for every $v \in \mathbb{R}^2$ and every $(x,t) \in \Omega$.

It is known~\cite[Theorem 7.17]{Lieberman1996} that if $k_1 \in
L^{\infty}(\Omega)^{2 \times 2}$, $k_2 \in L^{\infty}(\Omega)^{2}$,
$k_3 \in L^{\infty}(\Omega)$, the components of $k_1$ are Lipschitz in spacetime, $f
\in L^p(\Omega)$, and $u^0 \in W^{2,p}(\Omega^0)$ with $1<p< \infty$,
then the problem~(\ref{mb}) has a unique solution $u$ with $u(\cdot,t)
\in W^{2,p}(\Omega^t)$ and $\frac{\partial u}{\partial t}(\cdot,t) \in
L^p(\Omega^t)$ for every $0 \le t \le T$.  Here, $W^{s,p}$ denotes the
Sobolev space of differentiability $s \ge 0$ and integrability  $1 \le
p \le \infty$, and $L^p = W^{0,p}$ denotes the Lebesgue space of
integrability $1 \le p \le \infty$.  Later, we shall also denote $H^s
= W^{s,2}$, and we write $H^1_0(\Omega^t)$ for the space of functions
in $H^1(\Omega^t)$ with vanishing trace.  We denote the norm on $W^{s,p}(\Omega^t)$ by $\|\cdot\|_{s,p,\Omega^t}$ and the associated semi-norm by $|\cdot|_{s,p,\Omega^t}$.

\subsection{Equivalent Formulation of the Continuous Problem}

In the following, we derive an equivalent formulation of the moving-boundary problem~(\ref{mb}) that is well-suited for numerical discretization.  For reasons that will soon be made clearer, we restrict our attention to a temporal subinterval $(t^{n-1},t^n] \subset [0,T]$ for the remainder of this section.

\paragraph{Weak formulation}
A weak formulation of~(\ref{mb}) reads: Find $u(\cdot,t) \in \mathcal{V}(\Omega^t) := H^1_0(\Omega^t)$ such that
\begin{equation}
m^t(\dot{u},w) + a^t(u,w) = m^t(f,w) \;\;\; \forall w \in \mathcal{V}(\Omega^t) \label{weak_Eul}
\end{equation}
for every $t \in (t^{n-1},t^n]$, where the time-dependent bilinear forms $m^t$ and $a^t$ are given by
\begin{align*}
m^t(u,w) &= \int_{\Omega^t} uw \, dx \\
a^t(u,w) &= \int_{\Omega^t} \nabla_x w \cdot k_1 \nabla_x u + (k_2 \cdot \nabla_x u) w + k_3 u w \, dx.
\end{align*}
Here and throughout this paper, the dot notation denotes differentiation with respect to time while holding the remaining arguments to the function fixed.

\paragraph{Pulling back to a cylindrical domain}
Given any sufficiently smooth family of bijections $\{\varphi^{n,t} : \Omega^{t^{n-1}} \rightarrow \Omega^t \mid t \in (t^{n-1},t^n]\}$, equation~(\ref{weak_Eul}) may be recast on the cylindrical spacetime domain $\Omega^{t^{n-1}} \times (t^{n-1},t^n]$, since, by a change of variables,~(\ref{weak_Eul}) is equivalent to the statement
\begin{equation}
M^t(\dot{U},W) - B^t(U,W) + A^t(U,W) = M^t(F,W) \;\;\; \forall W \in (\varphi^{n,t})^*\mathcal{V}(\Omega^t) \label{weak_Lag}
\end{equation}
for every $t \in (t^{n-1},t^n]$, where
\begin{equation*}
(\varphi^{n,t})^*\mathcal{V}(\Omega^t) = \left\{ W : \Omega^{t^{n-1}} \rightarrow \mathbb{R} \mid W = w \circ \varphi^{n,t} \text{ for some } w \in \mathcal{V}(\Omega^t) \right\} 
\end{equation*}
is the space of functions in $\mathcal{V}(\Omega^t)$ pulled back to $\Omega^{t^{n-1}}$ by $\varphi^{n,t}$,
\begin{equation*}
\dot{U}(X,t) = \left.\frac{\partial}{\partial t}\right|_X U(X,t),
\end{equation*}
and
\begin{align*}
M^t(U,W) &= \int_{\Omega^{t^{n-1}}} UW \, |\nabla_X\varphi^{n,t}| \, dX \\
B^t(U,W) &= \int_{\Omega^{t^{n-1}}} \left((\nabla_X\varphi^{n,t})^{-\dagger} \nabla_X U \cdot V^{n,t}\right) W |\nabla_X\varphi^{n,t}| dX \\
A^t(U,W) &= \int_{\Omega^{t^{n-1}}} \Big[ \left((\nabla_X\varphi^{n,t})^{-\dagger}\nabla_X W\right) \cdot K_1 \left((\nabla_X\varphi^{n,t})^{-\dagger}\nabla_X U\right) \\
&\hspace{0.7in} + \left((\nabla_X\varphi^{n,t})^{-\dagger} \nabla_X U \cdot K_2\right) W + K_3 UW \Big]
\, |\nabla_X\varphi^{n,t}| \, dX,
\end{align*}
with $|\nabla_X\varphi^{n,t}|$ denoting the absolute value of the Jacobian determinant of $\varphi^{n,t}$ and $(\nabla_X\varphi^{n,t})^{-\dagger}$ denoting the inverse adjoint of $\nabla_X\varphi^{n,t}$.   Here, $K_i = k_i \circ \varphi^{n,t}$, $i=1,2,3$ and $F = f \circ \varphi^{n,t}$ are the Lagrangian counterparts of $k_1,k_2,k_3$, and $f$, and
\[
V^{n,t}(X) := \dot{\varphi}^{n,t}(X) = \left.\frac{\partial}{\partial t}\right|_X \varphi^{n,t}(X)
\]
is the \emph{material} or \emph{Lagrangian velocity}.  

The validity of the preceding change of variables will hold if, for instance,
\begin{equation}
t \mapsto \varphi^{n,t}\in C^1\left( (t^{n-1},t^n], W^{1,\infty}(\Omega^{t^{n-1}})^2 \right) \label{phi_regularity},
\end{equation}
and $(\varphi^{n,t})^{-1}\in W^{1,\infty}(\Omega^t)^2$ for $t\in
(t^{n-1},t^n]$.
Note that under these assumptions, $(\varphi^{n,t})^*\mathcal{V}(\Omega^t) = \mathcal{V}(\Omega^{t^{n-1}}) = H^1_0(\Omega^{t^{n-1}})$.

The presence of the term $B^t(U,W)$ in~(\ref{weak_Lag}) arises from the identity
\begin{equation}
\frac{\partial U}{\partial t}(X,t) = \frac{\partial u}{\partial t}(\varphi^{n,t}(X),t) + \nabla_x u(\varphi^{n,t}(X),t) \cdot v^{n,t}(\varphi^{n,t}(X)), \label{material_time_derivative}
\end{equation}
which relates the partial time derivative of $u$ to the \emph{material time derivative}
\begin{displaymath}
\frac{D u}{D
    t}(\varphi^{n,t}(X),t) := \frac{\partial U}{\partial t}\left(X,t\right)
\end{displaymath}
of $u$ via a term involving the \emph{spatial} or \emph{Eulerian velocity}
\begin{equation*}
v^{n,t}\left(\varphi^{n,t}(X)\right) = V^{n,t}(X).
\end{equation*}
Upon discretization, the term $B^t(U,W)$ corresponds precisely to the term $\mathbf{B}(t)\mathbf{u}(t)$ that the reader encountered earlier in~(\ref{eq:3}).

\paragraph{``Hybrid'' Eulerian formulation} 
A third equivalent statement of (\ref{weak_Eul}) and (\ref{weak_Lag})
is obtained by acknowledging that, by~(\ref{material_time_derivative}),
\begin{equation}
  \label{eq:1b}
  \frac{\partial u}{\partial t}(x,t) = \frac{D u}{D
    t}(x,t) - \nabla_x u(x,t) \cdot v^{n,t}(x).
\end{equation}
It then follows that (\ref{weak_Eul}) is equivalent to 
\begin{equation}
m^t\left(\frac{D u}{D
    t},w\right) - b^t( u,w)+ a^t(u,w) = m^t(f,w) \;\;\; \forall w \in \mathcal{V}(\Omega^t) \label{hybrid_Eul}
\end{equation}
for every $t \in (t^{n-1},t^n]$, where the time-dependent bilinear form $b^t$ is given by
\begin{equation}
  \label{eq:2b}
  b^t(u,w)=  \int_{\Omega^t} \nabla_x u \cdot v^{n,t}\; w\;dx.
\end{equation}
So, $u$ satisfies (\ref{weak_Eul}) if and only if it satisfies
(\ref{hybrid_Eul}) and if and only if $U$ satisfies (\ref{weak_Lag}).
The advantage of this formulation is that it involves simpler
expressions for the bilinear forms than those in (\ref{weak_Lag}),
and these simpler expressions will be convenient for the numerical implementation later.
Notice as well
that the material time derivative on $\partial \Omega^t$ is now a directional
derivative in a direction tangential to the spacetime boundary
$\partial \Omega$, in contrast to $\dot u$, which can only be defined
as a one-side derivative therein.

\section{Discretization} \label{sec:discretization}

\subsection{Spatial Discretization on Short Time Intervals}

At this point it is instructive to derive, in a systematic manner, the
general form of a finite element spatial discretization of~(\ref{mb})
obtained via Galerkin projection.  We begin by spatially discretizing
the weak formulation~(\ref{weak_Eul}) and proceed by pulling the
semidiscrete equations back to a cylindrical spacetime
domain, and by obtaining the ``hybrid'' Eulerian formulation of the
same semidiscrete equations. The utility of these three formulations
will be evident towards the end of this section.

\begin{figure}
\centering
\vspace{0.7in}
\includegraphics[scale=0.6,trim=2in 1.5in 2in 2in]{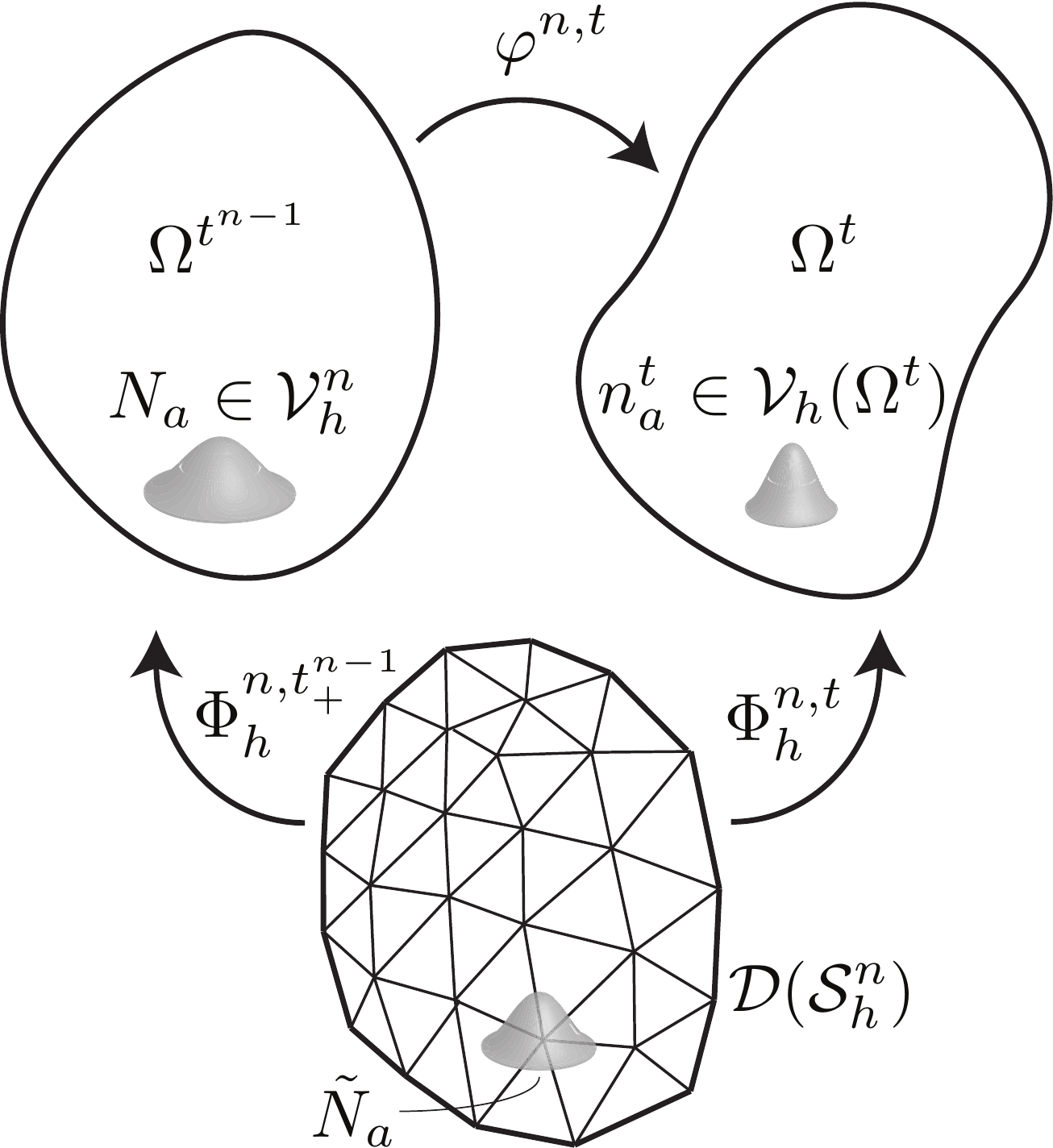}
\vspace{0.9in}
\caption{For each $t \in (t^{n-1},t^n]$, the map $\Phi_h^{n,t}$
  provides a bijection from a fixed reference triangulation
  $\mathcal{S}_h^n$ of a polygonal domain
  $\mathcal{D}(\mathcal{S}_h^n)$ to the moving domain $\Omega^t$.
  Depicted pictorially is a shape function $\tilde{N}_a$ on the
  reference triangulation and its pushforward to $\Omega^{t^{n-1}}$
  and $\Omega^t$, denoted $N_a$ and $n_a^t$, respectively. \label{fig:configurations}}
\end{figure}

\paragraph{Galerkin formulation}
A Galerkin projection of~(\ref{weak_Eul}) requires
choosing a finite-dimensional subspace $\mathcal{V}_h(\Omega^t)
\subset \mathcal{V}(\Omega^t)$ at each time $t$ and finding $u_h(t) \in \mathcal{V}_h(\Omega^t)$ such that
\begin{equation}
m^t(\dot{u}_h,w_h) + a^t(u_h,w_h) = m^t(f,w_h) \;\;\; \forall w_h \in \mathcal{V}_h(\Omega^t) \label{galerkin_Eul}
\end{equation}
for every $t \in (t^{n-1},t^n]$.  For concreteness, let us construct such a family of finite element spaces by fixing a reference triangulation $\mathcal{S}_h^n$ of a polygonal domain $\mathcal{D}(\mathcal{S}_h^n) \subset \mathbb{R}^2$ and constructing a family of continuous, bijective maps 
\begin{equation*}
\Phi_h^{n,t} : \mathcal{D}(\mathcal{S}_h^n) \rightarrow \Omega^t
\end{equation*}
that are differentiable in time and are affine on each triangle $K \in
\mathcal{S}_h^n$, except perhaps near the boundary, see Fig. \ref{fig:configurations}.  In informal language, the image of $\Phi_h^{n,t}$ provides a moving mesh that triangulates $\Omega^t$ for each $t \in (t^{n-1},t^n]$. Then, with $\{\tilde{N}_a\}_a$ denoting shape functions on the reference triangulation, we may set
\begin{equation}
\mathcal{V}_h(\Omega^t) = \mathrm{span}\{n_a^t\}_a \label{shapes_Eul}
\end{equation}
with
\begin{equation*}
n_a^t = \tilde{N}_a \circ (\Phi_h^{n,t})^{-1}
\end{equation*}
for each $t \in (t^{n-1},t^n]$.

\paragraph{Pulling back to a cylindrical domain}
We may pull back the semidiscrete equations~(\ref{galerkin_Eul}) to the cylindrical spacetime domain $\Omega^{t^{n-1}} \times (t^{n-1},t^n]$ with the aid of the bijections
\begin{equation} \label{phih}
\varphi^{n,t} := \Phi_h^{n,t} \circ (\Phi_h^{n,t^{n-1}_+})^{-1}.
\end{equation}
The resulting equivalent semidiscrete equation reads
\begin{equation}
M^t(\dot{U}_h,W_h) - B^t(U_h,W_h) + A^t(U_h,W_h) = M^t(F,W_h) \;\;\; \forall W_h \in (\varphi^{n,t})^*\mathcal{V}_h(\Omega^t) \label{galerkin_Lag}
\end{equation}
for every $t \in (t^{n-1},t^n]$.  

\paragraph{``Hybrid'' Eulerian formulation} Similarly, the discrete ``hybrid''
Eulerian formulation follows by taking advantage of
\eqref{material_time_derivative} to replace $\dot u_h$ in
(\ref{galerkin_Eul}), to get
\begin{equation}
  m^t\left(\frac{Du_h}{D t},w_h\right) -b^t(u_h,w_h)+ a^t(u_h,w_h) = m^t(f,w_h) \;\;\; \forall w_h \in \mathcal{V}_h(\Omega^t) \label{galerkin_hyb-Eul}
\end{equation}
for every $t \in (t^{n-1},t^n]$. 

\paragraph{Remark} We note  that (\ref{galerkin_Eul}),
(\ref{galerkin_Lag}), and (\ref{galerkin_hyb-Eul}) do not define three different
methods; they are three ways of writing precisely the same one.  That is, $u_h$
satisfies~(\ref{galerkin_Eul}) if and only if it satisfies
(\ref{galerkin_hyb-Eul}) and if and only if $U_h(t)=(\varphi^{n,t})^* u_h(t)$ satisfies~(\ref{galerkin_Lag}).

\paragraph{Finite element spaces}
Notice that~(\ref{galerkin_Lag}) is a discretization of~(\ref{weak_Lag}) with a particular choice of a finite element subspace of $\mathcal{V}(\Omega^{t^{n-1}})$, namely $(\varphi^{n,t})^*\mathcal{V}_h(\Omega^t)$.
The shape functions for this space are given by
\begin{align*}
N_a
&= n_a^t \circ \varphi^{n,t} \\
&= \tilde{N}_a \circ (\Phi_h^{n,t^{n-1}_+})^{-1},
\end{align*}
which are time-independent.

As a consequence, the material time derivative of
functions in $\mathcal V_h(\Omega^t)$ takes
a particularly simple form. Let 
\begin{displaymath}
  u_h(\varphi^{n,t}(X),t) = \sum_a u_a(t) n_a^t(\varphi^{n,t}(X)) =
  \sum_a u_a(t) N_a^t(X) = U_h(X,t).
\end{displaymath}
Then
\begin{equation}
  \label{eq:1}
  \frac{D u_h}{Dt}(\varphi^{n,t}(X),t) = \frac{\partial U_h}{\partial
    t}(X,t) = \sum_a \dot u_a(t) N_a^t(X) = \sum_a \dot u_a(t) n_a^t(\varphi^{n,t}(X)),
\end{equation}
since the shape functions $\{N_a\}_a$ do not depend on time. 

Since the map~(\ref{phih}) depends upon $h$, we make that dependence explicit by appending a subscript $h$ to $\varphi^{n,t}$ and all derived quantities ($v^{n,t}$, $V^{n,t}$, $M^t$, $A^t$, and $B^t$) in the remainder of this text.

\paragraph{Summary}
In summary, we have shown that if the semidiscrete equation~(\ref{galerkin_Eul}) is pulled back to the reference domain $\Omega^{t^{n-1}}$ through the use of a map
\begin{equation*}
\varphi^{n,t}_h = \Phi_h^{n,t} \circ (\Phi_h^{n,t^{n-1}_+})^{-1}, 
\end{equation*}
then the resulting semidiscrete equation~(\ref{galerkin_Lag}) involves a finite element space that does not change with time.  We may label that space $\mathcal{V}_h^n$ and write
\begin{equation}
M^t_h(\dot{U}_h,W_h) - B^t_h(U_h,W_h) + A^t_h(U_h,W_h) = M^t_h(F,W_h) \;\;\; \forall W_h \in \mathcal{V}_h^n \label{galerkin_Lag2}
\end{equation}
for every $t \in (t^{n-1},t^n]$.  
The shape functions for $\mathcal{V}_h^n$ are simply shape functions on the reference triangulation $\mathcal{S}_h^n$ pushed forward to $\Omega^{t^{n-1}}$:
\begin{equation*}
N_a = \tilde{N}_a \circ (\Phi_h^{n,t^{n-1}_+})^{-1}.
\end{equation*}

The utility of the above formulation is transparent.  Upon expanding
$U_h$ as a linear combination of shape functions, the
system~(\ref{galerkin_Lag2}) is a system of ordinary differential
equations for the coefficients of the expansion.  This is also evident
from the ``hybrid'' Eulerian formulation (\ref{galerkin_hyb-Eul}) upon
replacing the material time derivative by (\ref{eq:1}).  To this
system of ODEs we may apply a time integrator of choice to advance
from time $t^{n-1}$ to time $t^n$.

\subsection{Integration over Long Time Intervals}

In the preceding sections, we elected to restrict our attention to a temporal subinterval $(t^{n-1},t^n] \subset [0,T]$ and construct finite element subspaces of $\mathcal{V}(\Omega^t)$, $t \in (t^{n-1},t^n]$, using a smoothly varying triangulation of $\Omega^t$ given by the image of $\Phi_h^{n,t}$, $t \in (t^{n-1},t^n]$.  This decision allows for the use of different reference triangulations $\mathcal{S}_h^n$ on different temporal subintervals, simplifying the task of maintaining a nondegenerate triangulation of a domain undergoing large deformations.

To complete the picture and construct an algorithm for integration
over the interval $[0,T]$ of interest, we choose a partition
$0=t^0<t^1<\dots<t^N=T$ and make use of one last ingredient: a linear
projector $p^n_h$ onto $\mathcal{V}^n_h$ for each $n$.
For the definition of the algorithm, we require that the domain of definition of $p^n_h$ contains at least the space $\mathcal{U}^n_h$ given by
\begin{equation*}
\mathcal{U}^n_h =
\begin{cases}
\mathcal{V}(\Omega^0) &\mbox{if } n=1 \\
(\varphi_h^{n-1,t^{n-1}})_*\mathcal{V}_h^{n-1} + \mathcal{V}_h^n &\mbox{if } 1 < n \le N,
\end{cases}
\end{equation*}
where
\begin{equation*}
(\varphi_h^{n-1,t^{n-1}})_*\mathcal{V}_h^{n-1} = \left\{ w : \Omega^{t^{n-1}} \rightarrow \mathbb{R} \mid w \circ \varphi_h^{n-1,t^{n-1}} \in \mathcal{V}_h^{n-1} \right\} 
\end{equation*}
is the space of functions in $\mathcal{V}_h^{n-1}$ pushed forward to $\Omega^{t^{n-1}}$ by $\varphi_h^{n-1,t^{n-1}}$.  We assume the projector is surjective for each $n$; equivalently, $p_h^n \big|_{\mathcal{V}^n_h} = \mathrm{identity}$ for each $n$.  

Some examples of projectors are the orthogonal projector $p^n_{h,L^2}$ onto $\mathcal{V}_h^n$ with respect to the $L^2$-inner product, the orthogonal projector $p^n_{h,H^1}$ onto $\mathcal{V}_h^n$ with respect to the $H^1$-inner product, and the nodal interpolant $i^n_h$ onto $\mathcal{V}_h^n$; see~\cite[Chapter 1]{Ern2004} for details.
The appropriate projector depends on the problem being
approximated and the choice of temporal nodes $t^n$, triangulations $\mathcal{S}_h^n$, and maps $\Phi_h^{n,t}$. As we shall mention, $p_{h,L^2}^n$ is the projector best suited for use with the choices detailed in Section~\ref{sec:universal_meshes}.

With such a family of projectors at hand, a method for
integration over the full time interval $[0,T]$ is then summarized in
Algorithm \ref{algorithm_longtime}.

\begin{algorithm}
  \caption{General form of a time integrator for moving-boundary
    problems with a finite element discretization in space.}
  \label{algorithm_longtime}
\begin{algorithmic}[1]
\Require{Initial condition $u^0 \in \mathcal{V}(\Omega^0)$.}
\For{$n=1,2,\dots,N$}
\vspace{2pt}\State \begin{varwidth}[t]{\linewidth} Choose a reference triangulation $\mathcal{S}_h^n$ and a family of maps $\Phi_h^{n,t} : \mathcal{D}(\mathcal{S}_h^n) \rightarrow \Omega^t$, $t \in (t^{n-1},t^n]$. \end{varwidth}
\vspace{2pt}\State \begin{varwidth}[t]{\linewidth} Generate a finite-dimensional subspace $\mathcal{V}^n_h$ of the continuous solution space $\mathcal{V}(\Omega^{t^{n-1}})$ using shape functions on $\mathcal{S}_h^n$ composed with $(\Phi_h^{n,t^{n-1}_+})^{-1}$. \end{varwidth}
\vspace{2pt}\State \begin{varwidth}[t]{\linewidth} Project the current numerical solution (or the initial condition if $n=1$) onto $\mathcal{V}^n_h$ by setting 
\begin{equation*}
U_h(\cdot,t^{n-1}_+) = p^n_h u_h(\cdot,t^{n-1}), \label{projectedIC}
\end{equation*} 
where $u_h(\cdot,t^0)=u^0$ or, for $n>1$,
\begin{equation*}
u_h(x,t^{n-1}) = U_h((\varphi^{n-1,t^{n-1}}_h)^{-1}(x),t^{n-1})
\end{equation*}
is the pushforward of $U_h(\cdot,t^{n-1}) \in \mathcal{V}_h^{n-1} \subset \mathcal{V}(\Omega^{t^{n-2}})$ to $\Omega^{t^{n-1}}$. \end{varwidth} \label{projection_step}
\vspace{2pt}\State \begin{varwidth}[t]{\linewidth} Numerically integrate~(\ref{galerkin_Lag2}) over $(t^{n-1},t^n]$ with the projected initial condition $U_h(\cdot,t^{n-1}_+)$. \label{integration_step}  \end{varwidth}
\EndFor
\State \Return{$u_h(\cdot,t^N)$}
\end{algorithmic}
\end{algorithm}

\paragraph{Relationship to ALE} Let us emphasize that Algorithm \ref{algorithm_longtime} has been formulated with enough generality that it encompasses not only the method specific to this paper involving universal meshes (which is detailed in Section~\ref{sec:universal_meshes}) but also conventional ALE schemes.  In the case of an ALE scheme, the reference triangulation $\mathcal{S}_h^n$ is a triangulation of $\Omega^{t^{n-1}}$, the map $\varphi_h^{n,t}$ corresponds to a mesh motion derived from, e.g., solutions to the equations of linear elasticity, and the temporal nodes $t^n$ correspond to times at which remeshing is performed.  In the case of the method specific to this paper, we shall see in Section~\ref{sec:universal_meshes} that the reference triangulation $\mathcal{S}_h^n$ is a subtriangulation of a fixed background mesh, the map $\varphi_h^{n,t}$ induces deformations of triangles on the boundary of $\mathcal{S}_h^n$ while leaving the remaining triangles fixed, and the temporal nodes $t^n$ are spaced closely enough so that these deformations of boundary triangles remain well-behaved.

\subsection{Example: a Runge-Kutta Time-Integrator}
\label{sec:example:-runge-kutta}

We next exemplify how a time integrator of any given order can be
incorporated into step \ref{integration_step} of the algorithm. In this
case we consider  an $s$-stage Singly Diagonally Implicit Runge-Kutta (SDIRK) method of order $\le s$  as the time integrator~\cite{Burrage1980,Hairer2002}.  Such an integrator requires solving a sequence of $s$ systems of equations
\begin{equation}
M_h^{t_i}(U_i,W) = M_h^{t_i}\left(\sum_{j=0}^{i-1} \beta_{ij} U_j, W\right) + \gamma\Delta t G_h^{t_i}(U_i,F(t_i);W) \;\;\; \forall W \in \mathcal{V}_h^n \label{dirk}
\end{equation}
for $U_i \in \mathcal{V}_h^n$, $i=1,2,\dots,s$, where $U_0 = U_h(\cdot,t_0)$, $t_0 \in (t^{n-1},t^n]$, $t_i = \sum_{j=0}^{i-1} \beta_{ij} t_j + \gamma\Delta t$ for $0 < i \le s$, and 
\begin{equation*}
G^t_h(U,F;W) = M_h^t(F,W) - A_h^t(U,W) + B_h^t(U,W).
\end{equation*}
The time-$\Delta t$ advancement of $U_0$ is then given by $U_s$.  The coefficients $\gamma>0$ and $\beta_{ij} \in \mathbb{R}$, $i=1,2,\dots,s$, $j=0,1,\dots,i-1$, for various SDIRK methods are tabulated in~\ref{sec:appendix_sdirk}, Tables~\ref{tab:SDIRK2}-\ref{tab:SDIRK4}.  Pragmatically, implementing an SDIRK method amounts to computing $s$ ``backward-Euler'' steps, with the initial condition at the $i^{th}$ stage given by a linear combination of the solutions at the previous stages.

\subsection{Overview of Error Estimates}

In our companion paper~\cite{Gawlik2012b}, we derive error estimates
in the $L^2$-norm for the aforementioned method for the problem in \S \ref{sec:continuous-problem}.  Here, we give an
overview of the estimates for the case in which the finite element
spaces $\mathcal{V}_h^n$ consist of continuous functions made of
elementwise polynomials of degree $r-1$, where $r>1$ is an integer.  We begin by introducing some notation.

\paragraph{Notation}  
Let $u^n \in \mathcal{V}(\Omega^{t^n})$ denote the value of the exact solution $u$ at $t=t^n$, i.e. $u^n = u(\cdot, t^n)$, and let $u_h^{\Delta t,n} \in \mathcal{V}_h(\Omega^{t^n})$ denote the value of the fully discrete solution at $t=t^n$.  Finally, let $v_h : \Omega \rightarrow \mathbb{R}^2$ denote the vector field on the spacetime domain $\Omega$ whose restriction to each temporal slice is $v_h^{n,t}$, i.e. $v_h(\cdot,t) = v_h^{n,t}$ for $t \in (t^{n-1},t^n]$. 

From this point forward, the parameter $h$ denotes the maximum
diameter of a triangle belonging to $\{ K \in \mathcal{S}_h^n\mid 1
\le n \le N \}$. Additionally, let $\Delta t$ be the maximum time step
adopted while time-integrating over the interval $[0,T]$, namely, the
maximum time step employed by the time-integrator in
Line~\ref{integration_step} of Algorithm~\ref{algorithm_longtime}
among all intervals $(t^{n-1},t^n]$.  We remind the reader that the temporal nodes $t^n$ demarcate changes in the reference triangulation $\mathcal{S}_h^n$; hence, the time step adopted during integration over $(t^{n-1},t^n]$ is  less than or equal to $t^n-t^{n-1}$ for every $n$.  We assume that the time integrator employed during these intervals is stable and has a global truncation error of order $q \ge 1$ in the time step $\Delta t$. 

In Line~\ref{projection_step} of of Algorithm~\ref{algorithm_longtime}, the numerical solution is transferred, via a projection $p^n_h$, between two finite element spaces associated with differing triangulations of $\Omega^{t^{n-1}}$.  We denote by $\mathcal{R}_h^n \subseteq \Omega^{t^{n-1}}$ the region over which the two triangulations differ, and by $|\mathcal{R}_h^n|$ its (Lebesgue) measure.  We assume the projector is stable in the sense that there exists a constant $C_p$ independent of $h$ and $n$ such that
  \begin{equation*}
    \| p^n_h U\|_{0,2,\Omega^{t^{n-1}}} \le  C_p \| U\|_{0,2,\Omega^{t^{n-1}}}
  \end{equation*}
for all $U\in \mathcal{U}_h^n$.

\paragraph{General error estimate}
It is proven in~\cite{Gawlik2012b} that if the assumptions above are satisfied, the triangulations $\mathcal{S}_h^n$ are quasi-uniform, and the exact solution $u$ and the maps $\Phi_h^{n,t}$ are sufficiently regular, then an error estimate of the following form holds with constants $C_1(u,v_h,T)$, $C_2(u,v_h,T)$ and $C_3(u,T)$:
\begin{equation}
\|u_h^{\Delta t,N} - u^N\|_{0,2,\Omega^T} \le C_p^N \left( C_1(u,v_h,T) h^r + C_2(u,v_h,T) \Delta t^q + C_3(u,T) h^r \ell_{h,r}  \sum_{n=1}^N |\mathcal{R}_h^n|^{1/2} \right) \label{error_general}
\end{equation}
where
\begin{equation} \label{logh}
\ell_{h,r} =
\begin{cases}
\log (h^{-1}) &\mbox{if } r=2 \\
1 &\mbox{if } r > 2
\end{cases}
\end{equation}
and $C_1(u,v_h,T), C_2(u,v_h,T) \ge C(u,T)$ for some constant $C(u,T)>0$.

The content of this estimate is easily understood.  The error
committed by the method consists of three terms, amplified by the
$N^{th}$ power of the projector's stability constant: an error due to
spatial discretization of order at best $h^r$ (first term), an error due to
temporal discretization of order at best $\Delta t^q$ (second term), and an error
introduced by projecting between differing triangulations of the same
domain $\Omega^{t^n}$ at each temporal node $t^n$ (third term).  The coefficients
of the first two terms depend implicitly on $h$ through the choice of
the domain velocity $v_h$.  The precise scaling of these coefficients
with respect to $h$ depends upon the chosen mesh motion strategy and
is, of course, no better than $O(1)$ in $h$. We particularize this
estimate for the mesh motion strategy proposed here in \S \ref{sec:error-estim-univ}.

\paragraph{Discussion}

Error estimates that are specific to two categories of methods are immediately apparent from the general estimate~(\ref{error_general}).  The first category consists of classical ALE schemes with occasional remeshing -- that is, $N$ is independent of $h$ and $\Delta t$.  For methods of this type, the amplifier $C_p^N$ is of order unity (regardless of the choice of the projector), the summation of the mesh discrepancy volumes is of order unity, and the coefficients $C_1(u,v_h,T)$ and $C_2(u,v_h,T)$ can be bounded independently of $h$ for sufficiently regular mesh motion strategies.  The resulting error estimate reads
\begin{equation*}
\|u_h^{\Delta t,N} - u^N\|_{0,2,\Omega^T} \le C(u,T) (h^r + \Delta t^q + \ell_{h,r} h^r)  
\end{equation*}
for a constant $C(u,T)$ (it is straightforward to sharpen this estimate to $\|u_h^{\Delta t,N} - u^N\|_{0,2,\Omega^T} \le C(u,T) (h^r + \Delta t^q)$).

At the other extreme are methods for which the reference triangulation $\mathcal{S}_h^n$ is updated more frequently, e.g., at intervals proportional to $\Delta t$.   This strategy is, in fact, the one adopted in the method proposed in Section~\ref{sec:universal_meshes}.  For methods of this type, the proportionality between $N$ and $\Delta t^{-1}$ mandates the use of a projector with stability constant $C_p=1$ (such as the $L^2$-projector $p_{h,L^2}^n$).  However, the short time intervals between updates of the reference triangulation allow for the use of simple mesh motion strategies in which the nodal motions are independent, explicitly defined, and restricted to only nodes that lie on the moving boundary.  The resulting mesh discrepancy volumes $|\mathcal{R}_h^n|$ are of order $h$, and the coefficients $C_1(u,v_h,T)$ and $C_2(u,v_h,T)$ can be shown to be of order $\ell_{h,r} h^{-1/2}$ and of order unity, respectively, for suitable nodal motions.  The ultimate error estimate reads
\begin{equation}   \label{eq:5}
\|u_h^{\Delta t,N} - u^N\|_{0,2,\Omega^T} \le C(u,T) (\ell_{h,r} h^{r-1/2} + \Delta t^q + \ell_{h,r} h^{r+1/2} \Delta t^{-1})
\end{equation}
for a constant $C(u,T)$, which is suboptimal by half an order when $\Delta t \sim h$.  Note that in practice, we have observed in numerical experiments that the use of a projector with stability constant $C_p>1$ (such as the interpolation operator $i_h^n$) does not lead to a degradation of convergence beyond the existing half-order suboptimality, despite the theory's predictions.

\section{Universal Meshes} \label{sec:universal_meshes}

The algorithm presented in the preceding section requires at each temporal node $t^{n-1}$ the selection of a family of maps $\Phi_h^{n,t} : \mathcal{D}(\mathcal{S}_h^n) \rightarrow \Omega^t$, $t \in (t^{n-1},t^n]$, from a fixed polygonal domain $\mathcal{D}(\mathcal{S}_h^n)$ to the moving domain $\Omega^t$.  Here we present a means of constructing such maps using a single, \emph{universal mesh} that triangulates an ambient domain $\mathcal{D} \subset \mathbb{R}^2$ containing the domains $\{\Omega^t\}_{t=0}^T$ for all times $t \in [0,T]$.  Full details of the method are described in~\cite{Rangarajan2012a}.

The essence of the method is to triangulate $\mathcal{D}$ with a fixed
mesh $\mathcal{T}_h$ and to identify, for each time interval $(t^{n-1},t^n]$, a submesh $\mathcal{S}^n_h$ of $\mathcal{T}_h$ that approximates $\Omega^{t^{n-1}}$.  Triangles on the boundary of $\mathcal{S}^n_h$ are then deformed in such a way that the submesh conforms exactly to the moving domain $\Omega^t$ for all $t \in (t^{n-1},t^n]$.

The conditions under which a given triangulation $\mathcal{T}_h$ can
be so adapted to conform to a family of domains $\Omega^t$, $t \in
[0,T]$, are laid forth in~\cite{Rangarajan2012b,Rangarajan2012a}.  Briefly, the procedure is guaranteed to succeed if:
\begin{enumerate}[(i)]
\item $\Omega^t$ is $C^2$-regular for every $t$. \label{unicond1}
\item $\mathcal{T}_h$ is sufficiently refined in a neighborhood of $\partial\Omega^t$ for every $t$. \label{unicond2}
\item All triangles in $\mathcal{T}_h$ have angles bounded above by a constant $\vartheta < \pi/2$. \label{unicond3}
\end{enumerate}

The level of refinement requested by condition~(\ref{unicond2}) is dictated primarily by the minimum radius of curvature of $\partial \Omega^t$ among all times $t \in [0,T]$, which, roughly speaking, must be no less than a small multiple of the maximum element diameter.  This notion is made precise in~\cite{Rangarajan2012b}.  Note that condition~(\ref{unicond1}) precludes an application of the method in its present form to domains with corners.

\subsection{Construction of an Exactly Conforming Mesh} \label{sec:construction}

In detail, consider a triangulation $\mathcal{T}_h$ of $\mathcal{D}$ satisfying conditions~(\ref{unicond1}-\ref{unicond3}), with the parameter $h$ denoting the length of the longest edge in the triangulation. For a given domain $\Omega^t \subset \mathcal{D}$, $t \in [0,T]$, let $\phi^t : \mathcal{D} \rightarrow \mathbb{R}$ denote the signed distance function to $\partial \Omega^t$, taken to be positive outside $\Omega^t$ and negative inside $\Omega^t$.  Let $\pi^t : \mathcal{D} \rightarrow \partial\Omega^t$ denote the closest point projection onto $\partial\Omega^t$.  For $i=0,1,2,3$, let $\mathcal{T}_{h,i}^t$ denote the collection of triangles $K \in \mathcal{T}_h$ for which exactly $i$ vertices of $K$ do not lie in the interior of $\Omega^t$.

For a given subtriangulation $\mathcal{S}_h$ of $\mathcal{T}_h$, we make the distinction between $\mathcal{S}_h$, the list of vertices in the subtriangulation and their connectivities, and $\mathcal{D}(\mathcal{S}_h)$, the polygonal domain occupied by triangles in $\mathcal{S}_h$.  We write $K \in \mathcal{S}_h$ to refer to triangles $K \subseteq \mathcal{D}(\mathcal{S}_h)$ who have vertices in $\mathcal{S}_h$.

\begin{figure}[t]
\hspace{-1.5in}
\centering
\includegraphics[scale=0.4,trim=0in 1in 3.5in 2in]{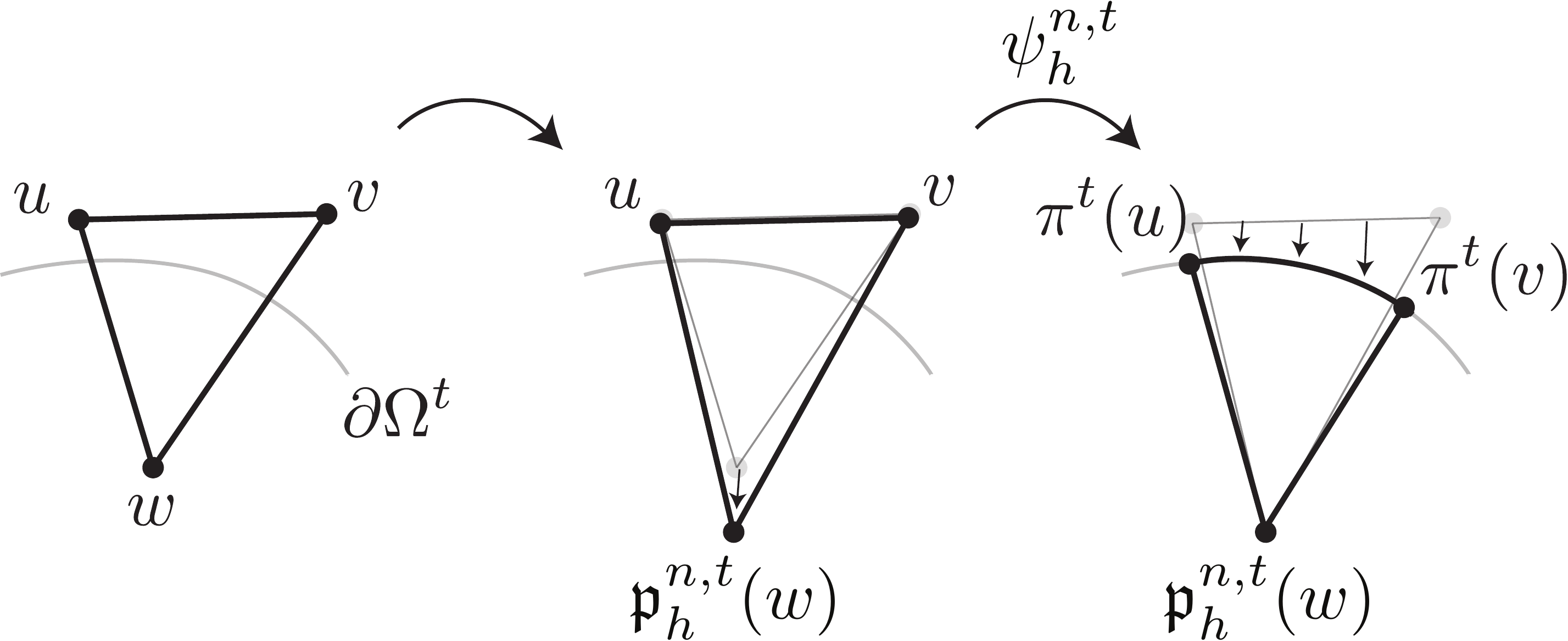}
\vspace{0.5in}
\caption{The action of $\Phi^{n,t}_h$ on a triangle $K \in \mathcal{T}_{h,2}^{t^{n-1}}$ comprises two steps: A relaxation step that moves $w$ away from the boundary, and a nonlinear blend map $\psi^{n,t}_h$ that maps the straight triangle to a curved one.}
\label{fig:blend}
\end{figure}

To construct a conforming mesh for $\Omega^t$ from the mesh $\mathcal{T}_h$, we choose
\begin{equation*}
\mathcal{S}^n_h = \mathcal{T}_{h,0}^{t^{n-1}} \cup \mathcal{T}_{h,1}^{t^{n-1}} \cup \mathcal{T}_{h,2}^{t^{n-1}}
\end{equation*}
as the reference subtriangulation for the domains $\Omega^t$, $t \in
(t^{n-1},t^n]$.  This subtriangulation is simply the set of triangles
in $\mathcal{T}_h$ with at least one vertex in $\Omega^{t^{n-1}}$.
The map $\Phi_h^{n,t} : \mathcal{D}(\mathcal{S}_h^n) \rightarrow
\Omega^t$ will then make use of three important mappings, described in
the following paragraphs, and illustrated in Fig. \ref{fig:blend}. The
{\it universal mesh map}, as described in \cite{Rangarajan2012a}, is
$\Phi_h^{n,t^{n-1}}$. 

\paragraph{Boundary evolution map}  The first is a \emph{boundary evolution map} $\gamma^{n,t}_h : \partial\mathcal{D}(\mathcal{S}_h^n) \rightarrow \partial\Omega^t$, which provides a correspondence between the piecewise linear boundary of $\mathcal{D}(\mathcal{S}_h^n)$ and the boundary of $\Omega^t$ for $t \in (t^{n-1},t^n]$, as in Fig.~\ref{fig:blend}.  The choice of $\gamma^{n,t}_h$ is not unique, although a simple choice is the closest point projection onto $\Omega^t$ composed with the closest point projection onto $\Omega^{t^{n-1}}$: 
\begin{equation}
\gamma^{n,t}_h = \pi^t \circ \pi^{t^{n-1}}\big|_{\partial\mathcal{D}(\mathcal{S}_h^n)}. \label{gamma1}
\end{equation}
By the regularity of the spacetime domain $\Omega$, this map is
well-defined for $h$ sufficiently small and $t$ sufficiently close to
$t^{n-1}$; see \cite{Rangarajan2012b}.

\paragraph{Relaxation map} The second is a \emph{relaxation map} $\mathfrak{p}^{n,t}_h$ that perturbs vertices lying both inside $\Omega^t$ and near $\partial \Omega^t$ in a direction away from $\partial \Omega^t$.  A simple choice of relaxation is the map
\begin{equation}
\mathfrak{p}^{n,t}_h(x) = 
\begin{cases}
x - \delta h \left(1+\frac{\phi^{t^{n-1}}(x)}{Rh}\right)\nabla\phi^{t^{n-1}}(x) & \text{if } -Rh < \phi^{t^{n-1}}(x) < 0 \\
x & \text{otherwise},
\end{cases} \label{relaxation}
\end{equation}
which moves vertices within a distance $Rh$ of
$\partial\Omega^{t^{n-1}}$ by an amount $\le \delta h$ in a direction
normal to the boundary, with $R > 1$ a small positive integer and
$(1+1/R)^{-1} \le \delta \le 1$.  It is proven
in~\cite{Rangarajan2012a} that for a straight boundary (or one of
small enough radius of curvature compared with the mesh size) such a map results in elements of
bounded quality at $t=t^{n-1}$ when
conditions~(\ref{unicond1}-\ref{unicond3}) hold.

Note that this choice of relaxation leaves relaxed vertices fixed over the duration of the interval $(t^{n-1},t^n]$.  We denote by $\mathfrak{p}^{n,t}_h(\mathcal{T}_h)$ the triangulation obtained by applying the relaxation $\mathfrak{p}^{n,t}_h$ to the vertices of $\mathcal{T}_h$ while preserving the mesh's connectivity.

\paragraph{Blend map} Finally, we will make use of a \emph{blend map} $\psi^{n,t}_h$ which takes a straight triangle $K \in \mathfrak{p}^{n,t}_h\left(\mathcal{T}_{h,2}^{t^{n-1}}\right)$ to a curved triangle that conforms exactly to the boundary.  The map we employ is proposed in~\cite{Rangarajan2012a}.  Letting $u,v,w$ denote the vertices of $K$, the blend map reads
\begin{align}
\psi^{n,t}_h(x) = \nonumber
&\frac{1}{2(1-\lambda_u)} [\lambda_v\gamma^{n,t}_h(\lambda_u u + (1-\lambda_u)v) + \lambda_u\lambda_w \gamma^{n,t}_h(u)] \nonumber \\ 
&+ \frac{1}{2(1-\lambda_v)} [\lambda_u\gamma^{n,t}_h((1-\lambda_v)u+\lambda_v v) + \lambda_v\lambda_w \gamma^{n,t}_h(v)] + \lambda_w w, \label{blend}
\end{align}
where $\lambda_u,\lambda_v,\lambda_w$ are the barycentric coordinates of $x \in K$.  Here, we have employed the convention the vertex $w$ is the unique vertex of $K$ lying inside $\Omega^{t^{n-1}}$.  It is not difficult to check that for fixed $t$, the blend map $\psi^{n,t}_h$ maps points $x$ lying on the edge $uv$ to their images under the boundary evolution map $\gamma^{n,t}_h$, preserves the location of the vertex $w$, and is affine on the edges $wu$ and $wv$.

\paragraph{Culmination} We now define  $\Phi^{n,t}_h$ over each triangle $K \in \mathcal{S}_h^n$ with vertices $u,v,w$ according to
\begin{equation}
\Phi^{n,t}_h(x) = 
\begin{cases}
\lambda_u \mathfrak{p}^{n,t}_h(u) + \lambda_v \mathfrak{p}^{n,t}_h(v) + \lambda_w \mathfrak{p}^{n,t}_h(w) & \text{if } K \in \mathcal{T}_{h,0}^{t^{n-1}} \\
\lambda_u \gamma^{n,t}_h(u) + \lambda_v \mathfrak{p}^{n,t}_h(v) + \lambda_w \mathfrak{p}^{n,t}_h(w) & \text{if } K \in \mathcal{T}_{h,1}^{t^{n-1}} \\
\psi^{n,t}_h(\lambda_u u + \lambda_v v + \lambda_w \mathfrak{p}^{n,t}_h(w))
& \text{if } K \in \mathcal{T}_{h,2}^{t^{n-1}},
\end{cases} \label{meshmap}
\end{equation}
where $\lambda_u,\lambda_v,\lambda_w$ are the barycentric coordinates of $x \in K$.  Once again, we have employed the convention that for triangles $K \in \mathcal{T}_{h,2}^{t^{n-1}}$, the vertex $w$ is the unique vertex of $K$ lying inside $\Omega^{t^{n-1}}$, and for triangles $K \in \mathcal{T}_{h,1}^{t^{n-1}}$, the vertex $u$ is the unique vertex of $K$ lying outside $\Omega^{t^{n-1}}$.

\paragraph{The domain evolution and its velocity}

It is now straightforward to record explicit expressions for the domain mapping $\varphi^{n,t}_h$ and its material velocity $V_h^{n,t}$.  By definition,
\begin{equation}
\varphi^{n,t}_h = \Phi^{n,t}_h \circ \left(\Phi^{n,t^{n-1}_+}_h \right)^{-1}. \label{domain_flow}
\end{equation}
The velocity field $V^{n,t}_h$ is then given by differentation with respect to time:
\begin{equation*}
V^{n,t}_h = \dot{\Phi}^{n,t}_h \circ \left(\Phi^{n,t^{n-1}_+}_h \right)^{-1}.
\end{equation*}
If the relaxation map $\mathfrak{p}^{n,t}_h$ is independent of time over $(t^{n-1},t^n]$ (as is the case for the choice~(\ref{relaxation})), this expression for $V^{n,t}_h$ is given explicitly by
\begin{equation}
V^{n,t}_h(X) = 
\begin{cases}
0 & \text{if } K \in \mathcal{T}_{h,0}^{t^{n-1}} \\
\lambda_u \dot{\gamma}^{n,t}_h(u) & \text{if } K \in \mathcal{T}_{h,1}^{t^{n-1}} \\
\begin{split}
& \frac{\lambda_v}{2(1-\lambda_u)} \dot{\gamma}^{n,t}_h (\lambda_u u + (1-\lambda_u) v)  + \frac{\lambda_u\lambda_w}{2(1-\lambda_u)} \dot{\gamma}^{n,t}_h(u)  \\
& \;\;\;+ \frac{\lambda_u}{2(1-\lambda_v)} \dot{\gamma}^{n,t}_h ((1-\lambda_v) u + \lambda_v v) + \frac{\lambda_v\lambda_w}{2(1-\lambda_v)} \dot{\gamma}^{n,t}_h(v) 
\end{split}
& \text{if } K \in \mathcal{T}_{h,2}^{t^{n-1}},
\end{cases} \label{vexact}
\end{equation}
where $\lambda_u,\lambda_v,\lambda_w$ are the barycentric coordinates of $(\Phi_h^{n,t^{n-1}_+})^{-1}(X) \in K$, with the conventional ordering of the vertices described earlier.  Formulas for the time derivative of $\pi^t$ (which are needed for the choice $\gamma_h^{n,t}=\pi^t \circ \pi^{t^{n-1}}$) in terms of local measures of the boundary's shape and velocity are given in~\ref{sec:appendix_pidot}.

\subsection{Alternative: Isoparametric Approximation of the Domain}
\label{sec:isop-appr-doma}

A convenient alternative to exact representations of the domain is to
adopt superparametric or isoparametric representations of the
domain. This entails approximating the map $\Phi^{n,t}_h$ (and hence
the domain $\Omega^t$) with a polynomial interpolant
\begin{equation}
\Phi^{n,t}_{h,\mathrm{approx}}(\tilde{X}) = \sum_a \tilde{M}_a(\tilde{X}) \Phi_h^{n,t}(\tilde{Y}_a) \label{Phi_approx}
\end{equation}
constructed from shape functions $\tilde{M}_a$ of a triangular Lagrange element (henceforth termed Lagrange shape functions) with corresponding degrees of freedom $\tilde{Y}_a$ on the reference triangulation $\mathcal{S}_h^n$.  In this way, expressions for the spatial derivatives of the corresponding shape functions
\begin{equation}
N_{a,\mathrm{approx}} = \tilde{N}_a \circ \left(\Phi^{n,t^{n-1}_+}_{h,\mathrm{approx}} \right)^{-1} \label{Naapprox}
\end{equation}
and
\begin{equation*}
n^t_{a,\mathrm{approx}} = \tilde{N}_a \circ \left(\Phi^{n,t}_{h,\mathrm{approx}} \right)^{-1}
\end{equation*}
involve only derivatives of the reference triangulation's shape functions $\tilde{N}_a$ and the Lagrange shape functions $\tilde{M}_a$, and not the gradients of the exact map $\Phi_h^{n,t}$:
\begin{align*}
&\nabla_X N_{a,\mathrm{approx}}(X) = \nabla_{\tilde{X}} \tilde{N}_a(\tilde{X}) \cdot \left(\nabla_{\tilde{X}} \Phi^{n,t^{n-1}_+}_{h,\mathrm{approx}}\right)^{-1} = \nabla_{\tilde{X}} \tilde{N}_a(\tilde{X}) \cdot \left(\sum_a \nabla_{\tilde{X}} \tilde{M}_a(\tilde{X}) \Phi_h^{n,t^{n-1}_+}(\tilde{Y}_a) \right)^{-1} \\
&\nabla_x n^t_{a,\mathrm{approx}}(x) = \nabla_{\tilde{X}} \tilde{N}_a(\tilde{X}) \cdot \left(\nabla_{\tilde{X}} \Phi^{n,t}_{h,\mathrm{approx}}\right)^{-1} = \nabla_{\tilde{X}} \tilde{N}_a(\tilde{X}) \cdot \left(\sum_a \nabla_{\tilde{X}} \tilde{M}_a(\tilde{X}) \Phi_h^{n,t}(\tilde{Y}_a) \right)^{-1}.
\end{align*}  
This, in turn, eliminates the need to compute gradients of the closest
point projection $\pi^t$. This, and other reasons detailed later, make
approximating the domain in this way more computationally convenient in practice. 

For completeness, we next detail the corresponding approximate
domain map 
\begin{equation*}
\varphi^{n,t}_{\mathrm{approx}} = \Phi^{n,t}_{h,\mathrm{approx}} \circ \left(\Phi^{n,t^{n-1}_+}_{h,\mathrm{approx}}\right)^{-1}
\end{equation*} 
and velocity fields, which take particularly simple
forms. In fact, with
\begin{equation*}
y_a(t) = \Phi_h^{n,t} (\tilde{Y_a})
\end{equation*}
denoting the trajectory of a degree of freedom $\tilde{Y}_a$ and
\begin{equation*}
M_a = \tilde{M}_a \circ \left(\Phi^{n,t^{n-1}_+}_{h,\mathrm{approx}}\right)^{-1}
\end{equation*}
denoting the pushforward of the Lagrange shape functions $\tilde{M}_a$ to $\Omega^{t^{n-1}}$, we have
\begin{align*}
\varphi^{n,t}_{h,\mathrm{approx}}(X) 
&= \Phi^{n,t}_{h,\mathrm{approx}} \left( \left(\Phi^{n,t^{n-1}_+}_{h,\mathrm{approx}}\right)^{-1}(X) \right) \\
&= \sum_a \tilde{M}_a\left( \left(\Phi^{n,t^{n-1}_+}_{h,\mathrm{approx}}\right)^{-1}(X) \right) \Phi_h^{n,t}(\tilde{Y}_a) \\
&= \sum_a M_a(X) y_a(t).
\end{align*}
The corresponding material and spatial velocity fields are thus
\begin{equation*}
V^{n,t}_{h,\mathrm{approx}}(X) = \sum_a M_a(X) \dot{y}_a(t)
\end{equation*}
and
\begin{equation*}
v^{n,t}_{h,\mathrm{approx}}(x) = \sum_a m_a^t(x) \dot{y}_a(t),
\end{equation*}
respectively, with $m_a^t = \tilde{M}_a \circ
\left(\Phi_{h,\mathrm{approx}}^{n,t}\right)^{-1}$. 

Introducing approximations of the domain requires some extra care in
the imposition of boundary conditions. In the example problem here,
homogeneous Dirichlet boundary conditions are imposed on the boundary
of the approximate domain. Therefore, the order of the Lagrange shape
functions $\tilde{M}_a$ should be high enough to ensure that the
errors introduced by approximating $\Omega^t$ with
$\Phi_{h,\text{approx}}^{n,t}(\mathcal{S}_h^n)$ converge to zero at
least as quickly as the error in the original spatial discretization
as $h \rightarrow 0$.  It is well-known~\cite{Brenner1994} that if the
shape functions $\tilde{N}_a$ are themselves Lagrange shape functions,
then it suffices to use Lagrange shape functions $\tilde{M}_a$ of
equal or higher degree for the approximation of the domain geometry.
Elements of this type are referred to as \emph{isoparametric} or
\emph{superparametric} elements, depending upon whether the functions
$\tilde{M}_a$ have equal or higher degree, respectively, than the
functions $\tilde{N}_a$.

\subsection{Example: A Complete Algorithm}
\label{sec:exampl-one-compl}

We now present an algorithm that takes advantage of the SDIRK method
of \S \ref{sec:example:-runge-kutta} for time integration and of the
isoparametric representation of the domain of \S
\ref{sec:isop-appr-doma}. For concreteness, we consider the case in which the partial differential operator $\mathpzc{a}(u) = -\Delta_x u$, so that 
\begin{equation*}
a^t(u,w) = \int_{\Omega^t} \nabla_x u \cdot \nabla_x w \, dx,
\end{equation*}
In what follows, we denote matrices and vectors with uppercase and lowercase boldface letters, respectively.  As shorthand notation, we denote by 
\begin{equation*}
(u,w)_{K^t} = \int_{K^t} u(x)w(x) \, dx
\end{equation*}
the inner product of two functions $u$ and $w$ over an element $K^t = \Phi_h^{n,t}(K)$, $K \in \mathcal{S}_h^n$.
 The algorithm is labeled Algorithm~\ref{alg:algorithm1}.

\begin{algorithm}
\caption{Time integration using a universal mesh with an $s$-stage
  SDIRK method}
\label{alg:algorithm1}
\begin{algorithmic}[1]
\Require Initial condition $u^0 \in H^1_0(\Omega^0)$.
\For{$n=1,2,\dots, N$}
\vspace{2pt}\State Identify triangles in
$\mathcal{T}_{h,i}^{t^{n-1}}$, $i=0,1,2$. Set $\mathcal{S}_h^n =
\mathcal{T}_{h,0}^{t^{n-1}} \cup \mathcal{T}_{h,1}^{t^{n-1}} \cup
\mathcal{T}_{h,2}^{t^{n-1}}$. \label{identifytriangles}
\vspace{2pt}\State Compute $\Phi_h^{n,t^{n-1}_+}(\tilde{Y}_a)$ for every degree of freedom $\tilde{Y}_a \in \mathcal{D}(\mathcal{S}_h^n)$ using~(\ref{meshmap}). \label{phicomp1}
\vspace{2pt}\State Set $\mathcal{V}_h^n =
\mathrm{span}\{N_{a,\mathrm{approx}}\}_a$, where
$\{N_{a,\mathrm{approx}}\}_a$ are the shape
functions~(\ref{Naapprox}). \label{setspace}
\vspace{2pt}\State \begin{varwidth}[t]{\linewidth} Project $u_h^{\Delta t,n-1} \in (\varphi_h^{n-1,t^{n-1}})_* \mathcal{V}_h^{n-1}$ (or $u^0$ if $n=1$) onto $\mathcal{V}_h^n$ using a projector $p_h^n$. Denote by $\mathbf{u}_0$ the vector of coefficients in the expansion 
\begin{equation*}
p_h^n u_h^{\Delta t,n-1} = \sum_a (\mathbf{u}_0)_a N_{a,\mathrm{approx}}.
\end{equation*}
\end{varwidth}
\For{$i=1,2,\dots,s$} 
\State Compute $\Phi_h^{n,t_i}(\tilde{Y}_a)$ for every degree of freedom $\tilde{Y}_a \in \mathcal{D}(\mathcal{S}_h^n)$ using~(\ref{meshmap}). \label{phicomp2}
\State With $K^{t_i} = \Phi_{h,\mathrm{approx}}^{t_i}(K)$ for every $K \in \mathcal{S}_h^n$, assemble\label{assemble}
\begin{align*}
\mathbf{M}_{ab} &= \sum_K \left(n_{b,\mathrm{approx}}^{t_i},\, n_{a,\mathrm{approx}}^{t_i}\right)_{K^{t_i}} \\
\mathbf{B}_{ab} &= \sum_K \left(v_{h,\mathrm{approx}}^{n,t_i} \cdot \nabla_x n_{b,\mathrm{approx}}^{t_i},\, n_{a,\mathrm{approx}}^{t_i}\right)_{K^{t_i}} \\
\mathbf{K}_{ab} &= \sum_K \left(\nabla_x n_{b,\mathrm{approx}}^{t_i},\, \nabla_x n_{a,\mathrm{approx}}^{t_i}\right)_{K^{t_i}} \\
\mathbf{f}_{a} &= \sum_K \left(f(t_i),\, n_{a,\mathrm{approx}}^{t_i}\right)_{K^{t_i}}
\end{align*}
\State With $\mathbf{u}_* = \sum_{j=0}^{i-1} \beta_{ij} \mathbf{u}_j$ and $\Delta t^n = t^n-t^{n-1}$, define 
\begin{align*}
\mathbf{A} &= \mathbf{M} + \gamma \Delta t^n  (\mathbf{K}-\mathbf{B}) \\
\mathbf{b} &= \mathbf{M} \mathbf{u}_* + \mathbf{f}
\end{align*}
\State For every degree of freedom $\tilde{Y}_a \notin \mathrm{int}(\mathcal{D}(\mathcal{S}_h^n))$, set
\begin{align*}
\mathbf{A}_{ab} &= \delta_{ab} \\
\mathbf{b}_{a} &= 0.
\end{align*} \label{applyBC}
\State Solve $\mathbf{A} \mathbf{u}_i = \mathbf{b}$ for $\mathbf{u}_i$.
\EndFor
\State Set $u_h^{\Delta t,n}(x) = \sum_a (\mathbf{u}_s)_a n_{a,\mathrm{approx}}^{t^n}(x)$.
\EndFor
\State \Return{$u_h^{\Delta t,N}$}
\end{algorithmic}
\end{algorithm}

\begin{figure}[h]
  \centering
  \subfigure[Construction of $\Phi_{h,\text{approx}}^{1,0_+}(\tilde
  K)$ for $\tilde K\in
  \mathcal T_{h,2}^0$]{\includegraphics[width=0.49\linewidth]{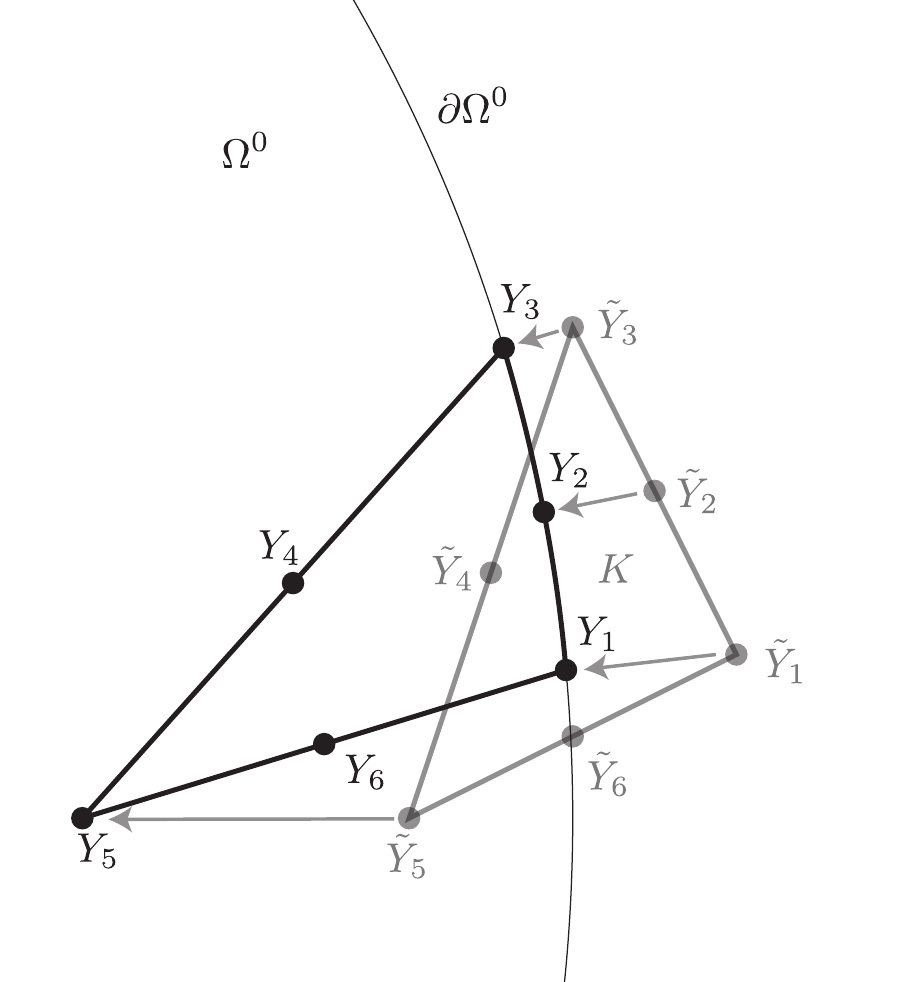}}
  \subfigure[Construction of $\Phi_{h,\text{approx}}^{1 ,t}(\tilde K)$
  for $\tilde K\in
  \mathcal T_{h,2}^0$, $0< t\le \Delta t$]{ \includegraphics[width=0.49\linewidth]{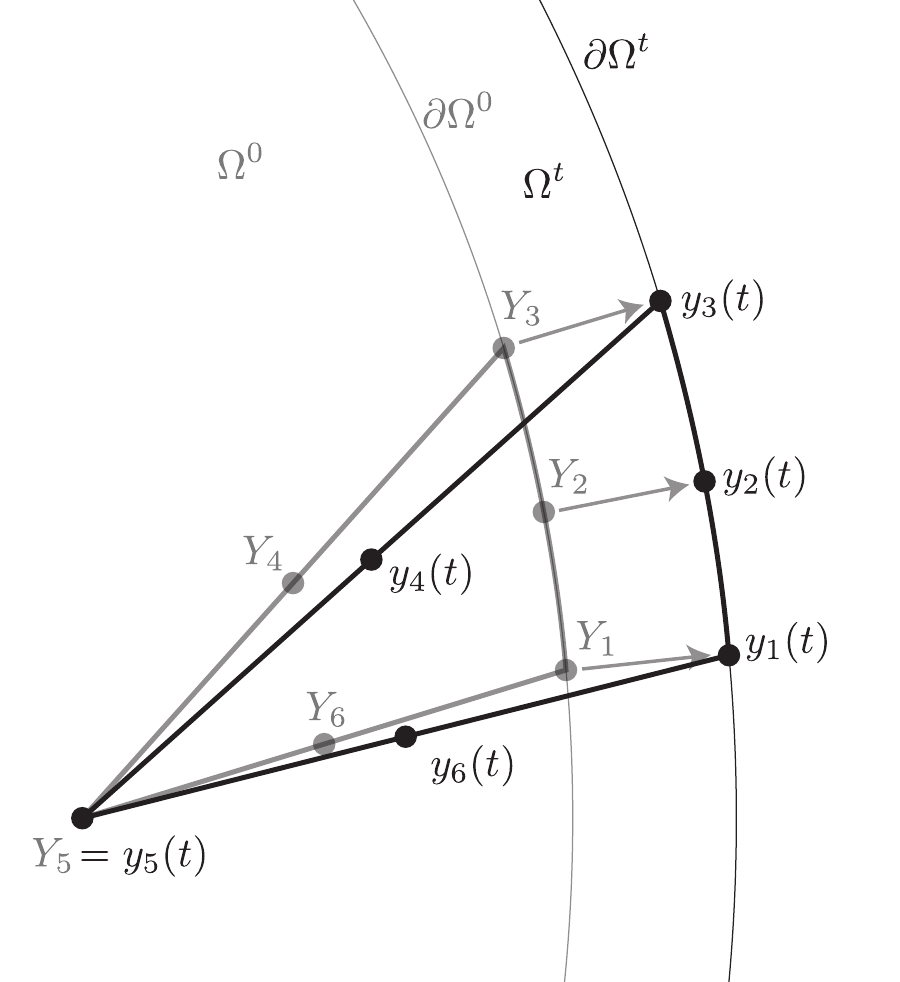}}
  \caption{Example of how the approximate evolving domain is accounted
    for in practice. See text in \S \ref{sec:exampl-one-compl} for the explanation.}
  \label{fig:exampletriangle}
\end{figure}

\paragraph{Implementation} We discuss some key steps of the algorithm
next, to show how the motion of the domain is accounted for in the
implementation of the algorithm, and how it affects the computation of
elemental quantities such as the mass matrix. For concreteness, in the
following it is useful to keep in mind a very simple example, such as
when the moving domain $\Omega^t$ is the circle centered at the origin
of radius $1+t$, for each small $t\ge 0$ (this is the geometry used to
draw Fig. \ref{fig:exampletriangle} later). Without loss of generality,
we discuss the case in which $n=1$, so that $t^{n-1}=0$. Finally, we will also
use the standard triangle $\hat K$, such as that with vertices
$(0,0)$, $(0,1)$, and $(1,0)$, which has traditionally been used in
finite element codes to perform quadrature.

In step \ref{identifytriangles} we identify triangles in
$\mathcal{T}_{h,i}^0$, for $i=0,1,2$, by labeling vertices of
triangles in the universal mesh according to whether they are inside
or outside $\Omega^0$. For example, Fig. \ref{fig:exampletriangle}
shows one triangle $\tilde K\in \mathcal T_{h,2}^{0}$. For this
example $\tilde K$ will be assumed to be quadratic and hence consists
of 6 nodes, with its nodes labeled by $\tilde Y_a$, $a=1,\ldots,6$.

In step \ref{phicomp1}, the positions $\{Y_a\}_a$ of these six nodes
in the mesh conforming to $\Omega_0$ are computed, and in general, of
all nodes in triangles intersecting $\partial \Omega_0$. This
computation involves computing the closest point projection for nodes
$\tilde Y_1, \tilde Y_2, \tilde Y_3$, and $\tilde Y_5$, moving the
first 3 nodes to their closest point projections, and moving $\tilde
Y_5$ along the normal to the boundary emanating from its closest point
projection, according to \eqref{relaxation}. Nodes $\tilde Y_4$ and
$\tilde Y_6$ are then mapped to the midpoints of segments $Y_3 Y_5$
and $Y_1Y_5$, respectively. These six nodes define the isoparametric
quadratic triangle $K=\Phi_{h,\text{approx}}^{1,0_+}(\tilde
K)$. Henceforth, {\it the construction of shape functions and quadrature
rules follow standard finite element procedures over isoparametric
elements}. For example, in this case,
$K=\hat \Psi(\hat K)$, where the isoparametric map is $\hat\Psi(\hat
X)=\sum_{a=1}^6 Y_a \hat N_a(\hat X)$ and $\{\hat
N_a\}_a$ denote the shape functions over $\hat K$. This is equivalent
to \eqref{Phi_approx}. 

In step \ref{setspace}, the shape functions over $K$ are
constructed. Because the map between $\hat K$ and $\tilde K$ is
affine, the shape function $\{N_{a,\text{approx}}\}_a$ can be
constructed over $\hat K$, namely, $N_{a,\text{approx}} (X
) =  \hat N_a(\hat\Psi^{-1}(X))$, for $X\in K$. This is, again, standard procedure for isoparametric elements.

As the boundary of the domain moves at each stage $i$ of the time
integration, the nodes $\{Y_a\}_a$ of triangle $K$ are deformed as
follows (step \ref{phicomp2}): Nodes $Y_1, Y_2$, $Y_3$ are mapped to their closest point
projections onto $\partial \Omega^{t_i}$, labeled $y_1, y_2$, and
$y_3$, respectively, node $Y_5$ remains where it is, so $y_5=Y_5$, and
nodes $Y_4$ and $Y_6$ are mapped to $y_4$ and $y_6$, the midpoints of
edges $y_3y_5$ and $y_1y_5$, respectively, see Fig. \ref{fig:exampletriangle}b. Shape functions over
triangle $K^{t_i}$ are formed in precisely the same way as those for
triangle $K$, in this case with nodal positions $\{y_a\}_a$. 

To assemble the system needed to solve  \eqref{dirk} at each stage of
the time integration (step \ref{assemble}) it is useful to notice that the elemental mass
matrix for each element is computed as
\begin{equation}
  \label{eq:4}
  \begin{aligned}
  \mathbf{M}_{ab}^K & = \int_{K^{t_i}} n^{t_i}_{a,\text{approx}}
  n^{t_i}_{b,\text{approx}}\;dx\\
  & = \int_{K} N_{a,\text{approx}} N_{b,\text{approx}} |\nabla_X
  \varphi_{h,\text{approx}}^{1,t_{i}}|\;dX\\
& = \int_{\tilde K} \tilde N_a \tilde N_b |\nabla_{\tilde X}
\Phi_{h,\text{approx}}^{1,t_i}|\; d\tilde X\\
  & = \int_{\hat K} \hat N_a\hat N_b |\nabla_{\hat X} \hat \Psi|\;d\hat X,
  \end{aligned}
\end{equation}
and similarly for the elemental contributions to the other terms of
\eqref{dirk}. Consequently, quadrature could be performed on any of
the triangles $K$, $\tilde K$, $\hat K$, or $K^{t_i}$, but for
convenience and following standard practice, we do it over the
standard element $\hat K$. 

Notice then that, in order to perform the quadrature over $\hat K$, it
is convenient to build the deformed mesh at time $t_i$, since it makes
the construction of $\hat \Psi$ straightforward. So, in this
time-integration scheme the deformed mesh is built $s$ times in a time
step.

An important practical matter we wish to highlight is the simplicity
of the data structures needed to implement our method.  In particular,
\emph{the connectivity of the universal mesh never changes} during
deformation -- only the nodal positions change.  As a consequence, the
sizes and sparsity structures of various discrete quantities (the
solution vector $\mathbf{u}$, the mass matrix $\mathbf{M}$, the
stiffness matrix $\mathbf{K}$, the convection matrix $\mathbf{B}$, and the forcing vector $\mathbf{f}$)
can be held fixed, even though differing subsets of degrees of freedom
may participate in the discrete equations at different intervals $(t^{n-1},t^n]$.  This
can be accomplished by simply imposing ``homogeneous Dirichlet boundary
conditions'' on the solution at degrees of freedom not belonging to the
subtriangulation $\mathcal{S}_h^n$.  In practice, this amounts to
replacing the corresponding rows of a particular matrix $\mathbf{A}$
with rows whose only nonzero entries are 1 on the diagonal, and
setting to zero the corresponding entries of a vector (see
step~\ref{applyBC} of Algorithm~\ref{alg:algorithm1}).  Note that
$\mathbf{A}$ is automatically asymmetric at the outset, so any
concerns of breaking symmetry via row replacement are irrelevant.

\subsection{Exact vs. Approximate Map: Cost Considerations}

The computational cost of evaluating the map $\Phi_h^{n,t}$ or its
approximant $\Phi_{h,\text{approx}}^{n,t}$ is dominated by the cost of
evaluating closest-point projections onto $\partial\Omega^t$.  In our
numerical experiments (which used $\Phi_{h,\text{approx}}^{n,t}$),
these calculations accounted for little more than $5-10\%$ of the
total run time of a typical simulation. 

Note that implementations that employ the exact map $\Phi_{h}^{n,t}$ require evaluations of the closest point projection and its gradient at \emph{quadrature points} in triangles $K \in \mathcal{T}_{h,1}^{t^{n-1}} \cup \mathcal{T}_{h,2}^{t^{n-1}}$, whereas implementations that employ the approximate map $\Phi_{h,\text{approx}}^{n,t}$ require evaluations only of the closest point projection (not its gradient) on those triangles' \emph{degrees of freedom}.  A counting argument reveals that the computational savings that accompany the use of $\Phi_{h,\text{approx}}^{n,t}$ over $\Phi_h^{n,t}$ are significant: For a polynomial interpolant $\Phi_{h,\mathrm{approx}}^{n,t}$ constructed from Lagrange elements of a fixed polynomial degree, it is not difficult to show that the use of $\Phi_{h,\text{approx}}^{n,t}$ over $\Phi_h^{n,t}$ reduces the computational cost (measured by number of closest point projection evaluations) by factors of 9, 9, and 5.2, respectively, for affine, quadratic, and cubic Lagrange elements, assuming the use of a quadrature rule that exactly computes entries of the elemental mass matrix on straight triangles.

\subsection{Error estimate for a universal mesh}
\label{sec:error-estim-univ}
We conclude this section by applying the error estimate~(\ref{error_general}) to the case in which the maps $\Phi_h^{n,t}$ are constructed from a universal mesh according to the algorithm in Section~\ref{sec:construction}.

It is shown in our forthcoming paper~\cite{Gawlik2012b} that for $u$ and $\Omega$ sufficiently regular, there exist
constants $B$, $c$, $C$, $\bar{C}_1(u,T)$, and $\bar{C}_2(u,T)$,
independent of $h$, $N$, and $\Delta t$, such that
\begin{align}
C_1(u,v_h,T) &\le \ell_{h,r} h^{-1/2}\bar{C}_1(u,T) \label{Border} \\
C_2(u,v_h,T) &\le \bar{C}_2(u,T) \label{Corder} \\
N &\le B T \Delta t^{-1} \label{Norder} \\
|\mathcal{R}_h^n| &\le C h, \;\;\; n=1,2,\dots,N\label{rorder}
\end{align}
as long as 
\begin{equation} \label{restriction}
c \, \max_{1 \le n \le N} (t^n - t^{n-1}) \le h
\end{equation}
and
\begin{equation*}
\max_{1 \le n \le N} (t^n-t^{n-1}) \le B \min_{1 \le n \le N} (t^n-t^{n-1}),
\end{equation*}
where $\ell_{h,r}$ is given by~(\ref{logh}).  We remind the reader that $\Delta t \le \min_{1 \le n \le N} (t^n-t^{n-1})$.

Relations~(\ref{Border}-\ref{rorder}) imply that the numerical solution obtained from the
proposed strategy has the accuracy stated in \eqref{eq:5}, that is, 
\begin{equation}
\|u_h^{\Delta t,N} - u^N\|_{0,2,\Omega^T} \le  C(u,T)\left(\ell_{h,r} h^{r-1/2} + \Delta t^q + \ell_{h,r} h^{r+1/2} \Delta t^{-1} \right). \label{error_uni}
\end{equation}
When $\Delta t$ and $h$ scale proportionately, this convergence rate is suboptimal by half an order (up to a logarithmic factor if $r=2$).

\paragraph{Time step restriction}

Notice that the text above included the restriction~(\ref{restriction}), which implies a restriction on the time step of the form $c \Delta t \le h$.   The necessity of such a restriction is made clear by noting that the mesh motion defined by $\varphi_h^{n,t}$, $t \in (t^{n-1},t^n]$, leaves all elements stationary except those with an edge on the moving boundary.  Imposing~(\ref{restriction}) with a suitable choice of $c$ ensures that the image under $\varphi_h^{n,t}$ of each such element has an aspect ratio that is bounded above and below uniformly in $h$ for all times $t \in (t^{n-1},t^n]$.  In particular, it ensures that no element collapses to a set of nonpositive measure at any time $t \in (t^{n-1},t^n]$.  Note that this restriction is intrinsic to the mesh motion strategy; it is a restriction that must be imposed in addition to any time step restriction needed to ensure stability of the particular time integrator chosen.

\paragraph{Explanation of estimate} Let us briefly describe how the
dependencies~(\ref{Border}-\ref{rorder}) arise.  To begin, consider
the regularity of the velocity field $v_h$.  The proposed strategy
employs a velocity field that is of order unity on the boundary of
$\Omega^t$ and decays to zero over a strip of neighboring elements,
i.e. a strip of width $O(h)$.  It follows that the velocity field
itself is everywhere of order unity, but its spatial gradient is of
order $h^{-1}$.  Moreover, the support of $v_h^{n,t}$ and its
derivatives (the strip of neighboring elements) has measure $O(h)$.
From these observations, a simple calculation reveals that the $L^2$
norm of $\nabla_x v_h^{n,t}$ is of order $h^{-1/2}$.  It is this fact
that contributes to the dependence of $C_1(u,v_h,T)$ on $h^{-1/2}$.
The factor $\ell_{h,r}$ arises because of the approximation properties
in the $L^\infty$-norm of piecewise linear finite element
spaces.

Consider next the relation~(\ref{rorder}).  The quantity $|\mathcal{R}_h^n|$ measures the discrepancy between two triangulations of the same domain $\Omega^{t^{n-1}}$, namely $\Phi_h^{n-1,t^{n-1}}(\mathcal{S}_h^{n-1})$ and $\Phi_h^{n,t^{n-1}_+}(\mathcal{S}_h^n)$.  By~(\ref{relaxation}) and~(\ref{meshmap}), these triangulations differ only in a neighborhood of $\partial\Omega^{t^{n-1}}$ that has measure $O(h)$.  A consequence of this fact is that the contribution to the error associated with projecting the solution onto a new finite element space at each temporal node $t^n$ is of order $Nh^{r+1/2} \sim h^{r+1/2} \Delta t^{-1}$, rather than of the order $Nh^r \sim h^r \Delta t^{-1}$ that one might expect if the triangulations had differed over the entire domain.

Notice that the terms $\ell_{h,r} h^{r-1/2}$ and $\ell_{h,r} h^{r+1/2} \Delta t^{-1}$ in~(\ref{error_uni}) are balanced when $\Delta t \sim h$.  Roughly speaking, the asymptotically unbounded gradient of $v_h^{n,t}$ introduces a half-order reduction in the spatial discretization error (ordinarily $h^r$), but the small support of $v_h^{n,t}$ mitigates the reduction of order introduced by the repeated projections (ordinarily a full order) to half an order.

\paragraph{Optimal balance}
Given the half-order reduction in error associated with the use of a velocity field $v_h^{n,t}$ whose support has measure $O(h)$, it is tempting to consider the possibility of employing mesh motions with more broadly supported velocity fields.  An immediate consequence of such a decision, however, is an increase in the error $h^r \ell_{h,r} |\mathcal{R}_h^n|^{1/2}$ associated with changing finite element spaces.  Indeed, if the triangulations $\Phi_h^{n-1,t^{n-1}}(\mathcal{S}_h^{n-1})$ and $\Phi_h^{n,t^{n-1}_+}(\mathcal{S}_h^n)$ differ over a region of measure $O(1)$, then $|\mathcal{R}_h^n|$ is of order 1 rather than of order $h$.  The resulting total error estimate is suboptimal by \emph{one} order rather than half an order when $N = O(\Delta t^{-1})$.  For this reason, the strategy as it has been presented at the outset can be said to provide an optimal balance between competing sources of error.

Note, however, that if the finite element spaces are changed less frequently (i.e. $N = O(1)$), then optimal order accuracy is, in principle, obtainable via the use of more broadly supported velocity fields.  This is precisely what is accomplished by conventional ALE schemes, which adopt global mesh motion strategies (in which all nodes of the mesh participate) and remesh occasionally (at temporal nodes $t^n$ whose separation in time is of order unity).  The price to be paid for such a decision, of course, is a reduction in the efficiency and robustness of the mesh motion strategy.  Conventional ALE mesh motions commonly require the solution of systems of equations (such as those of linear elasticity) for the positions of mesh nodes~\cite{Donea1983,Farhat1998,Johnson1994,Helenbrook2003}; in contrast, the nodal motions in our method are independent and explicitly defined, rendering the mesh motion strategy low-cost and easily parallelizable.  Second, our mesh motion strategy is robust in the sense that it enjoys provable bounds on the quality of the deformed mesh under suitable constraints on the time step and mesh spacing~\cite{Rangarajan2012a,Gawlik2012b}. 

\section{Numerical Examples} \label{sec:numerical_examples}	

In this section, we apply the proposed method to a modification of a classical
moving-boundary problem: Stefan's problem.  In our modification, the
evolution of the boundary is imposed through the exact solution,
instead of being computed. Our aim in this example is
is to illustrate the convergence rate of the method with respect to
the mesh spacing $h$ and time step $\Delta t$.

We begin by demonstrating, using a one-dimensional numerical test,
that the bound~(\ref{error_uni}) is sharp.  That is, the order of
accuracy of the method is suboptimal by half an order in the $L^2$
norm when $\Delta t$ and $h$ scale proportionately.  We observe,
however, that the suboptimal rate is difficult to detect from an
inspection of the total error $\|u_h^{N,\Delta
  t}-u^N\|_{0,2,\Omega^T}$, since the terms of suboptimal order
contributing to the total error are dominated by terms of optimal
order (for practical values of the mesh spacing $h$).  We follow with
a convergence test in two-dimensions, where, for the reason just
described, optimal rates are observed for the total error.

\subsection{The (Modified) One-Dimensional Stefan Problem with Prescribed Boundary Evolution}
\label{sec:one-dimens-stef}
Consider the following instance of the one-dimensional Stefan problem:  Find $u(x,t)$ and $s(t)$ such that
\begin{subequations}
\begin{alignat}{2}
\frac{\partial u}{\partial t} &= \frac{\partial^2 u}{\partial x^2} , &0 <  x < s(t), \; t \ge 1 \label{1dstefana} \\
\frac{ds}{dt} &= -\frac{\partial u}{\partial x}, &x=s(t) \label{1dstefanb} \\
u(0,t) &= e^t-1, &t \ge 1 \label{1dstefanc} \\
u(s(t),t) &= 0, &t \ge 1 \label{1dstefand} \\
u(x,1) &= e^{1-x}-1, &0 \le x \le 1 \label{1dstefane} \\
s(1) &= 1. & \label{1dstefanf}
\end{alignat} \label{1dstefan}
\end{subequations}
The exact solution is
\begin{align*}
u(x,t) &= e^{t-x}-1 \\
s(t) &= t.
\end{align*}
In this case, we treat
the boundary evolution as prescribed by supplying the exact evolution
$s(t)$, instead of solving for it by integrating
(\ref{1dstefanb}). 

We  computed the numerical solution $u_h^{\Delta t,N}$ using
a finite element space made of continuous elementwise-affine functions on a sequence of
uniform meshes with spacing $h=2^{-k}h_0$, $k=0,1,2,3$ over the time
interval $[1,T]$ with $h_0 = 1/4$ and $T=1+10^{-6}$ (the short time interval was
chosen, on the basis of numerical experiments, in order to detect the
suboptimal rate predicted by the theory).  The restriction of the algorithm to
a single spatial dimension is that specified in Algorithm
\ref{algorithm1d}, and is 
complemented with the choice $p^n_h = p^n_{h,L^2}$ for
the projection, relaxation parameters $\delta = 0.3$ and $R=3$, and the singly diagonally implicit Runge-Kutta (SDIRK)
scheme of order 2 given in Table~\ref{tab:SDIRK2} with a time step $\Delta t =
10^{-6} h/h_0$ for time-integration. 

Table~\ref{tab:L2stefan1d} presents the convergence of the method
measured at time $t=T$ in $L^2(\Omega^T)$.  The third column of the table suggests
that the total error $\|u_h^{\Delta t,N}-u^N\|_{0,2,\Omega^T}$
converges at an optimal rate $O(h^2)$.  However, columns 1 and 2
reveal that a piece of the error, namely the discrepancy between the
numerical solution $u_h^{\Delta t,N}$ and the nodal interpolant
$i_h^{t^N} u^N$ of the exact solution, decays at a suboptimal rate
$O(h^{3/2})$.  Since standard estimates from the theory of
interpolation give $\|u^N-i_h^{t^N} u^N\|_{0,2,\Omega^T} = O(h^2)$, it
follows from the inequality
\begin{equation*}
\|u_h^{\Delta t,N} - i_h^{t^N} u^N\|_{0,2,\Omega^T} \le \|u_h^{\Delta t,N} - u^N\|_{0,2,\Omega^T} + \|u^N - i_h^{t^N} u^N\|_{0,2,\Omega^T} 
\end{equation*}
that $\|u_h^{\Delta t,N} - u^N\|_{0,2,\Omega^T}$ must be decaying no faster than $O(h^{3/2})$.  However, the contribution to the error supplied by $u_h^{\Delta t,N} - i_h^{t^N} u^N$ is several orders of magnitude smaller than the remaining contribution, $i_h^{t^N} u^N - u^N$, explaining the apparent optimal rate observed for the total error.

\begin{table}
\centering
\caption{Convergence rates in the $L^2$-norm on $\Omega^T$ for the
  solution to the (modified) one-dimensional Stefan problem using a finite
  element space made of continuous
  elementwise-affine functions
  with a second-order implicit Runge-Kutta time integrator, see \S
  \ref{sec:one-dimens-stef}. Differences between the exact solution
  $u^N$, the numerical approximation $u_h^{\Delta t,N}$, and the nodal
interpolant of the exact solution $i_h^{t^N}u^N$ are shown in each
column. These values are used in \S \ref{sec:one-dimens-stef} to
illustrate that the expected theoretical convergence rate of $h^{3/2}$
is observed. Nevertheless, the slowly converging part is so small,
that the apparent convergence rate is $h^2$, as the third
column shows.}
\label{tab:L2stefan1d}
\newcolumntype{R}{>{\raggedleft\arraybackslash}X}%
\begin{tabularx}{470pt}{c|cc|cc|cc}
\hline\noalign{\smallskip}
$h_0/h$ & \mbox{\footnotesize $\|u_h^{\Delta t,N}-i_h^{t^N}u^N\|_{0,2,\Omega^T}$} &  Order & \mbox{\footnotesize $\|i_h^{t^N}u^N-u^N\|_{0,2,\Omega^T}$} & Order & \mbox{\footnotesize $\|u_h^{\Delta t,N}-u^N\|_{0,2,\Omega^T}$} & Order  \\
\noalign{\smallskip}\hline\noalign{\smallskip}
1 & 1.2e-11 & \multicolumn{1}{c|}{\hspace{0pt}-} & 4.4e-05   &  \multicolumn{1}{c|}{\hspace{0pt}-} & 4.4e-05  & \multicolumn{1}{c}{\hspace{0pt}-} \\
2 & 4.6e-12 &      1.42 & 1.1e-05 &      2.04 & 1.1e-05 &      2.04 \\
4 & 1.6e-12 &      1.49 & 2.6e-06 &      2.02 & 2.6e-06 &      2.02 \\
8 & 5.5e-13 &      1.56 & 6.5e-07 &      2.01 & 6.5e-07 &      2.01 \\
\noalign{\smallskip}\hline
\end{tabularx}
\end{table}

\subsection{The (Modified) Two-Dimensional Stefan Problem with Prescribed Boundary Evolution}

\begin{figure}
  \centering
  \includegraphics[width=0.65\textwidth]{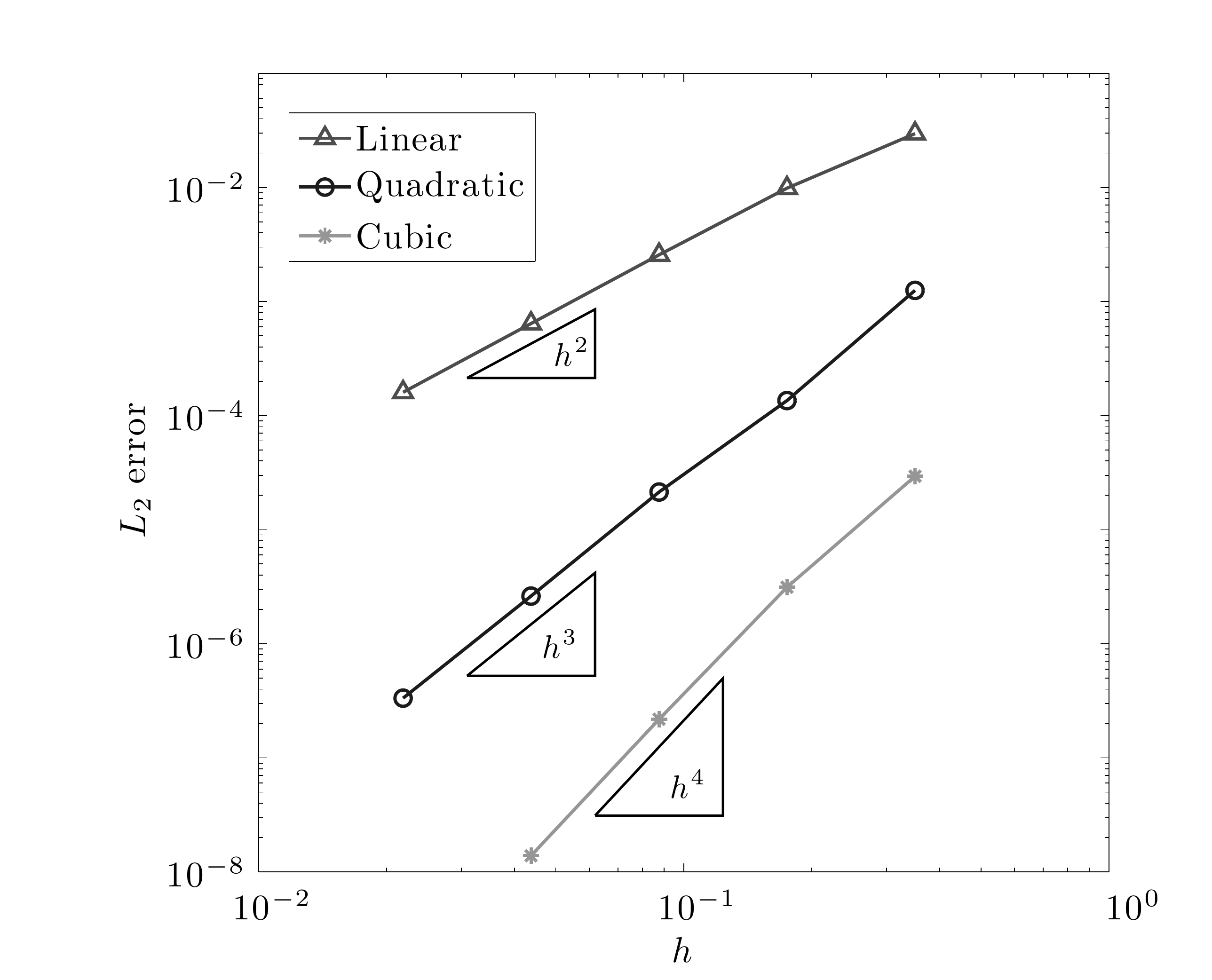}              
  \caption{(a) $L^2$-error $\|u_h^{\Delta t,N}-u^N\|_{0,2,\Omega^T}$
    at a fixed final time as a function of the mesh spacing $h$ for
    the (modified) two-dimensional Stefan problem~(\ref{2dstefan}) with
    prescribed boundary evolution.  The problem was solved using
    linear, quadratic, and cubic elements together with nodal
    interpolation as the projection operator, and second-, third-, and
    fourth-order implicit Runge-Kutta schemes, respectively, as time
    integrators, with $h \propto \Delta t$.}
  \label{fig:L2stefan}
\end{figure}

\begin{table}
\centering
\caption{Convergence rates in the $L^2$-norm on $\Omega^T$ for the
  solution to the (modified) two-dimensional Stefan problem \eqref{2dstefan} using linear, quadratic, and cubic elements together  with nodal
  interpolation as the projection operator, and second-, third-, and
  fourth-order implicit Runge-Kutta schemes, respectively, as time integrators. See Fig.~\ref{fig:L2stefan} for a graphical depiction of the same results.}
\label{tab:L2stefan}
\newcolumntype{R}{>{\raggedleft\arraybackslash}X}%
\begin{tabularx}{310pt}{c|cc|cc|cc}
\hline\noalign{\smallskip}
 &  \multicolumn{2}{c|}{Linear} & \multicolumn{2}{c|}{Quadratic} & \multicolumn{2}{c}{Cubic}  \\
$h_0/h$ & \multicolumn{1}{c}{Error} & \multicolumn{1}{c|}{Order} & \multicolumn{1}{c}{Error} & \multicolumn{1}{c|}{Order} & \multicolumn{1}{c}{Error} & \multicolumn{1}{c}{Order} \\
\noalign{\smallskip}\hline\noalign{\smallskip}
  1  &     3.0e-02  &  \multicolumn{1}{c|}{-}  &     1.3e-03  &  \multicolumn{1}{c|}{-}  &     2.9e-05  &  \multicolumn{1}{c}{-} \\ 
  2  &     9.8e-03  &        1.59  &     1.4e-04  &        3.21  &     3.1e-06  &        3.24 \\ 
  4  &     2.6e-03  &        1.94  &     2.1e-05  &        2.66  &     2.2e-07  &        3.84 \\ 
  8  &     6.4e-04  &        2.00  &     2.6e-06  &        3.03  &     1.4e-08  &        3.97 \\ 
 16  &     1.6e-04  &        2.00  &     3.3e-07  &        2.97  &         -  &         - \\ 
\noalign{\smallskip}\hline
\end{tabularx}
\end{table}

We consider now the following instance of the two-dimensional,
cylindrically symmetric Stefan problem with a circular boundary of
radius $\rho(t)$ centered at the origin. Find the scalar functions
$u(x,t)$ and $\rho(t)$ such that for all times $t\in[0,T]$,
\begin{subequations}
  \begin{alignat}{2}
    \frac{\partial u}{\partial t} - \Delta_x u & = f,\qquad\qquad &0 \le |x| < \rho(t) \label{2dstefana} \\
    \frac{d\rho}{dt} &= -\frac{\partial u}{\partial n}, &|x| = \rho(t)  \label{2dstefanb} \\
    u(x,t)& = 0, &|x| = \rho(t) \label{2dstefanc} \\
    u(x,0)& = J_0(r_0 |x|), &  \label{2dstefand} \\
    \rho(0) &= 1, & \label{2dstefane}
  \end{alignat} \label{2dstefan}
\end{subequations}
where $J_0$ is the zeroth-order Bessel function of the first kind, $r_0$ is the smallest positive root of $J_0$, and 
\begin{align*}
f(x,t) &= \frac{\alpha r_0^3 \beta(t)^2 |x|}{2\sigma(t)^3}  J_0'\left(\frac{r_0 |x|}{\sigma(t)} \right) \\
\sigma(t) &= \exp\left(\frac{\alpha(\beta(t)-1)}{2}\right) \\
\beta(t) &= \frac{1}{\alpha} \mathrm{Ei}^{-1}\left(\mathrm{Ei}(\alpha)-r_0^2 t e^\alpha\right) \\
\alpha &= \frac{2J_0'(r_0)}{r_0}.
\end{align*}
Here, $\mathrm{Ei}(z) = -\int_{-z}^\infty \frac{e^{-\zeta}}{\zeta} d\zeta$, the exponential integral.  The exact solution is
\begin{align*}
u(x,t) &= \beta(t) J_0\left(\frac{r_0 |x|}{\sigma(t)} \right) \\
\rho(t) &= \sigma(t).
\end{align*}

In our implementation, we treat the boundary evolution as prescribed
by supplying the exact evolution of the moving domain's radius
$\rho(t)$, instead of solving for it.  To study the convergence of the
method, the problem was solved using finite element spaces of
continuous functions that are affine, quadratic, and cubic over each
element (linear, quadratic, and
cubic Lagrange elements) together with nodal interpolation as the projection
operator, relaxation parameters $\delta=0.8$ and $R=3$, and singly diagonally implicit Runge-Kutta (SDIRK) schemes
of orders 2, 3, and 4, respectively, as the time integrators (see the
coefficients in Tables~\ref{tab:SDIRK2}-\ref{tab:SDIRK4}). The solution was
computed on a uniform mesh of equilateral triangles with a lowest resolution mesh spacing of $h_0=0.35$ and a time step $\Delta t=Th/h_0$, up to a final time $T=0.005$. 

Fig.~\ref{fig:L2stefan} displays the $L^2$-error of the numerical
solution as a function of the mesh spacing $h$ at $t=T$.  Optimal
convergence orders of 2, 3, and 4 are observed for the three schemes,
in agreement with the observations made in the one-dimensional test
case. Table \ref{tab:L2stefan} shows the same results.

\begin{figure}%
\centering
\subfigure[{$t=0$}]{%
\centering
\includegraphics[width=0.4\textwidth,trim=0.8in 0.85in 0.8in 2.7in,clip=true]{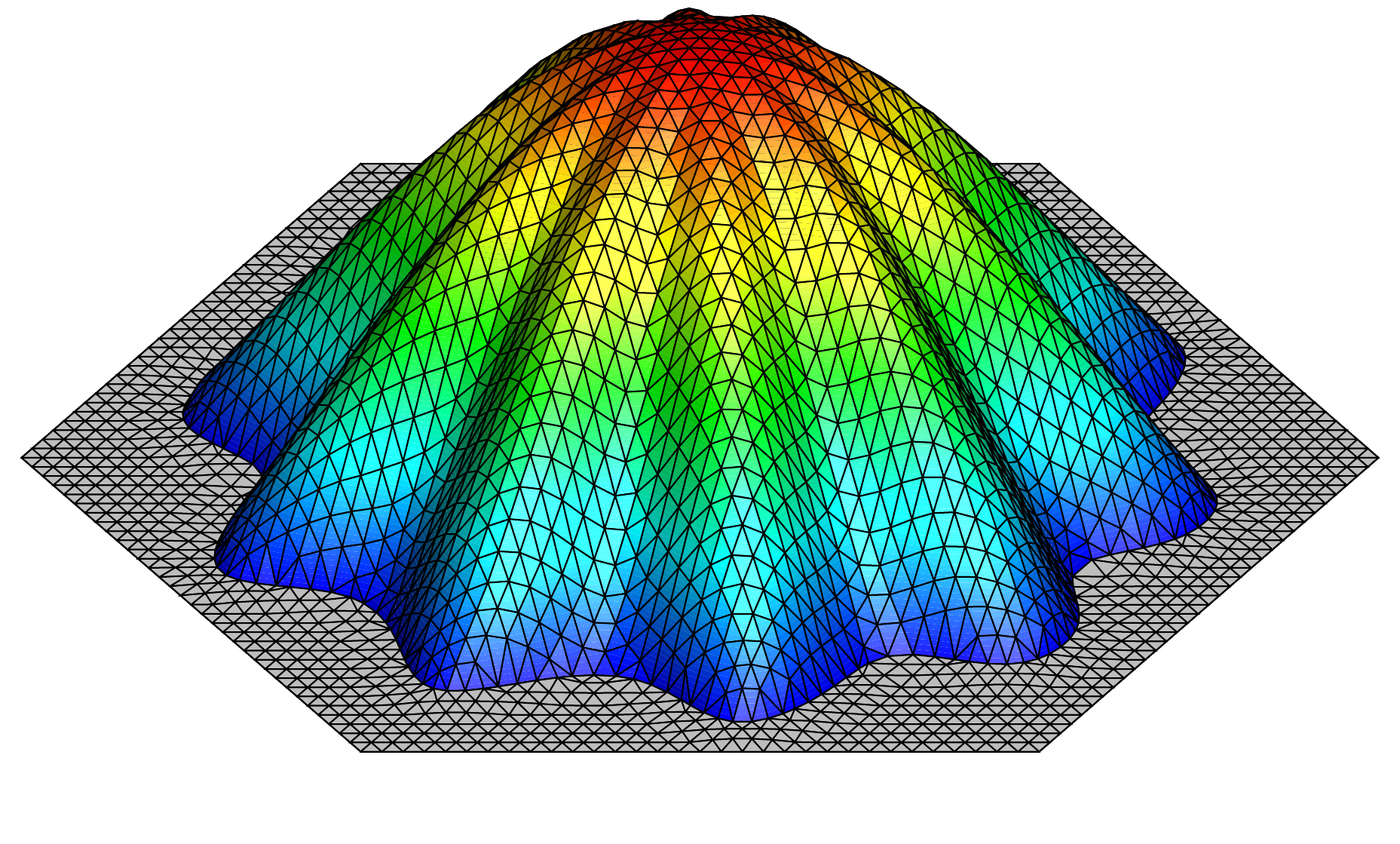}%
\label{fig:stefan000}%
}%
\subfigure[{$t=0.02$}]{%
\centering
\includegraphics[width=0.4\textwidth,trim=0.8in 0.85in 0.8in 2.7in,clip=true]{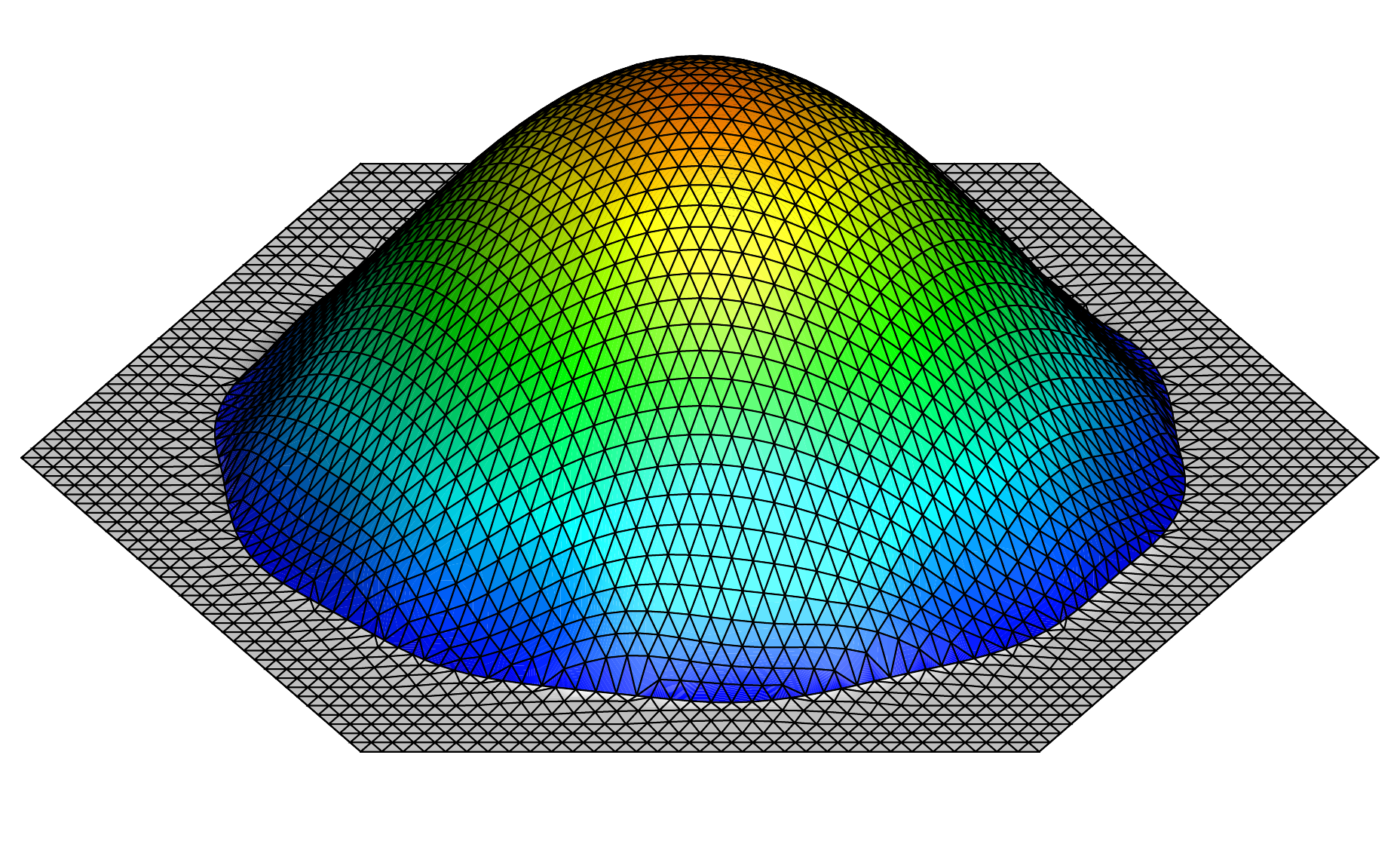}%
\label{fig:stefan031}%
}%
\\
\subfigure[{$t=0.04$}]{%
\centering
\includegraphics[width=0.4\textwidth,trim=0.8in 0.85in 0.8in 3.1in,clip=true]{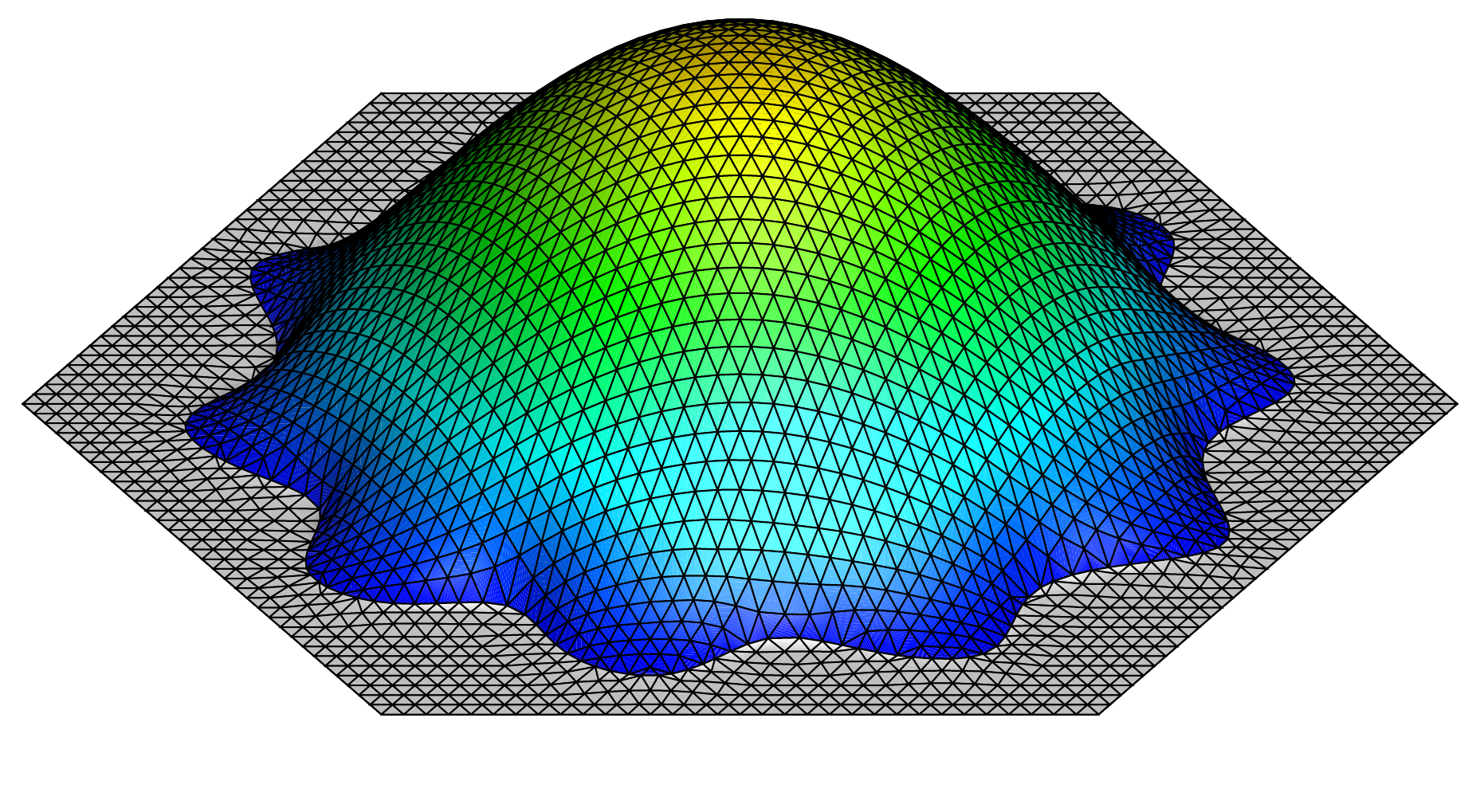}%
\label{fig:stefan063}%
}%
\subfigure[{$t=0.06$}]{%
\centering
\includegraphics[width=0.4\textwidth,trim=0.8in 0.85in 0.8in 3.1in,clip=true]{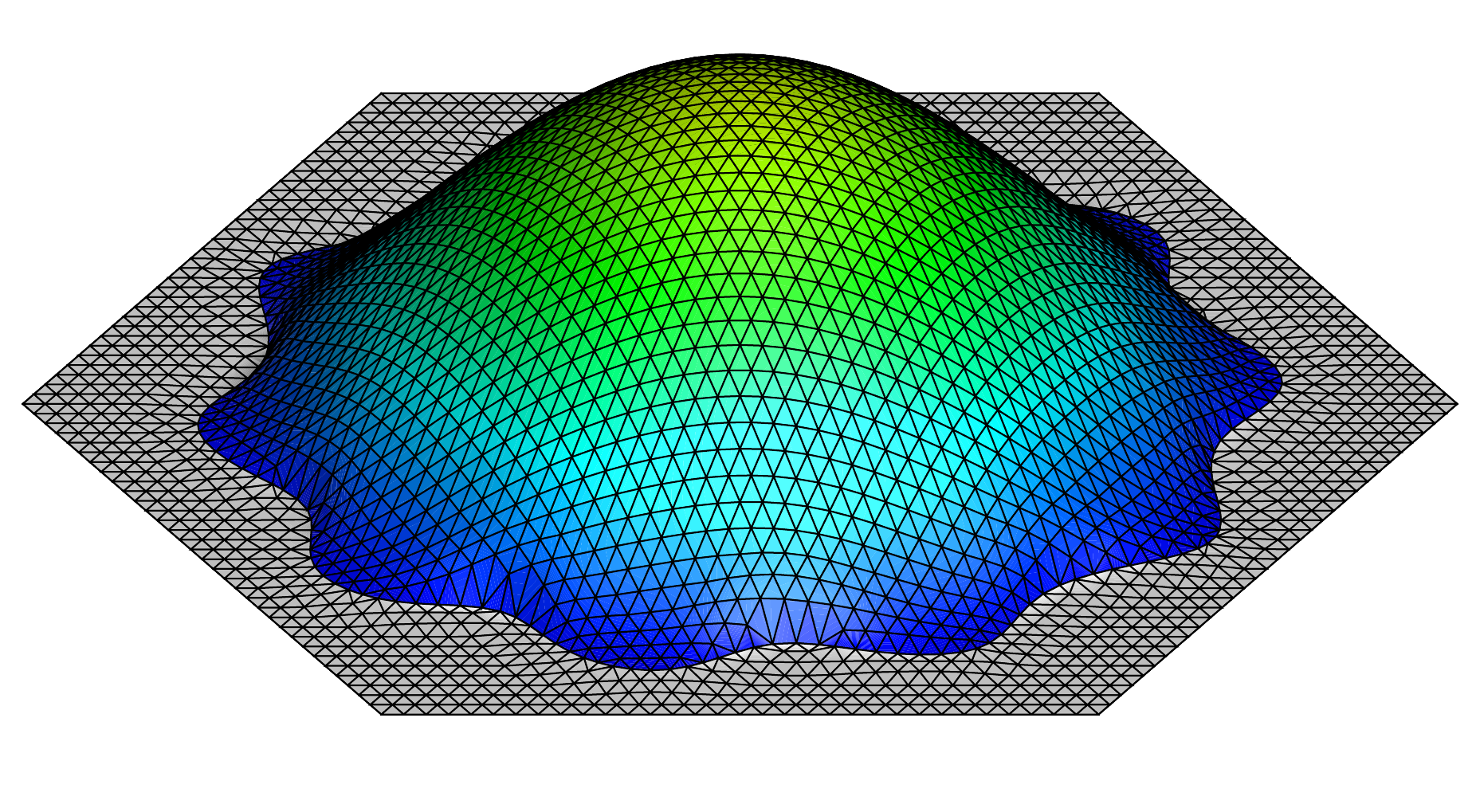}%
\label{fig:stefan095}%
}%
\caption{Solution to a prescribed-boundary variant of the Stefan problem in which the moving boundary is a sinusoidal perturbation of the unit circle.}
\label{fig:stefan_snapshots}
\end{figure}

\begin{figure}%
\centering
\subfigure[Universal mesh]{%
\centering
\includegraphics[width=0.45\textwidth,trim=0.8in 0.85in 0.8in 1.3in,clip=true]{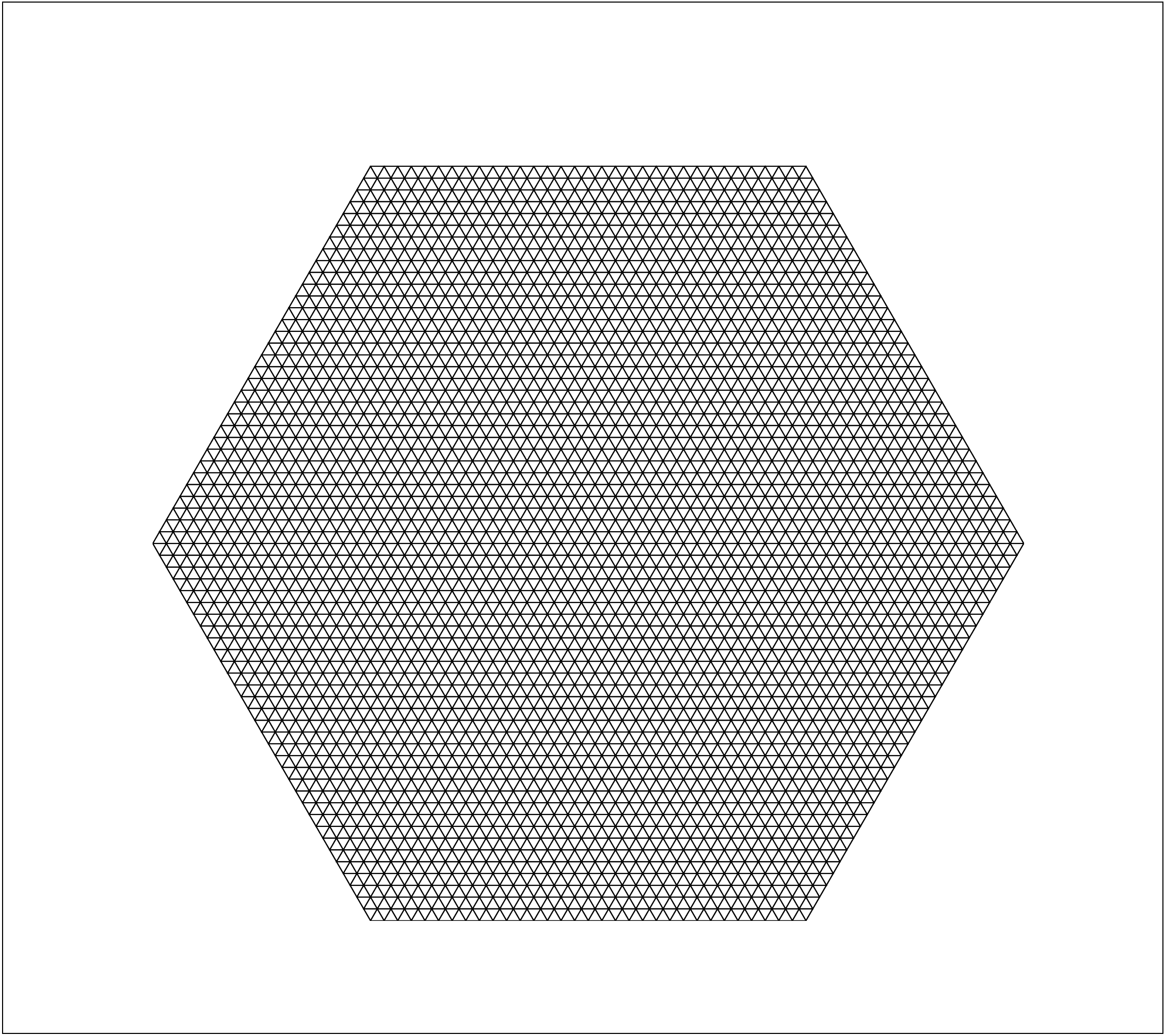}%
\label{fig:stefan_bgmesh}%
}%
\subfigure[{$t=0.06$}]{%
\centering
\includegraphics[width=0.45\textwidth,trim=0.8in 0.85in 0.8in 1.3in,clip=true]{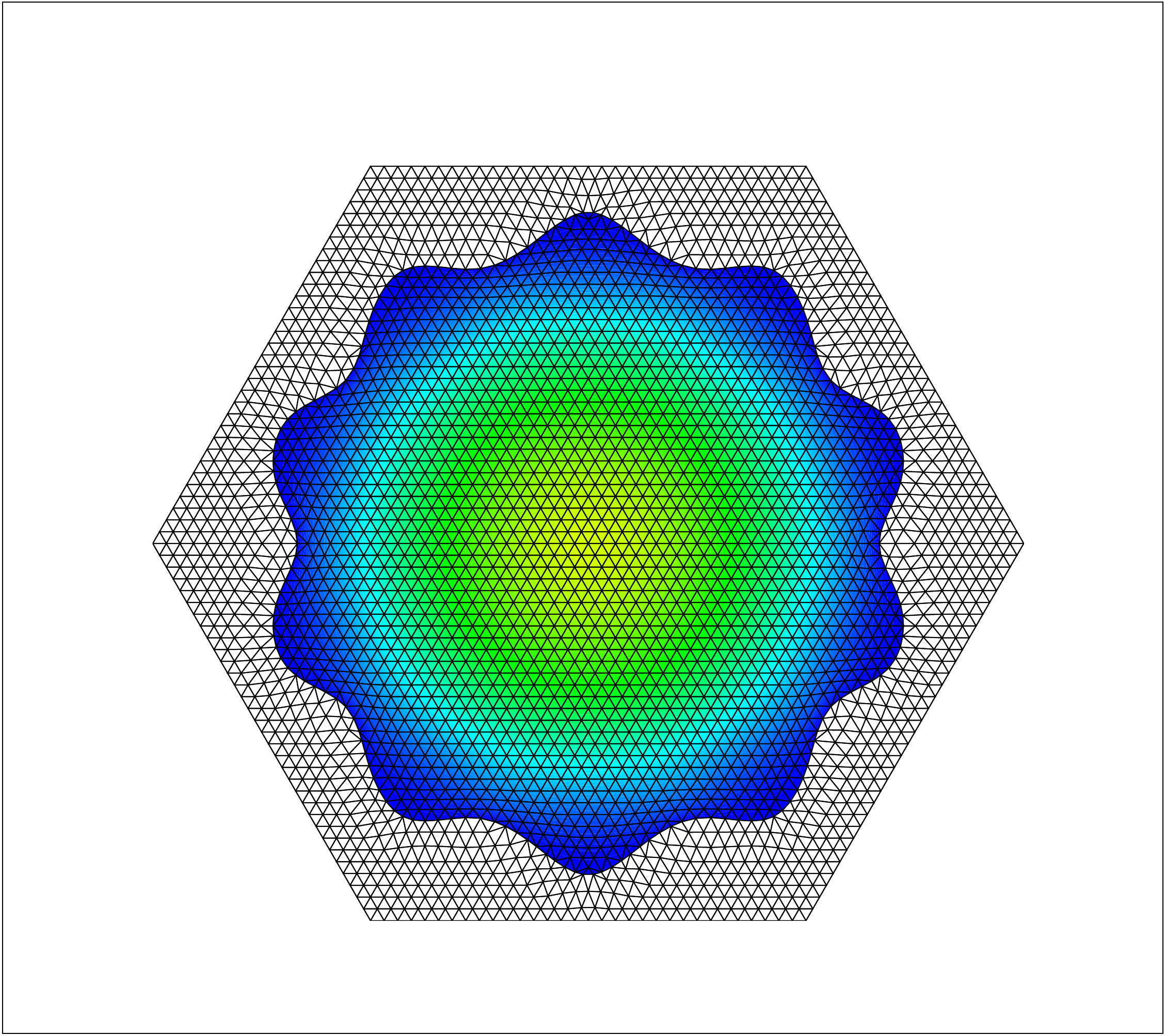}%
\label{fig:stefan2D095}%
}%
\caption{(a) Universal mesh adopted during the simulation depicted in Fig.~\ref{fig:stefan_snapshots}, and (b) its image under the universal mesh map at $t=0.06$, superposed with the contours of the solution.}
\label{fig:stefan_mesh}
\end{figure}

To illustrate the method on a second, more interesting example, we solved the partial differential equation~(\ref{2dstefana}) with homogeneous Dirichlet boundary conditions and initial condition
\begin{equation*}
u(x) = J_0\left(\frac{10 r_0 |x|}{10+\cos 10\theta} \right)
\end{equation*}
on a prescribed domain $\Omega^t$ whose boundary is given by a sinusoidal perturbation of the unit circle. Namely,
\begin{equation*}
\Omega^t = \left\{ x \bigm ||x| < 1 + \frac{1}{10} \cos 10\theta \cos 250 t \right\}
\end{equation*}
with $\theta = \tan^{-1}(x_2/x_1)$.  Fig.~\ref{fig:stefan_snapshots} shows snapshots of the solution, which was computed using quadratic Lagrange elements on a uniform mesh of equilateral triangles ($h = 0.04375$) together with nodal interpolation as the projection operator, relaxation parameters $\delta=0.8$ and $R=3$, and the third-order SDIRK scheme~(\ref{tab:SDIRK3}) with time step $\Delta t = 0.000625$.  The universal mesh and its image under the universal mesh map at an instant in time are shown in Fig.~\ref{fig:stefan_mesh}.

\section{Conclusion} \label{sec:conclusion}

We have presented a general framework for the design of high-order finite element methods for moving boundary problems with prescribed boundary evolution.  A key role in our approach was played by universal meshes, which combine the immunity to large mesh distortions enjoyed by conventional fixed-mesh methods with the geometric fidelity of deforming-mesh methods.  A given accuracy in space and time may be achieved by choosing an appropriate finite element space on the universal mesh and an appropriate time integrator for ordinary differential equations.  The order of accuracy of the resulting scheme is suboptimal by one half an order according to theory, although we observed in our numerical examples that terms of optimal order tend to dominate in practice.

Several aspects of this research motivate further study.  First, we have yet to address problems for which the boundary itself is an unknown, rather than prescribed.  An extension of the method to three dimensions is alluring.  Finally, developing an analogous strategy for domains with lower regularity, such as domains with corners, is an open problem.

\section{Acknowledgments}
This research was supported by the U.S. Department of Energy grant
DE-FG02-97ER25308; Department of the Army Research Grant, grant number:
W911NF-07- 2-0027; and NSF Career Award, grant number: CMMI-0747089.

\bibliographystyle{elsarticle-num}
\bibliography{freeboundary}

\appendix

\section{Singly Diagonally Implicit Runge Kutta time integrators.} \label{sec:appendix_sdirk}

\renewcommand{\arraystretch}{1.3}
\begin{table}[t]
\centering
\caption{Coefficients $\beta_{ij}$ for a $s=2$-stage SDIRK scheme of order 2. ($\gamma = 1-\sqrt{2}/2$)}
\label{tab:SDIRK2}
\begin{tabularx}{120pt}{r|rr}
$i$ $\backslash$ $j$ & 0 & 1 \\
\hline
1 & $1$ & $$ \\
2 & $-\sqrt{2}$ & $1+\sqrt{2}$ \\
\end{tabularx}
\end{table}

\begin{table}[t]
\centering
\caption{Coefficients $\beta_{ij}$ for a $s=3$-stage SDIRK scheme of order 3. ($\gamma=0.43586652150845899942$)}
\label{tab:SDIRK3}
\begin{tabularx}{430pt}{r|rrr}
$i$ $\backslash$ $j$ & 0 & 1 & 2 \\
\hline
1 & $1.00000000000000000$ & $$ & $$ \\
2 & $0.352859819860479140$ & $0.647140180139520860$ & $$ \\
3 & $-1.25097989505606042$ & $3.72932966244456977$ & $-1.47834976738850935$ \\
\end{tabularx}
\end{table}

\begin{table}[t]
\centering
\caption{Coefficients $\beta_{ij}$ for a $s=5$-stage SDIRK scheme of order 4. ($\gamma = 1/4$)}
\label{tab:SDIRK4}
\begin{tabularx}{200pt}{r|rrrrr}
$i$ $\backslash$ $j$ & 0 & 1 & 2 & 3 & 4 \\
\hline
1 & $1$ & $$ & $$ & $$ & $$ \\
2 & $-1$ & $2$ & $$ & $$ & $$ \\
3 & $-\frac{13}{25}$ & $\frac{42}{25}$ & $-\frac{4}{25}$ & $$ & $$ \\
4 & $-\frac{4}{17}$ & $\frac{89}{68}$ & $-\frac{25}{136}$ & $\frac{15}{136}$ & $$ \\
5 & $\frac{7}{3}$ & $-\frac{37}{12}$ & $-\frac{103}{24}$ & $\frac{275}{8}$ & $-\frac{85}{3}$ \\
\end{tabularx}
\end{table}

\renewcommand{\arraystretch}{1}

Tables~\ref{tab:SDIRK2}-\ref{tab:SDIRK4} record the coefficients $\gamma>0$ and $\beta_{ij} \in \mathbb{R}$, $i=1,2,\dots,s$, $j=0,1,\dots,i-1$ for a collection of SDIRK methods~(\ref{dirk}) of orders 2 through 4.

Note that the structure of the Runge-Kutta stages in~(\ref{dirk}) differs from the structure that is most familiar to Runge-Kutta practitioners~\cite{Hairer2002}.  The former structure, which is algorithmically better-suited for problems with time-dependent mass matrices, is obtainable from any L-stable SDIRK scheme as follows.  Let $a_{ij}$, $b_j$, and $c_j$, $i,j=1,2,\dots,s$, be the coefficients of an SDIRK scheme with Butcher tableaux
\begin{equation}\label{Butcher}
\begin{matrix}
c_1 &\vline & a_{11} & a_{12} & \cdots & a_{1s} \\
c_2 &\vline& a_{21} & a_{22} & \cdots & a_{2s} \\
\vdots &\vline & \vdots & \vdots & \ddots & \vdots \\
c_s &\vline & a_{s1} & a_{s2} & \cdots & a_{ss} \\
\hline
 &\vline & b_1 & b_2 & \cdots & b_s
\end{matrix}
\end{equation}
By definition, $a_{11}=a_{22}=\dots=a_{ss}$ and $a_{ij}=0$ for $j>i$.  Assume that the scheme is L-stable, i.e. $b_j = a_{sj}$, $j=1,2,\dots,s$.  Then the coefficients $\gamma$ and $\beta_{ij}$ in the formulation~(\ref{dirk}) are related to $a_{ij}$, $b_j$, and $c_j$ via
\begin{align*}
\gamma &= a_{11} \\
\beta_{ij} &= 
\begin{cases}
\delta_{ij} - a_{ij}^* &\mbox{ if } j>0 \\
\sum_{k=1}^i a_{ik}^* &\mbox{ if } j=0.
\end{cases}
\end{align*}
Here, $\delta_{ij}$ denotes the Kronecker delta and $a_{ij}^*$ is the $i,j$ entry of the matrix $\gamma A^{-1}$, where $A = (a_{ij})$.  The equivalence between~(\ref{dirk}) and the scheme defined by~(\ref{Butcher}) is proven in~\cite{Ying2009}.
\section{The closest point projection onto a moving curve and its time derivative.} \label{sec:appendix_pidot}

The following paragraphs derive explicit expressions for the time derivative of the closest point projection of a fixed point in space onto a moving curve.  Such expressions are needed in numerical implementations for the evaluation of~(\ref{vexact}) when the boundary evolution operator $\gamma^{n,t}_h$ is given by~(\ref{gamma1}).

Consider a moving curve $c^t \in \mathcal{C} := \{s \in C^2([0,1],\mathbb{R}^2) \mid s'(\theta) \neq 0 \; \forall \theta \in [0,1] \}$ whose velocity at any point $y=c^t(\theta) \in \mathrm{image}(c^t)$ is given by $v^t(y) = \dot{c}^t(\theta)$.  Let $\hat{n}^t(y)$, $\hat{t}^{\,t}(y)$, and $\kappa^t(y)$ denote the unit normal vector, unit tangent vector, and signed curvature at $y$, respectively, and let $\pi^t$ and $\phi^t$ denote the closest point projection onto $\mathrm{image}(c^t)$ and the signed distance function on $\mathbb{R}^2$, respectively, as in Section~\ref{sec:universal_meshes}.  Let $\tau$ denote the arclength parameter on $\mathrm{image}(c^t)$.  Henceforth, we employ the arclength parametrization and write $c^t(\tau)$ to denote the point on $\mathrm{image}(c^t)$ with arclength parameter $\tau$.

With respect to the arclength parametrization, the unit normal, unit tangent, and signed curvature satisfy the following relations at any point $y=c^t(\tau)$:
\begin{equation*}
\hat{t}^{\,t}(y) = \frac{\partial c^t}{\partial\tau}(\tau), \;\;\; \frac{\partial\hat{t}^{\,t}}{\partial\tau}(y) = \kappa^t(y)\hat{n}^t(y), \;\;\; \frac{\partial\hat{n}^t}{\partial\tau}(y) = -\kappa^t(y)\hat{t}^{\,t}(y).
\end{equation*}
Here, for a given function $f^t : \mathrm{image}(c^t) \rightarrow \mathbb{R}^k$, $k \in \{1,2\}$, we are abusing notation by writing
\begin{equation*}
\frac{\partial f^t}{\partial\tau}(y) := \left.\frac{\partial}{\partial\tau}\right|_{t} f^t(c^t(\tau))
\end{equation*}
for any $y = c^t(\tau) \in \mathrm{image}(c^t)$.  Likewise, we write
\begin{equation*}
\frac{\partial g^t}{\partial t}(x) = \left.\frac{\partial}{\partial t}\right|_x g^t(x)
\end{equation*}
for a function $g^t : \mathbb{R}^2 \rightarrow \mathbb{R}^k$, $k \in \{1,2\}$.

The closest point projection satisfies
\begin{equation}
x - \pi^t(x) = \phi^t(x) \hat{n}^t(\pi^t(x)) \label{npirelation}
\end{equation}
for any $x \in \mathbb{R}^2$ for which $\pi^t(x)$ is uniquely defined.  Another identity that will be of use momentarily concerns the normal velocity $v^t_n(y) := v^t(y) \cdot \hat{n}^t(y)$.  Namely,
\begin{equation}
\frac{\partial v_n^t}{\partial \tau}(y) 
= \hat{n}^t(y) \cdot \frac{\partial v^t}{\partial \tau}(y) - \kappa^t(y) \hat{t}^{\,t}(y) \cdot v^t(y) \label{dvndtau}
\end{equation}
for any $y \in \mathrm{image}(c^t)$ by the product rule.

\begin{proposition}
Suppose $\{c^t\}_{t \in [0,T]} \subset \mathcal{C}$ is a family of curves such that the map 
\begin{align*}
c : \{(\tau,t) :  0 \le \tau \le \mathrm{length}(\mathrm{image}(c^t)), 0 \le t \le T\} &\rightarrow \mathbb{R}^2 \\
(\tau,t) &\mapsto c^t(\tau)
\end{align*}
is of class $C^2$.  Let $x \in \mathbb{R}^2$ be a point for which $\pi^t(x)$ is uniquely defined and $\phi^t(x)\kappa^t(\pi^t(x)) < 1 $ for every $0 \le t \le T$.  Then 
\begin{equation}
\frac{\partial\pi^t}{\partial t}(x) = v_n^t(\pi^t(x)) \hat{n}^t(\pi^t(x)) + \sigma^t(x) \hat{t}^{\,t}(\pi^t(x))
\end{equation}
for every $0 \le t \le T$, where
\begin{equation}
\sigma^t(x) = \frac{\phi^t(x) \frac{\partial v_n^t}{\partial \tau} (\pi^t(x)) }{1-\phi^t(x)\kappa^t(\pi^t(x))}. \label{sigma}
\end{equation}
\end{proposition}
\proof Let $\bar{\tau}^t(x)$ denote the arclength parameter along $\mathrm{image}(c^t)$ assumed by $\pi^t(x)$; that is,
\begin{equation}
c^t(\bar{\tau}^t(x)) = \pi^t(x). \label{taubar1}
\end{equation}
Differentiating this relation with respect to time gives
\begin{equation}
v^t(\pi^t(x)) + \hat{t}^{\,t}(\pi^t(x)) \frac{\partial\bar{\tau}^t}{\partial t}(x) = \frac{\partial\pi^t}{\partial t}(x). \label{dtaubar1}
\end{equation}
On the other hand, relation~(\ref{npirelation}) implies that
\begin{equation}
\left( x - c^t(\bar{\tau}^t(x)) \right) \cdot \hat{t}^{\,t}(c^t(\bar{\tau}^t(x)) = 0 \label{taubar2}
\end{equation}
for every $t$.  Using the fact that $\hat{t}^{\,t}(c^t(\bar{\tau}^t(x)) = \frac{\partial c^t}{\partial \tau} (\bar{\tau}^t(x))$ has unit length, the time derivative of~(\ref{taubar2}) reads
\begin{equation}
-v^t(\pi^t(x)) \cdot \hat{t}^{\,t}(\pi^t(x)) - \frac{\partial\bar{\tau}^t}{\partial t} + \phi^t(x) \hat{n}^t(\pi^t(x)) \cdot \left( \frac{\partial v^t}{\partial \tau}(\pi^t(x)) + \kappa^t(\pi^t(x)) \hat{n}^t(\pi^t(x)) \frac{\partial\bar{\tau}^t}{\partial t} \right) = 0. \label{dtaubar2}
\end{equation}
Together, relations~(\ref{dtaubar1}) and~(\ref{dtaubar2}) provide enough information to solve for the normal and tangential components of $\frac{\partial\pi^t}{\partial t}(x)$.  

The normal component of $\frac{\partial\pi^t}{\partial t}(x)$ is obtained easily by dotting~(\ref{dtaubar1}) with $\hat{n}^t(\pi^t(x))$, resulting in 
\begin{equation*}
\frac{\partial\pi^t}{\partial t}(x) \cdot \hat{n}^t(\pi^t(x)) = v_n^t(\pi^t(x)).
\end{equation*}

To compute the tangential component $\sigma^t(x) := \frac{\partial\pi^t}{\partial t}s(x) \cdot \hat{t}^{\,t}(\pi^t(x))$, take the dot product of~(\ref{dtaubar1}) with $\hat{t}^{\,t}(\pi^t(x))$ and simplify~(\ref{dtaubar2}) to obtain the following system of equations in two unknowns $\sigma^t(x)$ and $\frac{\partial\bar{\tau}^t}{\partial t}$:
\begin{align*}
v^t(\pi^t(x)) \cdot \hat{t}(\pi^t(x)) + \frac{\partial\bar{\tau}^t}{\partial t}(x) &= \sigma^t(x) \\
-v^t(\pi^t(x)) \cdot \hat{t}(\pi^t(x)) - \frac{\partial\bar{\tau}^t}{\partial t}(x) + \phi^t(x) \hat{n}^t(\pi^t(x)) \cdot \frac{\partial v^t}{\partial \tau}(\pi^t(x)) + \kappa^t(\pi^t(x)) \phi^t(x) \frac{\partial\bar{\tau}^t}{\partial t}(x) &= 0.
\end{align*}
Solving this system and invoking~(\ref{dvndtau}) leads to~(\ref{sigma}).
\qed

\paragraph{Remark}
The restriction $\phi^t(x)\kappa^t(\pi^t(x)) < 1$ in the preceding proposition is mild.  In general, $\phi^t(x) \kappa^t(\pi^t(x)) \le 1$ whenever $\pi^t(x)$ is uniquely defined.  Indeed, since $|x-c^t(\tau)|^2$ is minimal at $\tau=\bar{\tau}^t(x)$, it follows that 
\[
0 \le \left.\frac{\partial^2}{\partial \tau^2} \right|_{\tau=\bar{\tau}^t(x)} |x-c^t(\tau)|^2 = 2(1-\phi^t(x) \kappa^t(\pi^t(x))).
\]  
The assumption of strict inequality rules
out degenerate cases in which $\left.\frac{\partial^2}{\partial
    \tau^2} \right|_{\tau=\bar{\tau}^t(x)} |x-c^t(\tau)|^2 = 0$.

\end{document}